\documentclass[12pt]{amsart}
\usepackage{amssymb,amsmath, amsthm,latexsym}
\usepackage{a4wide}

\theoremstyle{plain}

\newtheorem{lem}{Lemma}[section]
\newtheorem{theo}[lem]{Theorem}
\newtheorem{prop}[lem]{Proposition}

\parskip=\medskipamount
\font\k=cmr7

  \newcommand {\be}{\Re(\rho)}
  \newcommand {\gam}{\Im(\rho)}
  \newcommand {\cu}{\mbox{\k cus}}
  \newcommand {\di}{\mbox{\k dis}}
  \newcommand {\reg}{\mbox{\k reg}}
  \newcommand {\res}{\mbox{\k res}}
  \newcommand {\spec}{\mbox{\k spec}}
  \newcommand {\geo}{\mbox{\k geo}}
  \newcommand {\C}{{\mathbb C}}
  \newcommand {\N}{{\mathbb N}}
  \newcommand {\R}{{\mathbb R}}
  \newcommand {\Z}{{\mathbb Z}}
  \newcommand {\Q}{{\mathbb Q}}
  \newcommand {\A}{{\mathbb A}}
  \newcommand {\af}{{\mathfrak a}}
  \newcommand {\gf}{{\mathfrak g}}
  \newcommand {\mf}{{\mathfrak m}}
  \newcommand {\kf}{{\mathfrak k}}
  \newcommand {\tf}{{\mathfrak t}}
  \newcommand {\gl}{{\mathfrak gl}}
  \newcommand {\hf}{{\mathfrak h}}
  \newcommand {\ho}{{\mathfrak o}}

  \newcommand {\pg}{{\mathfrak p}}
  
  \newcommand {\Pg}{{\mathfrak P}}
 \newcommand {\mX}{{\mathfrak X}}
 \newcommand {\mM}{{\mathfrak M}}
 \newcommand {\mN}{{\mathfrak N}}
 \newcommand {\mO}{{\mathfrak O}}
  
\renewcommand {\H}{{\mathcal H}}
  
  \newcommand {\Co}{{\mathcal C}}
  \newcommand {\cO}{{\mathcal O}}
  \newcommand {\G}{{\mathcal G}}
\renewcommand {\L}{{\mathcal L}}
  \newcommand {\E}{{\mathcal E}}
  
 \newcommand {\cP}{{\mathcal P}}
 \newcommand {\cF}{{\mathcal F}}
 \newcommand {\cL}{{\mathcal L}}
 \newcommand {\cA}{{\mathcal A}}
\newcommand  {\cZ}{{\mathcal Z}}
\newcommand {\ba}{\backslash}
 \newcommand {\ov}{\overline}
  \newcommand {\bs}{{\bf s}}

\renewcommand{\Im}{\operatorname{Im}}
\renewcommand{\Re}{\operatorname{Re}}
\newcommand{\Tr}{\operatorname{Tr}}
\newcommand{\End}{\operatorname{End}}
\newcommand{\sign}{\operatorname{sign}}
\newcommand{\tr}{\operatorname{tr}}
\newcommand{\Id}{\operatorname{Id}}
\newcommand{\Hom}{\operatorname{Hom}}
\newcommand{\Ind}{\operatorname{Ind}}

\newcommand{\vol}{\operatorname{vol}}

\newcommand{\SL}{\operatorname{SL}}
\newcommand{\GL}{\operatorname{GL}}
\newcommand{\SO}{\operatorname{SO}}
\newcommand{\PSL}{\operatorname{PSL}}
\newcommand{\Ad}{\operatorname{Ad}}
\newcommand{\rO}{\operatorname{O}}
\newcommand{\supp}{\operatorname{supp}}

\begin{document}
%Topmatter
\title[Weyl's law  for $\SL_n$]
{\Large\bf Weyl's law for the cuspidal spectrum of  $\SL_n$}
\date{\today}

\author{Werner M\"uller}
\address{Universit\"at Bonn\\
Mathematisches Institut\\
Beringstrasse 1\\
D -- 53115 Bonn, Germany}
\email{mueller@math.uni-bonn.de}


\begin{abstract}
Let $\Gamma$ be a principal congruence subgroup of $\SL_n(\Z)$ and let 
$\sigma$ be an irreducible unitary representation of $\SO(n)$. Let
$N^\Gamma_{\cu}(\lambda,\sigma)$ be the counting function of the eigenvalues
of the Casimir operator acting in the space of cusp forms for $\Gamma$ which
transform under $\SO(n)$ according to $\sigma$. In this paper we prove that
the counting function $N^\Gamma_{\cu}(\lambda,\sigma)$ satisfies Weyl's law.
Especially, this implies that there exist infinitely many cusp forms for the
full modular group $\SL_n(\Z)$.
\end{abstract}

%End topmatter
\maketitle

\setcounter{section}{-1}
\section{Introduction}

Let $G$ be a connected reductive algebraic group over $\Q$ and let
$\Gamma\subset G(\Q)$ be an arithmetic subgroup. An important problem in the
theory of automorphic forms is the question of the existence and the 
construction of cusp forms for $\Gamma$. By Langlands' theory of Eisenstein
 series \cite{La}, cusp forms are the building blocks of the spectral 
resolution of the regular representation of $G(\R)$ in 
$L^2(\Gamma\ba G(\R))$. Cusp forms are also fundamental in number theory. 
Despite of their importance, very little is known about the existence of 
cusp forms in general.
In this paper we will address the question of existence of cusp forms 
for the group $G=\SL_n$. The main purpose of
this paper is to prove that cusp forms exist in abundance for congruence
subgroups of $\SL_n(\Z)$, $n\ge2$. 

To formulate our main result  we need to introduce some notation. For 
simplicity assume that $G$ is semisimple. 
Let $K_\infty$ be a maximal compact subgroup of $G(\R)$ and let 
$X=G(\R)/K_\infty$ be the associated Riemannian symmetric space. Let
$\cZ(\gf_\C)$
be the center of the unviersal enveloping algebra of the complexification of
the Lie algebra $\gf$ of $G(\R)$. Recall  that a cusp form  for 
$\Gamma$ in the sense of \cite{La} is a
smooth and $K_\infty$-finite function $\phi:\Gamma\ba G(\R)\to\C$ which is a
simultaneous eigenfunction of $\cZ(\gf_\C)$ and which satisfies 
$$\int_{\Gamma\cap N_P(\R)\ba N_P(\R)}\phi(nx)\;dn=0,$$
for all unipotent radicals $N_P$ of proper rational parabolic subgroups $P$
of $G$. We note that each cusp form $f\in C^\infty(\Gamma\ba G(\R))$
is rapidly decreasing on $\Gamma\ba G(\R)$ and hence square integrable. Let
$L^2_{\cu}(\Gamma\ba G(\R))$ be the closure of the linear span of all cusp
forms.
 Let $(\sigma,V_\sigma)$ be an irreducible unitary representation
of $K_\infty$. Set
$$L^2(\Gamma\ba G(\R),\sigma)=(L^2(\Gamma\ba G(\R))\otimes
V_\sigma)^{K_\infty}$$
and define $L^2_{\cu}(\Gamma\ba G(\R),\sigma)$ similarly.  Then 
$L^2_{\cu}(\Gamma\ba G(\R),\sigma)$ is the space of cusp forms with fixed
$K_ \infty$-type $\sigma$. Let $\Omega_{G(\R)}\in\cZ(\gf_\C)$ be the Casimir
 element of $G(\R)$. Then $-\Omega_{G(\R)}\otimes\Id$ induces a selfadjoint
 operator
$\Delta_\sigma$  in the Hilbert space $L^2(\Gamma\ba G(\R),\sigma)$ 
which is bounded from below. 
If $\Gamma$ is torsion free, 
$L^2(\Gamma\ba G(\R),\sigma)$ is isomorphic to
the space $L^2(\Gamma\ba X,E_\sigma)$ of square integrable sections of the
locally homogeneous vector bundle $E_\sigma$ associated to $\sigma$,
and $\Delta_\sigma=(\nabla^\sigma)^*\nabla^\sigma
-\lambda_\sigma\Id$, where $\nabla^\sigma$ is the canonical invariant
 connection and 
$\lambda_\sigma$ the Casimir eigenvalue of $\sigma$.
This shows that  $\Delta_\sigma$ is a second order elliptic differential
 operator. Especially, if $\sigma_0$ is the trivial representation,
then $L^2(\Gamma\ba G(\R),\sigma_0)\cong L^2(\Gamma\ba X)$ and 
$\Delta_{\sigma_0}$ equals the Laplacian $\Delta$ of $X$.  

The restriction of $\Delta_\sigma$ to the
subspace $L^2_{\cu}(\Gamma\ba G(\R),\sigma)$ has pure point spectrum
consisting of eigenvalues 
$\lambda_0(\sigma)<\lambda_1(\sigma)<\cdots$
of finite multiplicity. We call it
 the {\it cuspidal spectrum} of $\Delta_\sigma$. 
A convenient way of counting the number of cusp forms
for $\Gamma$ is to use their Casimir eigenvalues. 
For this purpose we introduce the counting function
$N^{\Gamma}_{\cu}(\lambda,\sigma)$, $\lambda\ge0$,  for the cuspidal spectrum
 of type $\sigma$
which  is defined as follows. Let $\E(\lambda_i(\sigma))$ be the eigenspace 
corresponding to the eigenvalue $\lambda_i(\sigma)$. Then
$$N^\Gamma_{\cu}(\lambda,\sigma)=\sum_{\lambda_i(\sigma)\le\lambda}
\dim\E(\lambda_i(\sigma)).$$
For non-uniform lattices $\Gamma$ the selfadjoint operator $\Delta_\sigma$
has a large continuous spectrum so that almost all of the eigenvalues of
$\Delta_\sigma$ will be embedded in the continous spectrum. This makes it very
difficult to study the cuspidal spectrum of $\Delta_\sigma$. 

The first results concerning the growth of the cuspidal spectrum are due to
 Selberg \cite{Se}. Let $H$ be the upper half-plane and let $\Delta$ be the 
hyperbolic Laplacian of $H$. Let $N_{\cu}^{\Gamma}(\lambda)$ be the counting
function of the cuspidal spectrum of $\Delta$. In this case the cuspidal 
eigenfunctions of
$\Delta$ are  called {\it Maass cusp forms}.
Using the trace formula, Selberg \cite[p.668]{Se} proved that for 
every congruence
subgroup $\Gamma\subset\SL_2(\Z)$, the counting function satisfies Weyl's law, 
i.e.
\begin{equation}\label{0.1}
N_{\cu}^{\Gamma}(\lambda)\sim\frac{\vol(\Gamma\ba H)}{4\pi}\lambda
\end{equation}
as $\lambda\to\infty$. In particular this implies that for congruence
 subgroups  of $\SL_2(\Z)$  there exist
as many Maass cusp forms as one can expect. On the other hand, it is 
conjectured
by Phillips and Sarnak \cite{PS} that for a non-uniform lattice
$\Gamma$ of $\SL_2(\R)$ whose Teichm\"uller space $T$ is non trivial and
different from the Teichm\"uller space corresponding to the once
punctured torus, a generic lattice $\Gamma\in T$ has only finitely
many Maass cusp forms. This indicates that the existence of cusp forms is very
subtle and may be related to the arithmetic nature of $\Gamma$. 

Let $d=\dim X$. 
In has been conjectured in \cite{Sa}  that for rank$(X)>1$
and $\Gamma$ an irreducible lattice
\begin{equation}\label{0.2}
\limsup_{\lambda\to\infty}\frac{N_{\cu}^{\Gamma}(\lambda)}{\lambda^{d/2}}=
\frac{\vol(\Gamma\ba X)}{(4\pi)^{d/2}\Gamma(d/2+1)},
\end{equation}
where $\Gamma(s)$ denotes the Gamma function. A lattice $\Gamma$ for which
(\ref{0.2}) holds is called by Sarnak {\it essentially cuspidal}. An analogous
conjecture was made in \cite[p.180]{Mu3} for the counting function
$N^\Gamma_{\di}(\lambda,\sigma)$  of the discrete spectrum of any Casimir 
operator $\Delta_\sigma$. This conjecture
states that for any arithmetic subroup $\Gamma$ and any $K_\infty$-type 
$\sigma$ 
\begin{equation}\label{0.2a}
\limsup_{\lambda\to\infty}\frac{N_{\di}^{\Gamma}(\lambda,\sigma)}
{\lambda^{d/2}}=
\dim(\sigma)\frac{\vol(\Gamma\ba X)}{(4\pi)^{d/2}\Gamma(d/2+1)}.
\end{equation}
Up to now these conjectures have been verified only in a few cases. 
Besides of Selberg's result, Weyl's law (\ref{0.2}) has been proved in the
following cases: 
For congruence subgroups of $G=\SO(n,1)$ by Reznikov \cite{Rez}, for 
 congruence subgroups of $G=R_{F/\Q}\SL_2$, where $F$ is a totally real
number field, by Efrat \cite[p.6]{Ef}, and  for $\SL_3(\Z)$ by St. Miller
 \cite{Mil}.

In this paper we will prove that each principal congruence 
subgroup $\Gamma$ of
$\SL_n(\Z)$, $n\ge 2$, is essentially cuspidal, i.e. Weyl's law holds
for $\Gamma$. Actually we prove the corresponding result for all 
$K_\infty$-types $\sigma$. Our main result is the following theorem.
\begin{theo}\label{th0.1}
For $n\ge2$ let $X_n=\SL_n(\R)/\SO(n)$. Let $d_n=\dim X_n$.
For every principal congruence subgroup  $\Gamma$ 
 of $\SL_n(\Z)$ and  every irreducible unitary representation $\sigma$ of
$\SO(n)$ we have
\begin{equation}\label{0.3}
N_{\cu}^{\Gamma}(\lambda,\sigma)\sim
\dim(\sigma)\frac{\vol(\Gamma\ba X_n)}
{(4\pi)^{d_n/2}\Gamma(d_n/2+1)}\lambda^{d_n/2}
\end{equation}
as $\lambda\to\infty$.
\end{theo}

If $n=2$ the asymptotic formula (\ref{0.3}) can be improved with a 
remainder term  of order $O(\sqrt{\lambda}/\log\lambda)$ \cite{Se}. It is
an interesting problem to obtain estimations for the remainder term for
all $n\ge2$. 

Next we reformulate Theorem \ref{0.1} in the ad\`elic language. 
Let $G=\GL_n$ regarded as an algebraic group over $\Q$. 
 Let $\A$ be the ring of ad\`eles of $\Q$.  Denote by $A_G$ the split
component of the center of $G$ and let $A_G(\R)^0$ be the component of 1 in
$A_G(\R)$. Let $\xi_0$ be the trivial character of 
$A_G(\R)^0$ and denote by $\Pi(G(\A),\xi_0)$ the set of equivalence classes of 
irreducible unitary representations of $G(\A)$ whose central character is
trivial on $A_G(\R)^0$.
Let $L^2_{\cu}(G(\Q)A_G(\R)^0\ba G(\A))$ be the subspace of cusp
forms in  $L^2(G(\Q)A_G(\R)^0\ba G(\A))$.  
Denote by $\Pi_{\cu}(G(\A),\xi_0)$ the subspace of all
$\pi$ in $\Pi(G(\A),\xi_0)$ which are equivalent to a subrepresentation of
the regular representation in $L^2_{\cu}(G(\Q)A_G(\R)^0\ba G(\A))$.
 By \cite{Sk} the multiplicity of any
$\pi\in\Pi_{\cu}(G(\A),\xi_0)$ in the space of cusp forms 
$L^2_{\cu}(G(\Q)A_G(\R)^0\ba G(\A))$ is one. 
 Let $A_f$ be the ring of finite ad\`eles. Any irreducible unitary
representation $\pi$ of $G(\A)$ can be written as
$\pi=\pi_\infty\otimes\pi_f,$
where $\pi_\infty$ and $\pi_f$ are irreducible unitary representations
of $G(\R)$ and $G(\A_f)$, respectively. Let $\H_{\pi_\infty}$ and $\H_{\pi_f}$
denote the Hilbert space of the representation $\pi_\infty$ and $\pi_f$,
respectively. Let $K_f$ be an open compact subgroup of $G(\A_f)$. Denote
by $\H_{\pi_f}^{K_f}$ the subspace of $K_f$-invariant vectors in $\H_{\pi_f}$.
Given an irreducible unitary representation $\sigma$ of $\rO(n)$,
 denote by $\H_{\pi_\infty}(\sigma)$ the
$\sigma$-isotypical subspace of $\H_{\pi_\infty}$. 
Let $G(\R)^1$ be the subgroup of all $g\in G(\R)$ with $|\det(g)|=1$. 
Given $\pi\in\Pi(G(\A),\xi_0)$, denote by
$\lambda_\pi$ the Casimir eigenvalue of the restriction of $\pi_\infty$
to $G(\R)^1$.  For
$\lambda\ge0$ let $\Pi_{\cu}(G(\A),\xi_0)_\lambda$ be the space of all
$\pi\in\Pi_{\cu}(G(\A),\xi_0)$ which satisfy $|\lambda_\pi|\le\lambda$.
Set $\varepsilon_{K_f}=1$, if $-1\in K_f$ and $\varepsilon_{K_f}=0$
otherwise. Then we have
\begin{theo}\label{th0.2}
Let $G=\GL_n$ and let $d_n=\dim\SL_n(\R)/\SO(n)$. Let $K_f$ be an open 
compact subgroup of $G(\A_f)$ and let $\sigma$ be an irreducible unitary 
representation of $\SO(n)$. Then
\begin{equation}\label{0.4}
\begin{split}
\sum_{\pi\in\Pi_{\cu}(G(\A),\xi_0)_\lambda}&\dim\bigl(\H_{\pi_f}^{K_f}\bigr)
\dim\bigl(\H_{\pi_\infty}(\sigma)\bigr)\\
&\sim\dim(\sigma)
\frac{\vol(G(\Q)A_G(\R)^0\ba G(\A)/K_f)}{(4\pi)^{d_n/2}\Gamma(d_n/2+1)}
(1+\varepsilon_{K_f})\lambda^{d_n/2}
\end{split}
\end{equation}
as $\lambda\to\infty$.
\end{theo}
This may be regarded as the ad\`elic version of Weyl's law for $\GL_n$.
A similar result holds if we replace $\xi_0$ by any unitary character of
$A_G(\R)^0$. If we specialize Theorem \ref{0.2} to the congruence subgroup
$K(N)$ which defines $\Gamma(N)$, we obtain Theorem \ref{0.1}.

Theorem \ref{th0.2} will be derived from the Arthur trace formula combined with
the heat equation method. The heat equation method is a very convenient way to
derive Weyl's law for the  counting function of the eigenvalues of the
 Laplacian on a compact Riemannian manifold \cite{Cha}. It is based on the 
study of the asymptotic behaviour of the trace of the heat operator. Our
approach is similar. We will
use the Arthur trace formula to compute the trace of the heat operator on the 
discrete spectrum and to determine its asymptotic behaviour as $t\to0$.

We will now describe our method in more detail.
Let $G(\A)^1$ be the subgroup of all $g\in G(\A)$ satisfying $|\det(g)|=1$. 
Then $G(\Q)$ is contained in $G(\A)^1$ and
the noninvariant trace formula of Arthur \cite{A1} is an identity
\begin{equation}\label{0.6}
\sum_{\chi\in\mX}J_\chi(f)=\sum_{\ho\in\mO}J_\ho(f),
\quad f\in C_c^\infty(G(\A)^1),
\end{equation}
between distributions on $G(\A)^1$. 
The left hand side is  the {\it spectral side} $J_{\spec}(f)$ 
and the right hand side the {\it geometric side} $J_{\geo}(f)$ of the 
trace formula. The distributions $J_\chi$
are defined in terms of truncated Eisenstein series. They are parametrized 
by the set  of cuspidal data $\mX$. 
The distributions $J_{\ho}$ are parametrized by semisimple 
conjugacy in $G(\Q)$ and are closely related to weighted orbital integrals on
$G(\A)^1$. 

For simplicity we consider only the case of the trivial $K_\infty$-type.
We choose a certain family of test functions $\widetilde\phi_t^1
\in C^\infty_c(G(\A)^1)$, depending on $t>0$, which at the infinite place 
are given by the heat kernel $h_t\in C^\infty(G(\R)^1)$ of the Laplacian 
on $X$, multiplied by a certain cutoff function $\varphi_t$, and which
at the finite places is given by the normalized characteristic function of
an open compact subgroup $K_f$ of $G(\A_f)$. Then we evaluate the spectral and
the geometric side at $\widetilde\phi_t^1$ and study their asymptotic 
bahaviour as $t\to0$. Let $\Pi_{\di}(G(\A),\xi_0)$ be the set of irreducible
 unitary representations of
$G(\A)$ which occur discretely in the regular representation of $G(\A)$ in
$L^2(G(\Q)A_G(\R)^0\ba G(\A))$. Given $\pi\in\Pi_{\di}(G(\A),\xi_0)$, let 
$m(\pi)$ denote the multiplicity with which $\pi$ occurs in 
$L^2(G(\Q)A_G(\R)^0\ba G(\A))$. Let $\H_{\pi_\infty}^{K_\infty}$ be the space 
of $K_\infty$-invariant vectors in $\H_{\pi_\infty}$.
 Comparing the asymptotic behaviour of the two
 sides of the trace formula, we obtain 
\begin{equation}\label{0.7}
\begin{split}
\sum_{\pi\in\Pi_{\di}(G(\A),\xi_0)}m(\pi)
&e^{t\lambda_{\pi}}\dim(\H_{\pi_f}^{K_f})\dim(\H_{\pi_\infty}^{K_\infty})\\
&\sim\frac{\vol(G(\Q)\ba G(\A)^1/K_f)}{(4\pi)^{d_n/2}}
(1+\varepsilon_{K_f})t^{-d_n/2}
\end{split}
\end{equation}
as $t\to0$, where the notation is as in Theorem \ref{th0.2}. Applying
Karamatas theorem \cite[p.446]{Fe}, we  obtain Weyl's law
for the discrete spectrum with respect to the trivial $K_\infty$-type. 
A nontrivial $K_\infty$-type can be treated in the same way. 
The discrete spectrum is the union of the
cuspidal and the residual spectrum. 
It follows from  \cite{MW} combined with Donnelly's estimation of the
cuspidal spectrum \cite{Do}, that the order of growth of the counting function
 of the residual spectrum for $\GL_n$ is at most  $O(\lambda^{(d_n-1)/2})$ as
$\lambda\to\infty$. This implies (\ref{0.4}). 
 
To study the asymptotic behaviour of the geometric side, we use the fine
 ${\ho}$-expansion
\cite{A10} 
\begin{equation}\label{0.8}
J_{\geo}(f)=\sum_{M\in\cL}\sum_{\gamma\in(M(\Q_S))_{M,S}}a^M(S,\gamma)
J_M(f,\gamma),
\end{equation}
which expresses the distribution $J_{\geo}(f)$ in terms of weighted  orbital 
integrals $J_M(\gamma,f)$.
Here $M$ runs over the set of Levi subgroups $\cL$ containing the Levi
component $M_0$ of the standard minimal parabolic subgroup $P_0$, $S$ is
a finite set of places of $\Q$, and $(M(\Q_S))_{M,S}$ is a certain set of
equivalence classes in $M(\Q_S)$. This reduces our problem to the 
investigation of weighted orbital integrals. The key result is that 
$$\lim_{t\to0}t^{d_n/2}J_M(\widetilde\phi_t^1,\gamma)=0,$$
unless  $M=G$ and $\gamma =\pm1$. 
The contributions to (\ref{0.8}) of the terms where $M=G$ and
 $\gamma=\pm1$ are easy to determine. Using
the behaviour of the heat kernel $h_t(\pm1)$ as $t\to0$, it follows that
\begin{equation}\label{0.9}
J_{\geo}(\widetilde\phi_t^1)\sim\frac{\vol(G(\Q)\ba G(\A)^1/K_f)}{(4\pi)^{d/2}}
(1+\varepsilon_{K_f})t^{-d_n/2}
\end{equation}
as $t\to0$. 

To deal with the spectral side, we 
use the results of \cite{MS}. Let $\Co^1(G(\A)^1)$ denote the space of 
integrable rapidly decreasing functions on $G(\A)^1$ (see \cite[\S 1.3]{Mu2}
for its definition). By Theorem 0.1 of \cite{MS}, the spectral side is 
absolutely
convergent for all $f\in\Co^1(G(\A)^1)$. Furthermore, it can be written as
a finite linear combination
\begin{equation*}
J_{\spec}(f)=\sum_{M\in\cL}\sum_{L\in\cL(M)}\sum_{P\in\cP(M)}
\sum_{s\in W^L(\af_M)_{\reg}}
 a_{M,s} J^L_{M,P}(f,s).
\end{equation*}
of distributions $J^L_{M,P}(f,s)$, where $\cL(M)$ is the set of Levi subgroups
containing $M$, $\cP(M)$ denotes the set of parabolic subgroups with Levi
component $M$ and  $W^L(\af_M)_{\reg}$ is a certain set of Weyl group elements.
Given $M\in\cL$,  the main ingredients of the distribution $J^L_{M,P}(f,s)$
are generalized logarithmic derivatives of the intertwining operators
$$M_{Q|P}(\lambda):\cA^2(P)\to\cA^2(Q),\quad P,Q\in\cP(M),\;
\lambda\in\af_{M,\C}^*,$$
acting between the spaces of automorphic forms attached to $P$ and $Q$,
 respectively. First of all, Theorem 0.1 of  \cite{MS} allows us to  replace
 $\widetilde\phi_t^1$ by a similar function $\phi_t^1\in\Co^1(G(\A)^1)$ which
 is given
as the product of the heat kernel at the infinite place and the normalized
characteristic function of $K_f$. 
Consider the distribution where $M=L=G$. Then $s=1$ and
\begin{equation}\label{0.10}
J^G_{G,G}(\phi_t^1)=\sum_{\pi\in\Pi_{\di}(G(\A),\xi_0)}m(\pi)
e^{t\lambda_{\pi}}\dim(\H_{\pi_f}^{K_f})\dim(\H_{\pi_\infty}^{K_\infty}).
\end{equation}
This is exactly the left hand side of (\ref{0.7}). Thus in order to prove
(\ref{0.7}) we need to show that for all proper Levi subgroups $M$, all $L\in
\cL(M)$, $P\in\cP(M)$ and $s\in W^L(\af_M)_{\reg}$, we have
\begin{equation}\label{0.11}
J^L_{M,P}(\phi_t^1,s)=O(t^{-(d_n-1)/2})
\end{equation} 
as $t\to0$. This is the key result where we really need that our group is
$\GL_n$. It relies on estimations of the logarithmic derivatives of
intertwining operators for $\lambda\in i\af_M^*$. Given 
$\pi\in\Pi_{\di}(M(\A),\xi_0)$, let 
$M_{Q|P}(\pi,\lambda)$ be the restriction
of the intertwining operator $M_{Q|P}(\lambda)$ to the subspace 
$\cA^2_\pi(P)$ of automorphic forms of type $\pi$. The intertwining 
operators can be normalized by certain meromorphic functions 
$r_{Q|P}(\pi,\lambda)$ \cite{A7}. Thus
$$M_{Q|P}(\pi,\lambda)=r_{Q|P}(\pi,\lambda)^{-1}N_{Q|P}(\pi,\lambda),$$
where $N_{Q|P}(\pi,\lambda)$ are the normalized intertwining operators. Using 
Arthur's theory of $(G,M)$-families \cite{A5}, our problem can be reduced to 
the estimation of derivatives of $N_{Q|P}(\pi,\lambda)$ and 
$r_{Q|P}(\pi,\lambda)$ on $i\af_M^*$. The derivatives of $N_{Q|P}(\pi,\lambda)$
can be estimated using Proposition 0.2 of \cite{MS}. Let 
$M=\GL_{n_1}\times\cdots\times\GL_{n_r}$. Then $\pi=\otimes_i\pi_i$ with 
$\pi_i\in\Pi_{\di}(\GL_{n_i}(\A)^1)$ and the normalizing factors 
$r_{Q|P}(\pi,\lambda)$ are given in terms of the Rankin-Selberg $L$-functions
 $L(s,\pi_i\times\widetilde\pi_j)$ and the corresponding $\epsilon$-factors
$\epsilon(s,\pi_i\times\widetilde\pi_j)$. So our problem is finally reduced to 
the estimation of the logarithmic derivative of Rankin-Selberg $L$-functions
on the line $\Re(s)=1$. Using the available knowledge of the analytic
properties of Rankin-Selberg $L$-functions together with standard methods of
analytic number theory, the necessary estimates can be derived. 

In the proof of Theorem \ref{th0.1} and Theorem \ref{th0.2} we have used the 
following key results which at present are only known for $\GL_n$:
1) The nontrivial bounds of the Langlands parameters of
local components of cuspidal automorphic representations 
\cite{LRS} which are needed in \cite{MS}, 2) The description of the residual
 spectrum given in \cite{MW},
3) The theory of the Rankin-Selberg $L$-functions \cite{JPS}.

The paper is organized as follows. In section 2 we prove some estimations 
for the heat kernel on a symmetric space. In section 3 we establish some
estimates for the growth of the discrete spectrum in general. We are
essentially using Donnelly's result \cite{Do} combined with the description
of the residual spectrum \cite{MW}. The main purpose of section 4 is to
prove estimates for the growth  of the number of poles of  Rankin-Selberg
$L$-functions in the critical strip. We use these results in section 5
to establish the key estimates for the logarithmic derivatives of 
normalizing factors. In section 6 we study the asymptotic behaviour of the
spectral side $J_{\spec}(\phi_t^1)$. Finally, in section 7 we study the
asymptotic behaviour of the geometric side, compare it to the asymptotic 
behaviour of the spectral  side and prove the main results.

{\bf Acknowledgment.} The author would like to thank W. Hoffmann, D.
Ramakrishnan and P. Sarnak for very helpful discussions on parts of this
paper. Especially Lemma \ref{l7.1} is due to W. Hoffmann.

\section{Preliminaries}
\setcounter{equation}{0}

\subsection{}

Fix a positive integer $n$ and let $G$ be the group $\GL_n$ considered as
algebraic group over $\Q$. By a parabolic subgroup of $G$ we will
always mean a parabolic subgroup which is defined over $\Q$. Let
$P_0$ be the subgroup of upper triangular matrices of $G$. The Levi subgroup 
$M_0$ of $P_0$ is the group of diagonal matrices in $G$. A
parabolic subgroup $P$ of $G$ is called standard, if $P\supset
P_0$. By a Levi subgroup we will mean a subgroup of $G$ which contains $M_0$
and is the Levi component of a parabolic subgroup of $G$ defined over $\Q$.
 If $M\subset L$ are Levi
subgroups, we denote the set of Levi subgroups of $L$ which contain $M$ by
${\cL}^L(M)$. Furthermore, let ${\cF}^L(M)$ denote the set of
parabolic subgroups of $L$ defined over $\Q$ which contain $M$, and let
${\cP}^L(M)$ be the set of groups in ${\cF}^L(M)$ for which $M$
is a Levi component. If $L=G$, we shall denote these sets by
${\cL}(M)$, ${\cF}(M)$ and ${\cP}(M)$. Write $\cL=\cL(M_0)$. Suppose that  
$P\in{\cF}^L(M)$. Then 
$$P=N_PM_P,$$
where $N_P$ is the unipotent radical of $P$ and $M_P$ is the unique Levi
component of $P$ which contains $M$.

Let $M\in\cL$ and denote by $A_M$  the split component of the center of $M$. 
Then $A_M$ is defined over $\Q$. Let $X(M)_{\Q}$ be the group of characters of
$M$ defined over $\Q$ and set 
\begin{equation*}
{\mathfrak a}_M=\mbox{Hom}(X(M)_{\Q},{\mathbb R}).
\end{equation*}
Then $\af_M$ is a real vector space whose dimension equals that of $A_M$. Its
dual space is
\begin{equation*}
{\mathfrak a}^*_M= X(M)_{\Q}\otimes \R.
\end{equation*}

Let $P$ and $Q$ be groups in $\cF(M_0)$ with $P\subset Q$. Then there is a
canonical surjection $\af_P\rightarrow \af_Q$ and a canonical injection ${\mathfrak a}^*_Q\hookrightarrow {\mathfrak a}^*_P$.
The kernel of the first map will be denoted by
${\mathfrak a}^Q_P.$ Then the dual vector space of
${\mathfrak a}^Q_P$ is ${\mathfrak a}^*_P /{\mathfrak a}^*_Q$.

Let $P\in{\cF}(M_0)$. We shall denote the roots of $(P,A_P)$ by $\Sigma_P$, 
 and the simple roots by $\Delta_P$. Note that for $\GL_n$ all roots are
reduced. They are elements in $X(A_P)_\Q$ and are canonically embedded in 
$\af_P^*$.

For any $M\in\L$ there exists a partition $(n_1,...,n_r)$ of $n$
such that 
$$M=\GL_{n_1}\times\cdots\times\GL_{n_r}.$$ 
Then $\af_M^*$
can be canonically identified with $(\R^r)^*$ and the Weyl group $W(\af_M)$
coincides with the group $S_r$ of permutations of the set
$\{1,...,r\}$.

\subsection{}

Let $F$ be a local field of characteristic zero. If $\pi$ is an admissible
representation of $\GL_m(F)$, we shall denote by $\widetilde\pi$ the
contragredient representation to $\pi$. 
Let $\pi_i$, $i=1,...,r$, be irreducible admissible
representations of the group $\GL_{n_i}(F)$. Then
$\pi=\pi_1\otimes\cdots\otimes\pi_r$ is an irreducible admissible
representation of
$$M(F)=\GL_{n_1}(F)\times\cdots\times\GL_{n_r}(F).$$ For
$\bs\in\C^r$ let $\pi_i[s_i]$ be the representation of
$\GL_{n_i}(F)$ which is defined by
$$\pi_i[s_i](g)=|\det(g)|^{s_i}\pi_i(g),\quad g\in\GL_{n_i}(F).$$
Let 
$$I_P^G(\pi,\bs)=
\Ind^{G(F)}_{P(F)}(\pi_1[s_1]\otimes\cdots\otimes\pi_r[s_r])$$ 
be the
induced representation and denote by $\H_P(\pi)$ the Hilbert space
of the representation $I_P^G(\pi,\bs)$. We refer to $\bs$ as the
continuous parameter of $I_P^G(\pi,\bs)$. Sometimes we will write 
$I^G_P(\pi_1[s_1],...,\pi_r[s_r])$ in place of $I^G_P(\pi,\bs)$.

\subsection{}

Let $\G$ be a locally compact topological group. Then we denote by $\Pi(\G)$ 
the set of equivalence classes of irreducible unitary representations of $\G$.

\subsection{}

Let  $M\in\cL$. Denote by $A_M(\R)^0$ the
component of 1 of $A_M(\R)$. Set
$$M_P(\A)^1=\bigcap_{\chi\in X(M)_\Q}\ker(|\chi|).$$
This is a closed subgroup of $M(\A)$, and $M(\A)$ is the direct product of
$M(\A)^1$ and $A_M(\R)^0$.

Given a unitary character $\xi$  of $A_M(\R)^0$,
denote by $L^2(M(\Q)\ba M(\A),\xi)$ the space of all measurable 
functions
$\phi$ on $M(\Q)\ba M(\A)$ such that 
$$\phi(xm)=\xi(x)\phi(m),\quad 
x\in A_M(\R)^0,\; m\in M(\A),$$
 and $\phi$ is square integrable on 
$M(\Q)\ba M(\A)^1$. Let $L^2_{\di}(M(\Q)\ba M(\A),\xi)$ denote the discrete
subspace of $L^2(M(\Q)\ba M(\A),\xi)$ and let
$L^2_{\cu}(M(\Q)\ba M(\A),\xi)$ be the subspace
of cusp forms in $L^2(M(\Q)\ba M(\A),\xi)$. The orthogonal complement of
$L^2_{\cu}(M(\Q)\ba M(\A),\xi)$ in
the discrete subspace is the resiudal subspace 
$L^2_{\res}(M(\Q)\ba M(\A),\xi)$.
Denote by  $\Pi_{\di}(M(\A),\xi)$, $\Pi_{\cu}(M(\A),\xi)$, 
and $\Pi_{\res}(M(\A),\xi)$ 
the subspace of all
$\pi\in\Pi(M(\A),\xi)$ which are equivalent to a subrepresentation of the 
regular representation of $M(\A)$ in $L^2(M(\Q)\ba M(\A),\xi)$, 
$L^2_{\cu}(M(\Q)\ba M(\A),\xi)$, and $L^2_{\res}(M(\Q)\ba M(\A),\xi)$,
respectively.

Let $\Pi_{\di}(M(\A)^1)$ be the subspace of all $\pi\in\Pi(M(\A)^1)$ which
are equivalent to a subrepresentation of the regular representation of 
$M(\A)^1$ in $L^2(M(\Q)\ba M(\A)^1)$. We denote by $\Pi_{\cu}(M(\A)^1)$ 
(resp. $\Pi_{\res}(M(\A)^1)$) the subspaces of all $\pi\in\Pi_{\di}(M(\A)^1)$ 
occurring in the cuspidal (resp. residual) subspace 
$L^2_{\cu}(M(\Q)\ba M(\A)^1)$ (resp. $L^2_{\res}(M(\Q)\ba M(\A)^1)$).

\subsection{} 

Let $P$ be a parabolic subgroup of $G$. We denote by $\cA^2(P)$ the space of
square integrable automorphic forms on 
$N_P(\A) M_P(\Q)A_P(\R)^0\ba G(\A)$ (see \cite[\S 1.7]{Mu2}).

Given $\pi\in\Pi_{\di}(M_P(\A),\xi_0)$, 
let $\cA^2_\pi(P)$ be the subspace of 
$\cA^2(P)$ of automorphic forms of type $\pi$ \cite[p.925]{A1}. Let $\pi\in
\Pi(M_P(\A)^1)$. We identify $\pi$ with a representation of $M_P(\A)$ which is
trivial on $A_P(\R)^0$. Hence we can define $\cA^2_\pi(P)$ for any 
$\pi\in\Pi(M_P(\A)^1)$.  It is a space of square integrable functions on
$N_P(\A)M_P(\Q)A_P(\R)^0\ba G(\A)$ such that for every $x\in G(\A)$, 
the function
$$\phi_x(m)=\phi(mx),\quad m\in M_P(\A),$$
belongs to the $\pi$-isotypical subspace of the regular representation of
$M_P(\A)$ in the Hilbert space $L^2(A_P(\R)^0M_P(\Q)\ba M_P(\A))$.

\section[Heat kernels]{Heat kernel estimates}
\setcounter{equation}{0}

In this section we shall prove some estimates for the heat kernel of the
Bochner-Laplace operator acting on sections of a homogeneous vector bundle
over a symmetric space.
Let $G$ be a connected semisimple algebraic group defined over $\Q$. Let
$K_\infty$ be a maximal compact subgroup of $G(\R)$ and
let $(\sigma,V_\sigma)$ be an irreducible unitary representation of
$K_\infty$ on a complex vector space $V_\sigma$. Let
$\widetilde E_\sigma=(G(\R)\times V_\sigma)/K_\infty$
be the associated homogeneous vector bundle over $X=G(\R)/K_\infty$. We equip
$\widetilde E_\sigma$ with the $G(\R)$--invariant Hermitian fibre metric which
 is
induced by the inner product in $V_\sigma.$ Let $C^\infty(\widetilde
E_\sigma),\;C^\infty_c(\widetilde E_\sigma)$ and
$L^2(\widetilde E_\sigma)$ denote the space of smooth sections, the space of
compactly supported smooth sections and the Hilbert space of square integrable
sections of $\widetilde E_\sigma,$ respectively. Then we have
\begin{equation} \label{1.1}
C^\infty (\widetilde E_\sigma)   
= (C^\infty(G(\R))\otimes V_\sigma )^{K_\infty},\quad
L^2(\widetilde E_\sigma)   = (L^2(G(\R))\otimes V_\sigma )^{K_\infty}
\end{equation}
and similarly for $C^\infty_c(\widetilde E_\sigma)$. 
Let $\Omega\in\cZ(\gf_\C)$ be the Casimir element of $G(\R)$ and let $R$
be the right regular representation of $G(\R)$ on $C^\infty(G(\R))$. Let 
$\widetilde\Delta_\sigma$ be the second order elliptic operator which is
induced by $-R(\Omega)\otimes\Id$ in  
$C^\infty(\widetilde E_\sigma)$.
Let $\widetilde\nabla^\sigma$ be the canonical connection on 
$\widetilde E_\sigma$,
and let $\Omega_K$ be the Casimir element of  $K_\infty$. Let $\lambda_\sigma
=\sigma(\Omega_K)$ be the Casimir eigenvalue of $\sigma$.
Then with respect to the identification (\ref{1.1}), we have
\begin{equation} \label{1.2}
(\widetilde\nabla^\sigma)^*\widetilde\nabla^\sigma=-R(\Omega)
\otimes \mbox{Id}+\lambda_\sigma\Id
\end{equation}
\cite[Proposition 1.1]{Mia},
and therefore
\begin{equation} \label{1.4}
\widetilde\Delta_\sigma=(\widetilde\nabla^\sigma)^*\widetilde\nabla^\sigma
-\lambda_\sigma\mbox{Id}.
\end{equation}
Hence $\widetilde\Delta_\sigma\colon C^\infty_c (\widetilde E_\sigma)
\to L^2(\widetilde E_\sigma)$ is essentially selfadjoint and bounded from 
below. We continue to denote its unique
selfadjoint extension by $\widetilde\Delta_\sigma$.
Let $\exp (-t\widetilde\Delta_\sigma)$ be the associated heat semigroup.  
The heat operator is a smoothing operator on $L^2(\widetilde
E_\sigma)$ which commutes with the representation of $G(\R)$ on 
$L^2(\widetilde E_\sigma)$.
Therefore, it is of the form
\begin{equation} \label{1.5}
( e^{-t\widetilde\Delta_\sigma}\varphi )(g)=
\int_{G(\R)}H_t^\sigma(g^{-1}g_1)(\varphi(g_1))dg_1, \quad
g\in G(\R),
\end{equation}
where $\varphi\in (L^2(G(\R))\otimes V_\sigma)^{K_\infty}$ and
$H_t^\sigma\colon G(\R)\to \mbox{End}(V_\sigma)$ is in $L^2\cap C^\infty$ and 
satisfies the covariance property
\begin{equation} \label{1.6}
H_t^\sigma(g)=\sigma(k)H_t^\sigma(k^{-1}gk')\sigma(k')^{-1},\quad
\mbox{for}\;g\in G(\R),\;k,k'\in K_\infty.
\end{equation}

In order to get estimates for $H_t^\sigma$, we proceed as in \cite{BM}
and relate $H_t^\sigma$ to the 
heat kernel of the Laplace operator of $G(\R)$ with respect to a left 
invariant metric on $G(\R)$. Let $\gf$ and $\kf$ denote the Lie algebras 
of $G(\R)$ and
$K_\infty$, respectively. Let $\gf=\kf\oplus \pg$ be the Cartan decomposition
and let $\theta$ be the corresponding Cartan involution.
 Let $B(Y_1,Y_2)$ be the Killing form of $\gf$. 
Set $\langle Y_1,Y_2\rangle=-B(Y_1,\theta Y_2)$, $Y_1,Y_2\in \gf$. 
By translation of $\langle\cdot,\cdot\rangle$ we get 
a left invariant Riemannian  metric on $G(\R)$. 
Let $X_1,\cdots, X_p$ be an orthonormal basis for
${\mathfrak p}$ with respect to $B|{\mathfrak p}\times{\mathfrak p}$ and
let
$Y_1,\cdots,Y_k$ be an orthonormal basis for ${\mathfrak k}$ with respect to
$-B|{\mathfrak k}\times{\mathfrak k}.$ Then we have
\begin{equation*}
\Omega=\sum^p_{i=1}X^2_i-\sum^k_{i=1}Y^2_i
\quad
\mbox{and}
\quad
\Omega_K=-\sum^k_{i=1}Y^2_i.
\end{equation*}
Let
\begin{equation} \label{1.7}
P=-\Omega+2\Omega_K=-\sum^p_{i=1}X^2_i-\sum^k_{i=1}Y^2_i.
\end{equation}
Then $R(P)$ is the Laplace operator $\Delta_{G}$ on $G(\R)$ with respect
 to the left
invariant metric defined above. The heat semigroup 
$e^{-t\Delta_{G}}$ is represented by a
smooth kernel $p_t,$ i.e.
\begin{equation} \label{1.8}
\left( e^{-t\Delta_{G}}f\right) (g)=
\int_{G(\R)} p_t(g^{-1}g')f(g')dg',\quad
f\in L^2(G(\R)),\;g\in G(\R),
\end{equation}
where $p_t\in C^\infty(G(\R))\cap L^2(G(\R))$. In fact, 
$p_t$ belongs to $L^1(G(\R))$ \cite{N} so that (\ref{1.8})
can be written as
\begin{equation*}
e^{-t\Delta_G}= R(p_t).
\end{equation*}
Let
\begin{equation*}
Q=\int_{K_\infty}R(k)\otimes \sigma (k)\;dk
\end{equation*}
be the orthogonal projection of $L^2(G(\R))\otimes V_\sigma$ onto its
$K_\infty$--invariant subspace $(L^2(G(\R))\otimes V_\sigma)^{K_\infty}$.
By (\ref{1.7}) we have
\begin{equation*}
\begin{split}
\widetilde\Delta_\sigma &   =
-Q(R(\Omega)\otimes\mbox{Id})Q\\
&= Q(R(P)\otimes\mbox{Id})Q-2Q(R(\Omega_K)\otimes \mbox{Id})Q\\
&=Q(\Delta_G\otimes\mbox{Id})Q-2\lambda_\sigma\Id_{L^2(\widetilde E_\sigma)}.
\end{split}
\end{equation*}
Hence, we get
\begin{equation*}
e^{-t\widetilde\Delta_\sigma} = Q(e^{-t\Delta_G} \otimes \mbox{Id})Q\cdot
e^{t2\lambda_\sigma}
\end{equation*}
which implies that
\begin{equation} \label{1.9}
H_t^\sigma(g)= e^{t2\lambda_\sigma} \int_{K_\infty}\int_{K_\infty} 
p_t(k^{-1}gk')\sigma(k{k'}^{-1})\;dk\;dk'.
\end{equation}
Let $\Co^1(G(\R))$ be Harish-Chandra's space of integrable rapidly decreasing
functions on $G(\R)$. Then (\ref{1.9}) can be used to show that
\begin{equation}\label{1.9a}
H_t^\sigma \in \left( \Co^1(G(\R))\otimes \mbox{End}(V_\sigma)\right)^{K_\infty
\times K_\infty}
\end{equation}
\cite[Proposition 2.4]{BM}.

Now we turn to the estimation of the derivatives of $H_t^\sigma$.
By (\ref{1.9}),  this problem can be reduced to the estimation of the 
derivatives of $p_t$. 
Let $\nabla$ denote the Levi-Civita
connection and $\rho(g,g')$ the geodesic distance of $g,g'\in G(\R)$ with 
respect
to the left invariant metric. Then all covariant derivatives of the 
 curvature tensor are bounded and the injectivity radius has a positive lower 
bound. Let $a=\dim G(\R)$, $l\in\N_0$ and $T>0$. Then it follows from  
Corollary 8 in \cite{CLY} that there exist $C,c>0$ such that
\begin{equation}\label{1.10}
\parallel\nabla^lp_t(g)\parallel\le C t^{-(a+l)/2}
\exp\left(-\frac{c\rho^2(g,1)}{t}\right)
\end{equation}
for all $0<t\le T$ and $g\in G(\R)$.  By (\ref{1.9}) and (\ref{1.10}) 
we get
\begin{equation}\label{1.11}
\begin{split}
\parallel\nabla^l H_t^\sigma(g)\parallel&\le 
e^{2t\lambda_\sigma}\int_{K_\infty}
\int_{K_\infty}
\parallel(\nabla^l p_t)(k^{-1}gk')\parallel dk dk'\\
&\le C t^{-(a+l)/2}\int_{K_\infty}\int_{K_\infty}
\exp\left(-\frac{c\rho^2(gk,k')}{t}\right) dk dk'
\end{split}
\end{equation}
for all $0<t\le T$. Choose the invariant Riemannian metric on $X$ which is 
defined by the restriction of the Killing form to $T_eX\cong\pg$. Then the 
canonical projection map $G(\R)\to X$ is a Riemannian submersion. Let $d(x,y)$
denote the geodesic distance on $X$. Then it follows that
$$\rho(g,e)\ge d(gK_\infty,K_\infty),\quad g\in G(\R).$$
Set $r(g)=d(gK_\infty,K_\infty)$, $g\in G(\R)$. 
Together with (\ref{1.11}) we get the following result.

\begin{prop}\label{p1.1}
Let $a=\dim G(\R)$,  $l\in \N_0$ and $T>0$. There exist $C,c>0$ such that
\begin{equation}\label{1.12}
\parallel\nabla^l H_t^\sigma(g)\parallel\le C t^{-(a+l)/2}
\exp\left(-\frac{cr^2(g)}{t}\right)
\end{equation}
for all $0<t\le T$ and $g\in G(\R)$. 
\end{prop}

We note that the exponent of $t$ on the right hand side of (\ref{1.12}) is not
optimal. Using the method of Donnelly \cite{Do2}, this estimate can be
improved for $l\le 1$. Indeed by Theorem 3.1 of \cite{Mu1} we have
\begin{prop}\label{p1.2}
Let $n=\dim X$ and $T>0$. There exist $C,c>0$ such that
\begin{equation}\label{1.12a}
\parallel\nabla^l H_t^\sigma(g)\parallel\le C t^{-n/2-l}
\exp\left(-\frac{cr^2(g)}{t}\right)
\end{equation}
for all $0<t\le T$, $0\le l\le1$, and $g\in G(\R)$. 
\end{prop}

We also need the asymptotic behaviour of the heat kernel on the diagonal.
It is described by the following lemma.
\begin{lem}\label{l1.2}
Let $n=\dim X$ and let $e\in G(\R)$ be the identity element. Then
$$\tr H^\sigma_t(e)=\frac{\dim(\sigma)}{(4\pi)^{n/2}}t^{-n/2}+O(t^{-(n-1)/2})$$
as $t\to0$.
\end{lem}
 
\begin{proof}
Note that for each $x\in X$, the injectivity radius at $x$ is infinite. Hence
we can construct a parametrix for the fundamental solution of the heat
equation for $\Delta_\sigma$ as in \cite{Do2}. Let $\epsilon>0$ and set
$$U_\epsilon=\{(x,y)\in X\times X\;\big|\; d(x,y)<\epsilon\}.$$
For any $l\in\N$
we define an approximate fundamental solution  $P_l(x,y,t)$ on $U_\epsilon$
by the formula
$$P_l(x,y,t)=(4\pi t)^{-n/2}\exp\left(\frac{-d^2(x,y)}{4t}\right)\left(
\sum_{i=0}^l\Phi_i(x,y)t^i\right),$$
where the $\Phi_i(x,y)$ are smooth sections of 
$E_\sigma\boxtimes E_\sigma^*$ over
$U_\epsilon\times U_\epsilon$ which are constructed recursively as in 
Theorem 2.26 of \cite{BGV}. In particular,
we have
$$\Phi_0(x,x)=\Id_{V_\sigma},\quad x\in X.$$
Let $\psi\in C^\infty(X\times X)$ be equal to 1 on $U_{\epsilon/4}$ and 0 on
$X\times X-U_{\epsilon/2}$. Set 
$$Q_l(x,y,t)=\psi(x,y) P_l(x,y,t).$$
If $l>n/2$, then the section $Q_l$ of $E_\sigma\boxtimes E_\sigma^*$ is a 
parametrix for the heat equation. Since $X$ is a Riemannian symmetric space,
we get
$$H^\sigma_t(e)=\Id_{V_\sigma}(4\pi t)^{-n/2}+O(t^{-(n-1)/2})$$
as $t\to 0$. This implies the lemma.
\end{proof}

\section{Estimations of the discrete spectrum}
\setcounter{equation}{0}

In this section we shall establish a number of facts concerning the growth of 
the discrete spectrum. 
Let 
$M=\GL_{n_1}\times\cdots\times\GL_{n_r},$
 $r\ge1$, and let 
$$M(\R)^1=M(\R)\cap M(\A)^1.$$
Then $M(\R)=M(\R)^1\cdot A_M(\R)^0$. Let $K_{M,\infty}
\subset M(\R)$ be  the standard maximal compact subgroup. Then $K_{M,\infty}$
is contained in $M(\R)^1$. Let
$$X_M=M(\R)^1/K_{M,\infty}$$
be the associated Riemannian symmetric space. Let $\Gamma_M\subset M(\Q)$ be
 an arithmetic subgroup and let $(\tau,V_\tau)$ be an
 irreducible unitary
representation of $K_{M,\infty}$ on $V_\tau$. Set
$$C^\infty(\Gamma_M\ba M(\R)^1,\tau):=
(C^\infty(\Gamma_M\ba M(\R)^1)\otimes V_\tau)^{K_{M,\infty}}.$$ 
If $\Gamma_M$ is 
torsion free, then $\Gamma_M\ba X_M$ is a Riemannian manifold and
the homogeneous vector bundle $\widetilde E_\tau$  over $X_M$, which is 
associated to $\tau$, can be pushed down to a vector bundle $E_\tau\to 
\Gamma_M\ba X_M$. Then $C^\infty(\Gamma_M\ba M(\R)^1,\tau)$ equals 
 $C^\infty(\Gamma_M\ba X_M,E_\tau)$, the
space of smooth sections of $E_\tau$. 
Define $C^\infty_c(\Gamma_M\ba M(\R)^1,\tau)$ and 
$L^2(\Gamma_M\ba M(\R)^1,\tau)$ similarly.
Let $\Omega_{M(\R)^1}$ be the Casimir element of $M(\R)^1$ and let 
$\Delta_\tau$ be the  operator  in 
$C^\infty(\Gamma_M\ba M(\R)^1,\tau)$ which is induced by 
$-\Omega_{M(\R)^1}\otimes\Id$. As unbounded operator in
$L^2(\Gamma_M\ba M(\R)^1,\tau)$ with domain 
$C^\infty_c(\Gamma_M\ba M(\R)^1,\tau)$,
 $\Delta_\tau$ is essentially selfadjoint. Let 
$L^2_{\cu}(\Gamma_M\ba M(\R)^1,\tau)$
be the subspace of cusp forms of $L^2(\Gamma_M\ba M(\R)^1,\tau)$. Then 
$L^2_{\cu}(\Gamma_M\ba M(\R)^1,\tau)$ is an invariant subspace of 
$\Delta_\tau$,
and $\Delta_\tau$ has pure point spectrum in this subspace consisting of
eigenvalues $\lambda_0<\lambda_1<\cdots$
of finite multiplicity. Let $\E(\lambda_i)$ be the eigenspace of
$\lambda_i$. Set
$$N_{\cu}^{\Gamma_M}(\lambda,\tau)
=\sum_{\lambda_i\le \lambda}\dim\E(\lambda_i).$$
Let $d=\dim X_M$ and let
$$C_d=\frac{1}{(4\pi)^{d/2}\Gamma(\frac{d}{2}+1)}$$
be Weyl's constant, where $\Gamma(s)$ denotes the Gamma function. Then Donnelly
\cite[Theorem 9]{Do} has established the following basic estimation of 
the counting function of the cuspidal spectrum.
\begin{theo}\label{th2.1} For every $\tau\in\Pi(K_{M,\infty})$ we have
$$ \limsup_{\lambda\to\infty} \frac{N_{\cu}^{\Gamma_M}(\lambda,\tau)}
{\lambda^{d/2}}\le C_d\dim(\tau)\vol(\Gamma_M\ba X_M).$$
\end{theo}
Actually, Donnelly proved this theorem only for the  case of a torsion free 
discrete
group. However, it is easy to extend his result to the general case.

We shall now reformulate this theorem in the representation theoretic 
context. Let $\xi_0$ be the trivial character of $A_M(\R)^0$ and let
$\pi\in\Pi(M(\A),\xi_0)$. 
Let $m(\pi)$ be the multiplicity with which $\pi$ occurs in the regular
representation of $M(\A)$ in $L^2(A_M(\R)^0M(\Q)\backslash M(\A)).$ Then
$\Pi_{\di}(M(\A),\xi_0)$ consists of all $\pi\in\Pi(M(\A),\xi_0)$ with
 $m(\pi)>0$.  Write
$$\pi=\pi_\infty\otimes\pi_f,$$
where $\pi_\infty\in\Pi(M(\R))$ and $\pi_f\in\Pi(M(\A_f)).$ Denote by 
${\H}_{\pi_\infty}$ (resp. $\H_{\pi_f}$) the Hilbert space of the representation $\pi_\infty$ (resp. $\pi_f$). 
Let $K_{M,f}$ be an open compact subgroup of $M(\A_f)$ and let 
$\tau\in\Pi(K_{M,\infty}).$ Denote by $\H_{\pi_\infty}(\tau)$ the 
$\tau$-isotypical subspace of 
${\H}_{{\pi}_{\infty}}$ and let $\H^{K_{M,f}}_{\pi_f}$ be the
subspace of $K_{M,f}$-invariant vectors in $\H_{\pi_f}.$ Denote by $\lambda
_\pi$ the Casimir eigenvalue of the restriction of $\pi_\infty$ to $M(\R)^1$.
Given $\lambda>0$, let
$$\Pi_{\di}(M(\A),\xi_0)_\lambda=\{\pi\in\Pi_{\di}(M(\A),\xi_0)\;\big|\;
 |\lambda_\pi|\le\lambda\}.$$
Define $\Pi_{\cu}(M(\A),\xi_0)_\lambda$ and $\Pi_{\res}(M(\A),\xi_0)_\lambda$ 
similarly. 
\begin{lem}\label{l2.2} Let $d=\dim X_M.$
For every open compact subgroup $K_{M,f}$ of $M(\A_f)$ and every $\tau\in\Pi
(K_{M,\infty})$ there 
exists $C>0$ such that 
$$\sum_{\pi\in\Pi_{\cu}(M(\A),\xi_0)_\lambda} m(\pi)\dim(\H^{K_{M,f}}_{\pi_f})
\dim(\H_{\pi_\infty}(\tau))\le C(1+\lambda^{d/2})$$
for $\lambda\ge0.$
\end{lem}
\begin{proof}
Extending the notation of \S1.4, we write $\Pi(M(\R),\xi_0)$ for the set of
representations in $\Pi(M(\R))$ whose central character is trivial on 
$A_M(\R)^0$. Given $\pi_\infty\in\Pi(M(\R),\xi_0)$, let $m(\pi_\infty)$ 
be the multiplicity with which $\pi_\infty$ occurs discretely in the regular 
representation of $M(\R)$ in $L^2(A_M(\R)^0M(\Q)\ba M(\A))^{K_{M,f}}$.
Then
\begin{equation}\label{2.1}
m(\pi_\infty)={\sum}^\prime_{\pi^\prime\in\Pi_{\cu}(M(\A),\xi_0)} m(\pi^\prime)
\dim(\H_{\pi_f^\prime}^{K_{M,f}}),
\end{equation}
where the sum is over all $\pi^\prime\in\Pi_{\di}(M(\A),\xi_0)$ such that the
 Archime\-dean 
component $\pi^\prime_\infty$ of $\pi^\prime$ equals $\pi_\infty$.

Let $\Pi_{\cu}(M(\R),\xi_0)$ be
the subset of all $\pi_\infty\in\Pi(M(\R),\xi_0)$ which are equivalent to an
irreducible subrepresentation of the regular representation of $M(\R)$
in the Hilbert space  $L^2_{\cu}(A_M(\R)^0M(\Q)\backslash M(\A))^{K_{M,f}}$.
Given  $\pi_\infty\in\Pi_{\cu}(M(\R),\xi_0)$, denote by 
$\lambda_{\pi_\infty}$ the Casimir
eigenvalue of the restriction of $\pi_\infty$ to $M(\R)^1$. For $\lambda\ge0$, 
let 
$$\Pi_{\cu}(M(\R),\xi_0)_\lambda=\{\pi_\infty\in\Pi_{\cu}(M(\R),\xi_0)\;
\big|\;|\lambda_{\pi_\infty}|\le \lambda\}.$$
Then by (\ref{2.1}), it suffices to show that for each 
$\tau\in\Pi(K_{M,\infty})$ there exists $C>0$ such that 
$$\sum_{\pi_\infty\in\Pi_{\cu}(M(\R),\xi_0)_\lambda}m(\pi_\infty)
\dim(\H_{\pi_\infty}(\tau)) \le C(1+\lambda^{d/2}).$$
To deal with this problem recall that there exist arithmetic subgroups
$\Gamma_{M,i}\subset M(\R),$ $i=1,\ldots,l,$ such that
$$M(\Q)\backslash M(\A)/K_{M,f}\cong \bigsqcup^l_{i=1}(\Gamma_{M,i}\backslash
M(\R))$$
(cf. \cite[section 9]{Mu1}). Hence
\begin{equation}\label{2.2}
L^2(A_M(\R)^0M(\Q)\backslash M(\A))^{K_{M,f}}\cong\bigoplus^l_{i=1} L^2(A_M
(\R)^0\Gamma_{M,i}\backslash M(\R))
\end{equation}
as $M(\R)$-modules.  
For each $i,$ $i=1,\ldots,l,$ and $\pi_\infty\in\Pi(M(\R))$
let $m_{\Gamma_{M,i}}(\pi_\infty)$ be the multiplicity with which $\pi_\infty$ 
occurs
discretely in the regular representation of $M(\R)$ in 
$L^2(A_M(\R)^0\Gamma_{M,i} \backslash M(\R)).$ Then
 $m(\pi_\infty)=\sum^l_{i=1}m_{\Gamma_{M,i}}(\pi_\infty)$ and
\begin{equation*}
\begin{split}
 \sum_{\pi_\infty\in\Pi_{\cu}(M(\R),\xi_0)_\lambda}&m(\pi_\infty)\dim
(\H_{\pi_\infty}(\tau))\\
&=\sum^l_{i=1}
\sum_{\pi_\infty\in\Pi_{\cu}(M(\R),\xi_0)_\lambda} m_{\Gamma_{M,i}}(\pi_\infty)
\dim(\H_{\pi_\infty}(\tau)).
\end{split}
\end{equation*}

The interior sum can be interpreted as follows. Fix $i$ and set 
$\Gamma_M:=\Gamma_{M,i}$.  
Let $\lambda_1<\lambda_2<\cdots$ be the eigenvalues of $\Delta_\tau$ in
the space of cusp forms $L^2_{\cu}(\Gamma_M\ba M(\R)^1,\tau)$ and let 
$\E(\lambda_i)$ be the eigenspace of $\lambda_i.$
By Frobenius reciprocity  it follows that
$$\dim\E(\lambda_i)=\sum_{-\lambda_{\pi_\infty}=\lambda_i}
m_{\Gamma_M}(\pi_\infty),$$
where the sum is over all $\pi_\infty\in\Pi_{\cu}(M(\R),\xi_0)$ such that the
Casimir eigenvalue $\lambda_{\pi_\infty}$ equals $-\lambda_i$. Hence we obtain
$$\sum_{\pi_\infty\in\Pi_{\cu}(M(\R),\xi_0)_\lambda}
m_{\Gamma_M}(\pi_\infty)\dim(\H_{\pi_\infty}(\tau))
=N_{\cu}^{\Gamma_M}(\lambda,\tau).$$
Combined with Theorem \ref{th2.1} the desired estimation follows.
\end{proof}

Next we consider the residual spectrum.
\begin{lem}\label{l2.3} Let $d=\dim X_M.$ For every open compact subgroup
$K_{M,f}$ of $M(\A_f)$ and every $\tau\in\Pi(K_{M,\infty })$ there exists
$C>0$ such that
$$\sum_{\pi\in\Pi_{\res}(M(\A),\xi_0)_\lambda}
m(\pi)\dim({\H}^{K_{M,f}}_{\pi_f})\dim(\H_{\pi_\infty}(\tau))
\le C(1+\lambda^{(d-1)/2})
$$
for $\lambda\ge0.$
\end{lem}
\begin{proof} We can assume that  $M=\GL_{n_1}\times\cdots \times\GL_{n_r}$.
Let $K_{M,f}$ be an open compact subgroup of $M(\A_f)$. There exist open
compact subgroups $K_{i,f}$ of $\GL_{n_i}(\A_f)$ such that 
$K_{1,f}\times\cdots\times K_{r,f}\subset K_{M,f}.$
Thus we can replace $K_{M,f}$ by $K_{1,f}\times\cdots\times K_{r,f}$.
Next observe that $K_{M,\infty}=\rO(n_1)\times\cdots\times\rO(n_r)$ and 
therefore, $\tau$ is given as $\tau=\tau_1\otimes\cdots\otimes\tau_r$, where
each $\tau_i$ is an irreducible unitary representation of $\rO(n_i)$.
Finally note that every
$\pi\in\Pi(M(\A),\xi_0)$ is of the form $\pi=\pi_1\otimes\cdots\otimes\pi_r$.
Hence we get $m(\pi)=\prod_{i=1}^r m(\pi_i)$ and
$$\dim\bigl(\H_{\pi_f}^{K_{M,f}}\bigr)=\prod_{i=1}^r\dim\bigl(\H_{\pi_{i,f}}
^{K_{i,f}}\bigr),\quad \dim\bigl(\H_{\pi_\infty}(\tau)\bigr)
=\prod_{i=1}^r\dim\bigl(\H_{\pi_{i,\infty}}(\tau_i)\bigr).$$
This implies immediately that it suffices to consider a single factor.

With the analogous notation the proof of the proposition is reduced to the
following problem. For  $m\in\N$ set $X_m=\SL_m(\R)/\SO(m)$ and 
$d_m=\dim X_m$. Then we need to show that for every open compact subgroup
$K_{m,f}$ of $\GL_m(\A_f)$ and every $\tau\in \Pi(\rO(m))$ there exists 
$C>0$ sucht that
$$\sum_{\pi\in\Pi_{res}(\GL_m(\A),\xi_0)_\lambda}m(\pi)\dim(\H^{K_{m,f}}
_{\pi_f})\dim(\H_{\pi_\infty}(\tau))\le C(1+\lambda^{(d_m-1)/2})$$
for $\lambda\ge 0.$ 
To deal with this problem  recall the description of the residual 
spectrum of $\GL_m$ by 
M{\oe}glin and Waldspurger \cite{MW}. Let
 $\pi\in\Pi_{\mbox{\k res}}(\GL_m(\A))$ 
and suppose that $\pi$ is trivial on $A_{\GL_m}(\R)^0.$ 
There exist $k|m$, a standard parabolic subgroup $P$ of $\GL_m$ of type
$(l,\ldots,l)$, $l=m/k$, and a cuspidal automorphic representation $\rho$ of
 $\GL_l$ which is trivial on $A_{\GL_l}(\R)^0$, such that $\pi$ is
equivalent to the unique irreducible quotient $J(\rho)$ of the induced
 representation
$$I^{\GL_m(\A)}_{P(\A)} (\rho[(k-1)/2]\otimes\cdots\otimes\rho[(1-k)/2]).$$
Here $\rho[s]$ denotes the representation $g\longmapsto \rho(g)|\det g|^s,$ 
$s\in\C$.
At the Archimedean place, the corresponding induced representation 
$$I^{\GL_m}_P(\rho_\infty,k):=I^{\GL_m(\R)}_{P(\R)} (\rho_\infty[(k-1)/2]
\otimes\cdots\otimes\rho_\infty[(1-k)/2])$$
has also a unique irreducible quotient $J(\rho_\infty)$. Comparing the 
definitions, we get $J(\rho)_\infty=J(\rho_\infty)$.  
Hence the Casimir eigenvalue of 
$\pi_\infty=J(\rho)_\infty$ equals the Casimir
eigenvalue of $J(\rho_\infty)$
 which in turn coincides with the Casimir eigenvalue
of the induced representation $I^{\GL_m}_P(\rho_\infty,k)$. 
Let $\lambda_\rho$ be the Casimir eigenvalue of $\rho_\infty$. Then it follows
that there exists $C>0$ such that
$|\lambda_\pi-k\lambda_\rho|\le C$
for all $\pi\in\Pi_{\res}(\GL_m(\A),\xi_0)$. Using  the main
theorem of \cite[p.606]{MW} it follows that it suffices to fix $l|m$, 
$l<m$, and to estimate
\begin{equation}\label{2.3}
\sum_{\rho\in\Pi_{\cu}(\GL_l(\A),\xi_0)_\lambda}m(\rho)\dim\bigl(\H_{J(\rho)_f}^
{K_{m,f}}\bigr)\dim\bigl(\H_{J(\rho)_\infty}(\tau)\bigr).
\end{equation}
First note that by \cite{Sk}, we have $m(\rho)=1$ for all $\rho\in
\Pi_{\cu}(\GL_l(\A),\xi_0)$. So it remains to estimate the dimensions. We
 begin with
the infinite place.  Observe that
$\dim(\H_{J(\rho)_\infty}(\tau))=\dim(\tau)[J(\rho_\infty)|_{\rO(m)}:\tau].$
Thus in order to estimate $\dim(\H_{J(\rho)_\infty}(\tau))$ it
suffices to estimate the multiplicity $[J(\rho_\infty)|_{\rO(m)}:\tau]$. Since $J(\rho_\infty)$
is an irreducible quotient of $I^{\GL_m}_P(\rho_\infty,k)$, we have
$$[J(\rho_\infty)|_{\rO(m)}:\tau]\le 
[I^{\GL_m}_P(\rho_\infty,k)|_{\rO(m)}:\tau].$$
Let $K_{l,\infty}=\rO(l)\times\cdots\times\rO(l)$. Using Frobenius 
reciprocity as in \cite[p.208]{Kn}, we obtain
\begin{equation*}
\begin{split}
[I^{\GL_m}_P(\rho_\infty,k)&|_{\rO(m)}:\tau]\\
&=\sum_{\omega\in\Pi(K_{l,\infty})}
[(\rho_\infty\otimes\cdots\otimes\rho_\infty)|_{K_{l,\infty}}:\omega]
\cdot[\tau|_{K_{l,\infty}}:\omega].
\end{split}
\end{equation*}
Finally note that $\omega=\omega_1\otimes\cdots\otimes\omega_k$ with 
$\omega_i\in\Pi(\rO(l))$. Therefore we have
$$[(\rho_\infty\otimes\cdots\otimes\rho_\infty)|_{K_{l,\infty}}:\omega]=
\prod_{i=1}^k[\rho_\infty|_{\rO(l)}:\omega_i].$$
At the finite places we proceed in an analogous way. This implies
that there exist an open compact subgroup $K_{l,f}$ of $\GL_l(\A_f)$ and
$\omega_1,...,\omega_p\in\Pi(\rO(l))$ such that (\ref{2.3}) is bounded
from above by a constant times
$$\sum_{i=1}^p\left(\sum_{\rho\in\Pi_{\cu}(\GL_l(\A),\xi_0)_\lambda}
m(\rho)\dim\bigl(\H_{\rho_f}^{K_{l,f}}\bigr)
\dim\bigl(\H_{\rho_\infty}(\omega_i)\bigr)\right)^k.
$$
By Lemma \ref{l2.2} this term is bounded by a constant times 
$(1+\lambda^{d_l/2})^k$, where $d_l=l(l+1)/2-1$.
Since $m=k\cdot l$ and $k>1$, we have 
$$d_lk=\frac{l(l+1)k}{2}-k\le \frac{m(m+1)}{2}-2=d_m-1.$$
This proves the desired estimation in the 
case of $M=\GL_m$, and as explained above, this suffices to prove the lemma.
\end{proof}
Combining Lemma \ref{l2.2} and Lemma \ref{l2.3}, we obtain
\begin{prop}\label{p2.4} Let $d=\dim X_M$. For every open compact subgroup
$K_{M,f}$ of $M(\A_f)$ and every $\tau\in\Pi(K_{M,\infty})$ there exists
$C>0$ such that
$$\sum_{\pi\in\Pi_{\di}(M(\A),\xi_0)_\lambda}
m(\pi)\dim(\H_{\pi_f}^{K_{M,f}})\cdot\dim(\H_{\pi_\infty}(\tau))
\le C(1+\lambda^{d/2})$$
for $\lambda\ge0$. 
\end{prop}
Next we restate Proposition \ref{p2.4} in terms of dimensions of spaces 
of automorphic forms. Let $P\in\cP(M)$ and 
let $\cA^2(P)$ be the space of square integrable automorphic forms on 
$N_P(\A)M_P(\Q)A_P(\R)^0\ba G(\A)$. Given
$\pi\in\Pi_{\di}(M(\A),\xi_0)$ let $\cA^2_\pi(P)$ be the subspace of 
$\cA^2(P)$ of automorphic forms of type $\pi$ \cite[p.925]{A1}. 
Let $K_\infty$ be the standard
maximal compact subgroup of $G(\R)$. 
Given  an open compact subgroup $K_f$ of $G(\A_f)$ and
$\sigma\in\Pi(K_\infty)$, let ${\cA}_\pi(P)_{K_f}$ denote the subspace of 
$K_f$-invariant automorphic forms in ${\cA}^2_\pi(P)$ and let ${\cA}^2_\pi(P)
_{K_f,\sigma}$ be the $\sigma$-isotypical subspace of ${\cA}^2_\pi(P)_{K_f}$.

\begin{prop}\label{p2.5} Let $d=\dim X_M$. For every open compact subgroup
$K_f$ of $G(\A_f)$ and every $\sigma\in\Pi(K_\infty)$ there exists $C>0$ such
that
$$\sum_{\pi\in\Pi_{\di}(M(\A),\xi_0)_\lambda}
\dim{\cA}^2_\pi(P)_{K_{f},\sigma}\le C(1+\lambda^{d/2})$$
for $\lambda\ge0$.
\end{prop}
\begin{proof} Let $\pi\in\Pi_{\di}(M(\A),\xi_0)$.  Let $\H_P(\pi)$
be the Hilbert space of the induced representation $I^{G(\A)}_{P(\A)}(\pi).$ 
There is a canonical isomorphism
\begin{equation}\label{2.4}
j_P:\H_P(\pi)\otimes \Hom_{M(\A)}(\pi,I^{M(\A)}_{M(\Q)A_M(\R)^0}
(\xi_0))\to \overline{\cA}^2_\pi(P),
\end{equation}
which intertwines the induced representations.
Let $\pi=\pi_\infty\otimes\pi_f.$ Let $\H_P(\pi_\infty)$ (resp. $\H_P(\pi_f)$)
be the Hilbert space of the induced representation
$I^{G(\R)}_{P(\R)}(\pi_\infty))$ 
(resp. $I^{G(\A_f)}_{P(\A_f)}(\pi_f)$). 
Denote by $\H_P(\pi_\infty)_\sigma$ the $\sigma$-isotypical subspace of 
$\H_P(\pi_\infty)$ and by $\H_P(\pi_f)^{K_f}$ the subspace of $K_f$-invariant 
vectors
of $\H_P(\pi_f).$ Then it follows from (\ref{2.4}) that 
\begin{equation}\label{2.5}
\dim{\cA}^2_\pi(P)_{K_f,\sigma}=m(\pi)\dim(\H_P(\pi_f)^{K_f})\dim(\H_P
(\pi_\infty)_\sigma).
\end{equation}
Using Frobenius reciprocity as in \cite[p.208]{Kn} we get
$$[I^{G(\R)}_{P(\R)}(\pi_\infty)|_{K_\infty}:\sigma]=
\sum_{\tau\in\Pi(K_{M,\infty})}[\pi_\infty |_{K_{M,\infty}}:\tau]\cdot
[\sigma|_{K_{M,\infty}}:\tau].$$
Hence we get
\begin{equation}\label{2.6}
\dim(\H_P(\pi_\infty)_\sigma)\le\dim(\sigma)\sum_{\tau\in\Pi(K_{M,\infty})}
\dim(\H_{\pi_\infty}(\tau))[\sigma|_{K_{M,\infty}}:\tau].
\end{equation}
Next we consider $\pi_f=\otimes_{p<\infty}\pi_p$. Replacing $K_f$ by a subgroup
of finite index if necessary, we can assume that $K_f=\Pi_{p<\infty}K_p$. 
For any $p<\infty$, denote by $\H_P(\pi_p)$ the Hilbert space of the induced
 representation 
$I^{G(\Q_p)}_{P(\Q_p)}(\pi_p)$. Let $\H_P(\pi_p)^{K_p}$ be the subspace
of $K_p$-invariant vectors. Then
$\dim\H_P(\pi_p)^{K_p}=1$ for alomost all $p$ and
$$\H_P(\pi_f)^{K_f}\cong\bigotimes_{p<\infty}\H_P(\pi_p)^{K_p}.$$
Furthermore
\begin{equation}\label{2.7}
\begin{split}
I^{G(\Q_p)}_{P(\Q_p)}(\pi_p)^{K_p}&=\left(I^{G(\Z_p)}_{P(\Z_p)}(\pi_p)
\right)^{K_p} \\
&\hookrightarrow \bigoplus_{G(\Z_p)/K_p}I^{K_p}_{K_p\cap P}(\pi_p)^{K_p} \\
&  \cong\bigoplus_{G(Z_p)/K_p}\pi^{K_p\cap P}_p.
\end{split}
\end{equation}
Let $K_{M,f}=K_f\cap M(\A_f)$. Using (\ref{2.5})--(\ref{2.7}), it
follows that in order to prove the proposition, it suffices to fix
$\tau\in\Pi(K_{M,\infty})$ and to estimate
$$\sum_{\pi\in\Pi_{\di}(M(\A),\xi_0)_\lambda}m(\pi)\dim(\H_{\pi_f}^{K_{M,f}})
\dim(\H_{\pi_\infty}(\tau)).$$
The proof is now completed applying Proposition \ref{p2.4}.
\end{proof}

Finally we consider the analogous statement of Lemma \ref{l2.3} 
 at the Archimedean place. For simplicity we consider only the case
$M=G$. Let  $K_\infty$ be the standard maximal compact subgroup of $G(\R)$. 
Let $\Gamma\subset G(\Q)$ be an arithmetic subgroup and 
$\sigma\in\Pi(K_\infty)$. Then the discrete subspace 
$L^2_{\di}(\Gamma\ba G(\R)^1,\sigma)$ of $\Delta_\sigma$ decomposes as
$$L^2_{\di}(\Gamma\ba G(\R)^1,\sigma)=
 L^2_{\cu}(\Gamma\ba G(\R)^1,\sigma)\oplus 
L^2_{\res}(\Gamma\ba G(\R),\sigma),$$
where $L^2_{\res}(\Gamma\ba G(\R)^1,\sigma)$ is the subspace which 
corresponds to the residual spectrum of $\Delta_\sigma$. Let
$$L^2_{\res}(\Gamma\ba G(\R)^1,\sigma)=\bigoplus_{i}\E_{\res}(\lambda_i)$$
be the decomposition into eigenspaces of $\Delta_\sigma$. For $\lambda\ge0$ set
$$N_{\res}^\Gamma(\lambda,\sigma)
=\sum_{\lambda_i\le\lambda}\dim\E_{\res}(\lambda_i).$$
\begin{prop}\label{p2.6} Let $d=G(\R)^1/K_\infty$. Let $\Gamma\subset G(\Q)$ 
be an
arithmetic subgroup. For every 
$\sigma\in\Pi(K_{\infty})$ there 
exists $C>0$ such that
$$N_{\res}^\Gamma(\lambda,\sigma)\le C(1+\lambda^{(d-1)/2})$$ 
for $\lambda\ge0.$
\end{prop}
\begin{proof} First assume that $\Gamma\subset \SL_n(\Z)$.  Let
$\Gamma(N)\subset \Gamma$ be a congruence subgroup. Then
\begin{equation}\label{2.8} 
N_{\res}^{\Gamma}(\lambda,\sigma)\le N_{\res}
^{\Gamma(N)}(\lambda,\sigma).
\end{equation}
Let
$$N=\Pi_p p^{r_p},\quad r_p\ge 0.$$
Set
$$K_p(N)=\{k\in \GL_n(\Z_p)\mid k\equiv 1\;\mbox{mod}\; p^{r_p}\Z_p\}$$
and
\begin{equation}\label{2.9}
K(N)=\Pi_{p<\infty}K_p(N).
\end{equation}
Then $K(N)$ is an open compact subgroup of $G(\A_f)$ and
\begin{equation}\label{2.10}
A_G(\R)^0G(\Q)\backslash G(\A)/K(N)\cong\bigsqcup_{(\Z/N\Z)^*}
(\Gamma(N)\backslash  \SL_n(\R))
\end{equation}
(cf. \cite{A9}). Hence
$$L^2_{\mbox{\k res}}(A_G(\R)^0G(\Q)\backslash G(\A))^{K(N)}\cong
\bigoplus_{(\Z/N\Z)^*}L^2_{\mbox{\k res}}(\Gamma(N)\backslash \SL_n(\R))$$
as $\SL_n(\R)$-modules.
Put $M=G$ in Lemma \ref{l2.3}. Then 
$$\sum_{\pi\in\Pi_{res}(G(\A),\xi_0)_\lambda} m(\pi)\dim(\H^{K_f}_{\pi_f})
\dim(\H_{\pi_\infty}(\sigma))=\varphi(N) N_{\res}^{\Gamma(N)}
(\lambda,\sigma),$$
where $\varphi(N)=\#[(\Z/N\Z)^*]$.
Hence by Lemma \ref{l2.3} it follows that there exists $C>0$ such that
$$N_{\res}^{\Gamma(N)}(\lambda,\sigma)\le C(1+\lambda^{(d-1)/2}).$$
This proves the proposition for $\Gamma\subset \SL_n(\Z)$. Since an 
arithmetic subgroup $\Gamma\subset G(\Q)$ is commensurable with $G(\Z)$,
the general case can be easily reduced to this one.

\end{proof}

\section{Rankin-Selberg $L$-functions}
\setcounter{equation}{0}

The main purpose of this section is  to prove  estimates for the
number of  zeros of Rankin-Selberg $L$-functions. We shall consider the 
Rankin-Selberg $L$-functions over an arbitrary number field, although in the
present paper we shall use them  only in the case of $\Q$.
We begin with the description of the local $L$-facors. 

Let $F$ be a local field of characteristic zero. 
Recall that any irreducible admissible representation of $\GL_m(F)$ is given
as a Langlands quotient: There exist a standard parabolic subgroup $P$ of type
$(m_1,...,m_r)$, discrete series representations $\delta_i$ of $\GL_{m_i}(F)$
and complex numbers $s_1,...,s_r$ satisfying $\Re(s_1)\ge \Re(s_2)\ge\cdots\ge \Re(s_r)$ 
such that
\begin{equation}\label{4.1n}
\pi=J^{\GL_m}_P(\delta_1[s_1]\otimes\cdots\otimes\delta_r[s_r]),
\end{equation}
where the representation on the right is the unique irreducible quotient of
 the induced representation
$I^{\GL_m}_P(\delta_1[s_1]\otimes\cdots\otimes\delta_r[s_r])$ \cite[I.2]{MW}.
 Furthermore
any irreducible generic representation $\pi$ of $\GL_m(F)$ is equivalent to
a fully induced representation 
$I^{\GL_m}_P(\delta_1[s_1]\otimes\cdots\otimes\delta_r[s_r])$. If  $\pi$
is generic and unitary, it follows from the classification of the unitary
 dual of 
$\GL_m(F)$ that the parameters $s_i$ satisfy
\begin{equation}\label{4.2n}
|\Re(s_i)|<1/2,\quad i=1,...,r.
\end{equation}
Suppose that $\pi$ is given as Langlands quotient of the form (\ref{4.1n}).
 Then the $L$-function satisfies
\begin{equation}\label{4.3n}
L(s,\pi)=\prod_{j}L(s+s_j,\delta_j)
\end{equation}
\cite{J}. Furthermore, suppose that $\pi_1$ and $\pi_2$ are irreducible
admissible representations of
$G_1=\GL_{m_1}(\R)$ and $G_2=\GL_{m_2}(\R)$, respectively. Let 
\begin{equation*}
\pi_i\cong J^{\GL_{n_{i}}}_{P_i}(\tau_{i1}[s_{i1}],\ldots,\tau_{ir_i}
[s_{ir_{i}}])
\end{equation*}
be the Langlands parametrizations of $\pi_i$, $i=1,2$.
Then it follows from the multiplicativity of the local Rankin-Selberg
$L$-factors \cite[(9.4)]{JPS}, \cite{Sh6} that
\begin{equation}\label{4.4n}
L(s,\pi_1\times\pi_2)= \prod^{r_1}_{i=1}\prod^{r_2}_{j=1}L
(s+s_{1i}+s_{2j},\tau_{1i}\times\tau_{2j}).
\end{equation}
This reduces the description of the local $L$-factors to the 
square-inte\-grable case. Now we distinguish three cases according to the 
type of the field.

{\bf 1.} $F$ {\it non-Archimedean}

Let $\cO_F$
denote the ring of integers of $F$ and $\Pg$ the maximal ideal of $\cO_F$.
Set $q=\cO_F/\Pg$. The square-integrable case can be further reduced to the
supercuspidal one. Finally for supercuspidal representations 
the $L$-factor is given by an 
elementary polynomial in $q^{-s}$. For details see \cite{JPS}
 (see also \cite{MS}). If we put 
together all steps of the reduction, we get the following result.
Let $\pi_1$ and $\pi_2$ be  irreducible admissible  representations  of
$\GL_{n_1}(F)$ and $\GL_{n_{2}}(F)$, resprectively. Then
 there is a polynomial 
$P_{\pi_1,\pi_2}(x)$ of degree at most $n_1\cdot n_2$ with
$P_{\pi_1,\pi_2}(0)=1$ such that
$$L(s,\pi_1\times\pi_2)=P_{\pi_1,\pi_2}(q^{-s})^{-1}.$$
In the special case 
where $\pi_1$ and $\pi_2$ are unitary and generic the $L$-factor has the
following special form.
\begin{lem}\label{l3.1} Let $\pi_1$ and $\pi_2$ be irreducible unitary generic
representations of $\GL_{n_1}(F)$ and $\GL_{n_2}(F)$, respectively. There
exist complex numbers $a_{i}$, $i=1,...,n_1\cdot n_2$,
with $|a_{i}|<q$ such that

\begin{equation}\label{3.1}
L(s,\pi_1\times\pi_2)=\prod^{n_1\cdot n_2}_{i=1}
(1-a_{i}q^{-s})^{-1}.
\end{equation}
\end{lem}

\begin{proof} 
Let $\delta_1$ and $\delta_2$ be square-integrable 
representations of
$\GL_{d_{1}}(F)$ and $\GL_{d_{2}}(F)$, respectively. As explained above
 there is a polynomial $P_{\delta_1,\delta_2}(x)$ of degree at most 
$d_1\cdot d_2$ with $P_{\delta_1,\delta_2}(0)=1$ such that 
$$L(s,\delta_1\times\delta_2)=P_{\delta_1,\delta_2}(q^{-s})^{-1}.$$
By (6) of \cite{JPS},
p. 445, $L(s,\delta_1\times\delta_2)$ is holomorphic in the half-plane
$\Re(s)>0$. Hence $P_{\delta_1,\delta_2}(x)$ has no zeros in the unit disc.
Thus there exist complex numbers $b_{i}$ with $|b_{i}|<1$ 
such that
\begin{equation}\label{4.6n}
L(s,\delta_1\times\delta_2)=\prod_{i=1}^{d_1\cdot d_2}(1-b_{i}q^{-s})^{-1}.
\end{equation}
Now let $\pi_1$ and $\pi_2$ be unitary and generic. Then 
$L(s,\pi_1\times\pi_2)$ can be written as a product of the form (\ref{4.4n})
 and by (\ref{4.2n}) the parameters $s_{ij}$ satisfy 
$|\Re(s_{ij})|<1/2$, $i=1,2$, $j=1,...,r_i$. Combined with (\ref{4.6n}) the 
lemma follows.
\end{proof}

If $F$ is Archimedean the $L$-factors are defined in terms of the $L$-factors
attached to semisimple representations of the Weyl group $W_F$ by means of the
Langlands correspondence \cite{La1}. The structure of the $L$-factors
 are described, for example, in \cite[\S3]{MS}. We briefly recall the result.

\medskip
{\bf 2.} $F=\R$.

First note that $\GL_m(\R)$ does not have square-integrable representations if
$m\ge3$.  To describe the 
principal $L$-factors in the remaining cases $d=1$ and $d=2$, we define Gamma 
factors by
\begin{equation}\label{3.5}
\Gamma_\R(s)=\pi^{-s/2}\Gamma\left(\frac{s}{2}\right),\quad \Gamma_\C(s)=2(2\pi)^{-s}\Gamma(s).
\end{equation}
In the case $d=1$, the unitary representations of $\GL_1(\R)=\R^\times$ are
of the form $\psi_{\epsilon,t}(x)=\sign^{\epsilon}(x)|x|^t$ with 
$\epsilon\in\{0,1\}$ and
$t\in i\R$. Then
$$L(s,\psi_{\epsilon,t})=\Gamma_\R(s+t+\epsilon).$$
For $k\in\Z$ let $D_k$ be the $k$-th discrete series representation of 
$\GL_2(\R)$ with the same infinitesimal character as the $k$-dimensional
representation. 
Then the unitary square-integrable representations of $\GL_2(\R)$  are unitary
 twists of  $D_k$, $k\in\Z$,   for
 which the $L$-factor is given by
$$L(s,D_k)=\Gamma_\C(s+|k|/2).$$
Let $\psi_\epsilon=\sign^\epsilon$, $\epsilon\in\{0,1\}$.
Then up to twists by unramified characters the following list describes  the 
Rankin-Selberg $L$-factors in the square-integrable case:
\begin{equation}\label{3.7}
\begin{aligned}
L(s,D_{k_1}\times D_{k_2})&=\Gamma_\C(s+|k_1-k_2|/2)\cdot
\Gamma_\C(s+|k_1+k_2|/2)\\
L(s,D_k\times\psi_\epsilon)&=L(s,\psi_\epsilon\times D_k)=\Gamma_\C(s+|k|/2)\\
L(s,\psi_{\epsilon_1}\times\psi_{\epsilon_2})&= \Gamma_\R((s+\epsilon_{1,2})), 
\end{aligned}
\end{equation}
where $0 \leq \epsilon_{1,2}\leq 1$ with $ \epsilon_{1,2}\equiv
\epsilon _1 + \epsilon _2 \mbox{ mod }2$. 

\medskip

{\bf 3.} $F=\C.$

There exist square-integrable representations of $\GL_k(\C)$  only 
if $k=1$. For $r\in\Z$ let $\chi_r$  be the character of 
$\GL_1(\C)=\C^\times$ which is given by $\chi(z)=(z/\ov{z})^{r}$, $z\in\C^*$.
Then
\begin{equation}\label{4.9n}
L(s,\chi_r)=\Gamma_\C(s+|r_i|/2).
\end{equation}
If $\chi_{r_1}$ and $\chi_{r_2}$ are two characters as above, then we have
$$L(s,\chi_{r_1}\times\chi_{r_2})=\Gamma_\C(s+|r_1+r_2|/2).$$
Up to twists by unramified characters, these are all possibilities for the 
$L$-factors in the square-integrable case.

To summarize we obtain the following  description of the local $L$-factors in 
the complex case.
Let $\pi$ be an irreducible unitary representation of $\GL_m(\C)$. It  
is given by a Langlands quotient of the form
$$\pi=J^{\GL_m}_B(\chi_1[s_1]\otimes\cdots\otimes\chi_m[s_m]),$$
where $B$ is the standard Borel subgroup of $\GL_m$ and the 
$\chi_i$'s are characters of $\GL_1(\C)=\C^\times$ which are defined by
$\chi(z)=(z/\ov z)^{r_i}$,  $r_i\in\Z$, $i=1,\ldots,m$. Then
\begin{equation}\label{4.10n}
L(s,\pi)=\prod^m_{i=1}\Gamma_\C(s+s_i+|r_i|/2).
\end{equation}
Let $\pi_1$ and $\pi_2$ be irreducible unitary representations of 
$\GL_{m_1}(\C)$ and $\GL_{m_2}(\C),$ respectively.
Let $B_i\subset\GL_{m_i}$ be the standard Borel subgroup. There exist 
characters $\chi_{ij}$ of $\C^\times$ of the form
$\chi_{ij}(z)=(z/\ov z)^{r_{ij}}$,
$r_{ij}\in\Z$, and complex numbers $s_{ij}$, 
$i=1,\ldots,m_1$, $j=1,\ldots,m_2$, satisfying
$$\Re(s_{i1})\ge\cdots\ge\Re(s_{im_i}),\quad |\Re(s_{ij}|<1/2,$$
 such that 
\begin{equation}\label{4.11n}
\pi_i=J^{\GL_{m_i}}_{B_i}(\chi_{i1}[s_{i1}]\otimes\cdots\otimes\chi_{im_i}
[s_{im_i}]),\quad i=1,2.
\end{equation}
Then the Rankin-Selberg $L$-factor is 
given by
\begin{equation}\label{4.12n}
L(s,\pi_1\times\pi_2)
=\prod^{m_1}_{i=1}\prod^{m_2}_{j=1}
\Gamma_\C(s+s_{1i}+s_{2j}+|r_{1i}+r_{2j}|/2).
\end{equation}

If $F=\R$, the $L$-factors have a similar form.

The description of the $L$-factors in the Archimedean case can be
unified in the following way. By the duplication formula of the Gamma function
we have
\begin{equation}\label{4.13n}
\Gamma_\C(s)=\Gamma_\R(s)\Gamma_\R(s+1).
\end{equation}
Let $F$ be Archimedean. Set $e_F=1$, if $F=\R$, and $e_F=2$, if $F=\C.$ Let
$\pi\in\Pi(\GL_m(F)).$ Then it follows from (\ref{4.13n})
and the definition of the $L$-factors, that there exist complex numbers 
$\mu_j(\pi),$ $j=1,\ldots,me_F,$ such that
\begin{equation}\label{4.14n}
L(s,\pi)=\prod^{me_F}_{j=1}\Gamma_\R(s+\mu_j(\pi)).
\end{equation}
The numbers $\mu_j(\pi)$ are determined by the Langlands parameters of $\pi.$
Similarly, if $\pi_i\in\Pi(\GL_{m_i}(F)),$ $i=1,2,$ it follows from the
description of the Rankin-Selberg $L$-factors that there exist complex numbers
$\mu_{j,k}(\pi_1\times\pi_2)$ such that
\begin{equation}\label{3.12}
L(s,\pi_1\times\pi_2)=\prod_{j,k}\Gamma_\R(s+\mu_{j,k}(\pi_1\times\pi_2)).
\end{equation}

\begin{lem}\label{l3.2}
Let $F$ be Archimedean. There exists $C>0$ such that
$$\sum_{j,k}|\mu_{j,k}(\pi_1\times\pi_2)|^2\le C\left( \sum_i|\mu_i(\pi_1)|^2+
\sum_j|\mu_j(\pi_2)|^2\right)$$
for all generic $\pi_i\in\Pi(\GL_{m_i}(F)),$ $i=1,2.$ 
\end{lem}

\begin{proof} First consider the case $F=\C$. Let $\pi_1$ and $\pi_2$
be irreducible unitary generic representations of $\GL_{m_1}(\C)$ and
$\GL_{m_2}(\C),$ respectively.
Write  $\pi_i$ as Langlands quotient of the form (\ref{4.11n}).
Using (\ref{4.10n}) and (\ref{4.12n}) together with (\ref{4.13n}), it follows
that it suffices to prove that there exist $C>0$ such that
\begin{equation*}
\sum_{j,k}|s_{1j}+s_{2k}+|r_{1j}+r_{2k}|/2|^2 
\le C\sum_{i,j}|s_{ij}+|r_{ij}|/2|^2
\end{equation*}
for all generic $\pi_i\in\Pi(\GL_{m_i}(\C))$, $i=1,2$. This follows 
immediately, if we use the fact that the parameters $s_{ij}$ satisfy
$|\Re(s_{ij})|<1/2$ and the $r_{ij}$'s are integers.

The proof in the case $F=\R$ is essentially the same. We only have to check
the different possible cases for the $L$-factors as listed above.
\end{proof}

Next we consider the global Rankin-Selberg $L$-functions. Let $E$ be a number
filed  and let $\A_E$ be the ring of ad\`eles of $E$. Given $m\in\N$, let
$\Pi_{\di}(\GL_m(\A_E))$ and $\Pi_{\cu}(\GL_m(\A_E))$ be defined in the same
way as in the case of $\Q$ (see \S1.4).
Recall that the Rankin-Selberg $L$-function attached to a pair of
automorphic representations $\pi_1$ of $\GL_{m_1}(\A_E)$ and $\pi_2$ of 
 $\GL_{m_2}(\A_E)$ is defined by the Euler product
\begin{equation}\label{3.13}
L(s,\pi_1\times\pi_2)=\prod_v L(s,\pi_{1,v}\times\pi_{2,v}),
\end{equation}
where $v$ runs over all places of $E$. The Euler product is known to converge 
in a certain half-plane $\Re(s)>c$. Suppose that $\pi_1$ and $\pi_2$ are
unitary cuspidal automorphic representations of 
$\GL_{m_1}(\A_E)$ and $\GL_{m_2}(\A_E)$, respectively.
Then $L(s,\pi_1\times\pi_2)$ has the following basic properties:
\begin{enumerate}
\item[i)]  The Euler product $L(s,\pi_1\times\pi_2)$ converges absolutely for
 all $s$ in the half-plane $\Re(s)>1$.
\item[ii)] $L(s,\pi_1\times\pi_2)$ admits a meromorphic continuation to the
entire complex plane with at most simple poles at 0 and 1. 
\item[iii)] $L(s,\pi_1\times\pi_2)$ is of order one and is
bounded in vertical strips outside of the poles.
\item[iv)] It satisfies a functional equation of the form
\begin{equation}\label{3.14}
L(s,\pi_1\times\pi_2)=\epsilon(s,\pi_1\times\pi_2)L(1-s,\widetilde\pi_1\times
\widetilde\pi_2)
\end{equation}
with
\begin{equation}\label{3.15}
\epsilon(s,\pi_1\times\pi_2)=W(\pi_1\times\pi_2)(D^{m_1m_2}_EN(\pi_1\times
\pi_2))^{1/2-s},
\end{equation}
where $D_E$ is the discriminant of $E$, $W(\pi_1\times\pi_2)$ is a complex
number of absolute value $1$ and $N(\pi_1\times\pi_2)\in\N$. 
\end{enumerate}

The absolute convergence of the Euler product (\ref{3.13})
in the half-plane $\Re(s)>1$ was proved in \cite{JS1}.
The functional equation is established in
\cite[Theorem 4.1]{Sh1}  combined with \cite[Prop. 3.1]{Sh3} and 
\cite[Theorem 1]{Sh3}. See also \cite{Sh5} for the general case. The location 
of the poles has been determined in the appendix of \cite{MW}. Property iii)
was proved in \cite[p.280]{RS}.

Now let $\pi_1\in\Pi_{\di}(\GL_{m_1}(\A_E))$ and 
$\pi_2\in\Pi_{\di}(\GL_{m_2}(\A_E))$. Using the 
description of the residual spectrum for $\GL_n$  
 \cite{MW}, $L(s,\pi_1\times\pi_2)$ can be expressed  in terms of
Rankin-Selberg $L$-functions attached to cuspidal automorphic
representations. Indeed, by \cite{MW}
 there exist $k_i\in\N$ with $k_i|m_i$, parabolic
subgroups $P_i$ of $G_i=\GL_{m_i}$ of type $(d_i,...,d_i)$, $d_i=m_i/k_i$, and
unitary cuspidal automorphic representations $\delta_i$ of $\GL_{d_i}(\A_E)$ 
such that 
\begin{equation}\label{3.16}
\pi_i=J^{G_i}_{P_i}(\delta_i[(k_i-1)/2]\otimes\cdots\otimes\delta_i[(1-k_i)/2]),
\end{equation}
where  the right hand side denotes the unique irreducible quotient of the
 induced representation
$I^{G_i}_{P_i}(\delta_i[(k_i-1)/2]\otimes\cdots\otimes\delta_i[(1-k_i)/2])$.
Set $k=k_1+k_2-2$. Then it follows from (\ref{4.4n}) that
\begin{equation}\label{3.17}
L(s,\pi_1\times\pi_2)=\prod_{i=0}^{k_1-1}\prod_{j=0}^{k_2-1}
L(s+k/2-i-j,\delta_1\times\delta_2).
\end{equation}
Using this equality and i)--iv) above,  we deduce immediately the corresponding
properties satisfied by $L(s,\pi_1\times\pi_2)$. Especially, 
$L(s,\pi_1\times\pi_2)$ satisfies a functional equation of the form 
(\ref{3.14}) with an $\epsilon$-factor similar to (\ref{3.15}).

We shall now investigate the logarithmic derivatives of the Rankin-Selberg
$L$-functions. First we
need to introduce some notation. Let $\pi_i\in\Pi_{\di}(\GL_{m_i}(\A_E))$,
 $i=1,2$. For each Archimedean place 
$w$ of $E$ let $\mu_{j,k}(\pi_{1,w} \times\pi_{2,w})$ $j=1,\ldots,r_w,$ 
$k=1,\ldots,h_w,$
be the parameters attached to $(\pi_{1,w},\pi_{2,w})$ by means of (\ref{3.12}).
Set 
\begin{equation}\label{3.18}
c(\pi_1\times\pi_2)=\sum_{w|\infty}\sum_{j,k}|\mu_{j,k}(\pi_{1,w}\times
\pi_{2,w})|.
\end{equation}
Let $N(\pi_1\times\pi_2)$ be the integer that is determined by the 
$\epsilon$-factor as in (\ref{3.15}).
 Set 
\begin{equation}\label{3.19}
\nu(\pi_1\times\pi_2)=D_E^{m_1m_2}N(\pi_1\times\pi_2)(2+c(\pi_1\times\pi_2)).
\end{equation}
We call $\nu(\pi_1\times\pi_2)$ the level of $(\pi_1,\pi_2)$. Given 
$\pi\in\Pi(\GL_m(\A_E))$, set
$$\pi_\infty=\otimes_{v|\infty}\pi_v,\quad \pi_f=\otimes_{v<\infty}\pi_v.$$

\begin{lem}\label{l3.3}
For every $\epsilon>0$ there exists $C>0$ such that
$$\bigg |\frac{L^\prime(s,\pi_{1,f}\times\pi_{2,f})}{L(s,\pi_{1,f}\times
\pi_{2,f})} \bigg |\le C$$
for all $s$ in the half-plane $\Re(s)\ge2+\epsilon$ and all $\pi_i\in\Pi_{\cu}
(\GL_{m_i}(\A_E)),$ $i=1,2.$
\end{lem}
\begin{proof}
Let $\pi_i\in\Pi_{\cu}(\GL_{m_i}(\A_E))$, $i=1,2$, and let $v<\infty$. By 
\cite{Sk},
 $\pi_{1,v}$ and $\pi_{2,v}$ are unitary generic
representations. Hence by Lemma \ref{l3.1} there exist complex numbers
$a_{i}(v),$ $i=1,\ldots,m_1\cdot m_2,$ with
\begin{equation}\label{3.20}
|a_{i}(v)|<N(v)
\end{equation}
such that
$$L(s, \pi_{1,v}\times\pi_{2,v})=
\prod_{i=1}^{m_1\cdot m_2}(1-a_{i}(v)N(v)^{-s})^{-1}.$$
Suppose that $\Re(s)>1.$ By (\ref{3.20}) we have $|a_{i}(v)/N(v)^s|<1.$
Hence, taking the logarithmic derivative, we get
\begin{equation*}
\begin{split}
\frac{L^\prime(s,\pi_{1,v}\times\pi_{2,v})}{L(s,\pi_{1,v}\times\pi_{2,v})}
& = -\sum_{i}\frac{a_{i}(v)\log N(v)}{N(v)^s(1-a_{i}(v)N(v)^{-s})}\\
& = -\log N(v)\sum_{i}\sum^\infty_{k=1}\frac{a_{i}(v)^k}{N(v)^{sk}}.
\end{split}
\end{equation*}
Suppose that $\sigma=\Re(s)>1$. Then by (\ref{3.20}) we get
$$\bigg | \frac{L^\prime(s,\pi_{1,v}\times\pi_{2,v})}{L(s,\pi_{1,v}\times
\pi_{2,v})}\bigg |\le m_1m_2\sum^\infty_{k=1}\frac{\log N(v)}{N(v)^
{(\sigma-1)k}}.$$
Let $\zeta_E(s)$ be the Dedekind zeta function of $E$. Let $\epsilon>0$ and
set $\sigma=2+\epsilon$. Then for $\Re(s)\ge\sigma$ we get
$$\bigg |\frac{L^\prime(s,\pi_{1,f}\times\pi_{2,f})}{L(s,\pi_{1,f}\times
\pi_{2,f})}\bigg |\le m_1m_2  \bigg |
\frac{\zeta^\prime_E(\sigma-1)}{\zeta_E(\sigma-1)}\bigg |.$$
\end{proof}

\begin{lem}\label{l3.4} For every $\epsilon>0$ there exists $C>0$ such that
$$\bigg |\frac{L^\prime(s,\pi_{1,\infty}\times\pi_{2,\infty})}
{L(s,\pi_{1,\infty}\times\pi_{2,\infty})}\bigg |\le C(1+\log(|s|+c(\pi_1\times
\pi_2)))$$
for all $s$ with $\Re(s)\ge 1+\epsilon$ and all 
$\pi_i\in\Pi_{\cu}(\GL_{m_i}(\A_E))$, $i=1,2$.
\end{lem}
\begin{proof}
Let $w|\infty.$ By (\ref{3.12}) we have
\begin{equation}\label{3.21}
L(s,\pi_{1,w}\times\pi_{2,w})=\prod_{j,k}\Gamma_\R(s+\mu_{j,k}(\pi_{1,w}
\times\pi_{2,w})).
\end{equation}
Since $\pi_{1,w}$ and $\pi_{2,w}$ are unitary and generic, the complex 
numbers $\mu_{j,k}(\pi_{1,w}\times\pi_{2,w})$ satisfy
\begin{equation}\label{3.22}
\Re(\mu_{j,k}(\pi_{1,w}\times\pi_{2,w}))>-1.
\end{equation}
Now recall that by Stirlings formula
$$\frac{\Gamma^\prime}{\Gamma}(s)=\log s+O(|s|^{-1})$$
is valid as $|s|\to\infty,$ in the angle $-\pi+\delta<\arg s<\pi-\delta,$ for 
any fixed $\delta>0.$ Hence
\begin{equation}\label{3.23}
\frac{\Gamma^\prime_{\R}(s)}{\Gamma_{\R}(s)} (s)=-\frac{1}{2}\log\pi+\log
s+O(|s|^{-1})
\end{equation}
holds in the same range of $s.$ Let $\epsilon>0.$ Using 
(\ref{3.21}), (\ref{3.22}) and 
(\ref{3.23}), it follows that there exists $C>0$ such that
$$\bigg|\frac{L^\prime(s,\pi_{1,w}\times\pi_{2,w})}{L(s,\pi_{1,w}\times 
\pi_{2,w})}\bigg|
\le C+\sum_{j,k}\log(|s|+|\mu_{j,k}(\pi_{1,w}\times\pi_{2,w})|).$$
for all $w|\infty$, all $s$ with $\Re(s)\ge 1+\epsilon$ and all 
$\pi_i\in\Pi_{\cu}(\GL_m(\A_E)),$ $i=1,2.$ This implies
the lemma.
\end{proof}
Let $\pi_i\in\Pi_{\di}(\GL_{m_i}(\A_E))$, $i=1,2$, and $T>0$, be given. Denote
by
$N(T;\pi_1,\pi_2)$ the number of zeros of $L(s,\pi_1\times\pi_2),$ counted
with multiplicity, which are contained in the disc of radius $T$ centered at 0.

\begin{prop}\label{p3.5}
There exists $C>0$ such that
$$N(T;\pi_1,\pi_2)\le C T\log(T+\nu(\pi_1\times\pi_2))$$ for all
$T\ge1$ and all $\pi_i\in\Pi_{\di}(\GL_m(\A_E)),$ $i=1,2.$
\end{prop}
\begin{proof} By (\ref{3.17}) we can assume that $\pi_1$ and $\pi_2$ are
unitary cuspidal automorphic representations. Set
\begin{equation}\label{3.24}
\Lambda(s)=s^a(1-s)^a\left(D_E^{m_1m_2}N(\pi_1\times\pi_2)\right)^{s/2}L(s,\pi_1\times\pi_2),
\end{equation}
where $a$ denotes the order of the pole of $L(s,\pi_1\times\pi_2)$ at 
$s=1.$ Note that $a$ can be at most $1.$ Since $\pi_i$ is unitary, we have
$\widetilde\pi_i=\overline{\pi_i},$ $i=1,2.$ Hence by (\ref{3.14}), it follows
that $\Lambda(s)$ satisfies the functional equation
\begin{equation}\label{3.25}
\Lambda(s)=W(\pi_1\times\pi_2)\left(D_E^{m_1m_2}N(\pi_1\times\pi_2)\right)^{1/2}\overline{\Lambda(1-\overline s)}.
\end{equation}
Since $L(s,\pi_1\times\pi_2)$ is of order one, $\Lambda(s)$ is an entire
function of order one and hence, it admits a representation as a Weierstra\ss{}
product of the form
$$\Lambda(s)=e^{A+Bs}\prod_\rho(1-s/\rho)e^{s/\rho},$$
where $A,B\in\C$ and the product runs over the set of zeros of $\Lambda(s).$
We note that for $s=\sigma+i T$
\begin{equation}\label{3.26}
\Re\sum_\rho\frac{1}{s-\rho}=\sum_{\rho}\frac{\sigma-\be}{
(\sigma-\be)^2+(\gam-T)^2}
\end{equation}
and this series is convergent since $\Lambda(s)$ is of oder one. Taking the 
real part of the logarithmic derivative of $\Lambda(s),$ and applying the
functional equation
(\ref{3.25}) to the right hand side, we get
\begin{equation*}
\begin{split}
\Re(B)+\Re\sum\frac{1}{\rho}+\Re\sum_\rho\frac{1}{s-\rho}=&-\Re(\overline B)
-\Re\sum\frac{1}{\rho}\\
&+\Re\sum_\rho\frac{1}{s-(1-\overline\rho)}.
\end{split}
\end{equation*}
Now observe that by 
(\ref{3.25}), $\rho$ is a zero of $\Lambda(s)$ if and only 
if $1-\overline\rho$ is a zero of $\Lambda(s).$ Hence the two sums involving
$s$ are equal, as they run over the same set of zeros. It follows that
\begin{equation}\label{3.27}
\Re(B)=-\Re(\sum_\rho\frac{1}{\rho}).
\end{equation}
Together with (\ref{3.24}) this leads to
\begin{equation*}
\begin{split}
\Re\sum_\rho\frac{1}{s-\rho} =\Re\frac{\Lambda^\prime(s)}{\Lambda(s)}&=
\frac{a}{s}+\frac{a}{s-1}+\frac{1}{2}\log (D_E^{m_1m_2}N(\pi_1\times \pi_2)) \\
& +\frac{L^\prime(s,\pi_{1,\infty}\times\pi_{2,\infty})}
{L(s,\pi_{1,\infty}\times\pi_{2,\infty})}
+\frac{L^\prime(s,\pi_{1,f}\times\pi_{2,f})}{L(s,\pi_{1,f}
\times\pi_{2,f})}.
\end{split}
\end{equation*}
Let $\epsilon>0,$ and set $\sigma=2+\epsilon.$ 
By Lemma \ref{l3.3}, Lemma \ref{l3.4} and (\ref{3.26})
it follows that there exists $C>0$ such that
\begin{equation}\label{3.28}
\begin{split}
\sum_{\rho}\frac{\sigma-\be}{(\sigma-\be)^2+(\gam-T)^2}
&\le \frac{1}{2}\log(D_E^{m_1m_2} N(\pi_1\times\pi_2))\\
&\quad+C(1+\log(|T|+c(\pi_1\times\pi_2))\\
&\le C_1\log(|T|+\nu(\pi_1\times\pi_2))
\end{split}
\end{equation}
for all $T\in\R$ and $\pi_i\in\Pi_{\cu}(\GL_{m_i}(\A_E)),$ $i=1,2.$ Let $T>0.$
Then it follows from (\ref{3.28}) that
\begin{equation*}
\begin{split}
N(T+1;\pi_1,\pi_2)-&N(T;\pi_1,\pi_2)\\
&\le 2(3+\epsilon)
\sum_{\rho} \frac{\sigma-\be}{(\sigma-\be)^2+(\gam-T)^2}\\
&  \le
C\log(T+\nu(\pi_1\times\pi_2))
\end{split}
\end{equation*}
for all $\pi_i\in\Pi_{\cu}(\GL_{m_i}(\A_E)),$ $i=1,2.$ This implies the 
proposition.
\end{proof}

\section{Normalizing factors}
\setcounter{equation}{0}

In this section we consider the global normalizing factors of  
intertwining operators. Our main purpose is to estimate certain integrals 
involving the logarithmic derivatives of the  normalizing factors. 
The behaviour of these integrals is crucial for the estimation of the 
spectral side. From now on we assume that the ground field is $\Q$.
 Denote by $\A$ the ring of ad\`eles of $\Q$.

Let $M\in{\L}$. Then there exists a partition $(n_1,\ldots,n_r)$ of 
$n$ such that
$$M=\GL_{n_1}\times\cdots\times\GL_{n_r}.$$
Let $Q,P\in\cP(M).$ Let $v$ be a place of $\Q$ and let $\pi_v\in\Pi(M(\Q_v)).$
Associated to $P,Q$ and $\pi_v$ is the local intertwining operator
$$J_{Q|P}(\pi_v,\lambda),\quad \lambda\in\af_{M,\C}^*,$$  
between the induced representations $I_P(\pi_{v,\lambda})$ and 
$I_Q(\pi_{v,\lambda})$, which is defined by an integral over $N_Q(\Q_v)\cap
N_{\ov P}(\Q_v)$, and hence depends upon a choice of Haar measure on this 
group. By \cite{A7} there exist meromorphic 
 functions $r_{Q|P}(\pi_v,\lambda)$, $\lambda\in\af^*_{M,\C}$, such that the
normalized local intertwining operators
$$R_{Q|P}(\pi_v,\lambda)=r_{Q|P}(\pi_v,\lambda)^{-1}J_{Q|P}(\pi_v,\lambda)$$
satisfy the conditions of Theorem 2.1 of \cite{A7}. 
There is a general construction of normalizing factors which works for any
 group \cite{A7}, \cite{CLL}. For $\GL_n$, however, the intertwining operators
can be normalized by $L$-functions \cite[\S4]{A7}, \cite[p.87]{AC}. The
normalizing factors are given as
\begin{equation}\label{4.1}
r_{Q|P}(\pi_v,\lambda)=\prod_{\alpha\in\sum_P\cap\sum_{\overline Q}}
r_\alpha(\pi_v,\lambda(\check{\alpha})),
\end{equation}

where $r_\alpha(\pi_v,s)$ is a meromorphic function of one complex variable
and $\Sigma_P$ (resp. $\Sigma_{\ov Q}$)  denotes the roots of $(P,A_M)$ (resp.
$(\ov Q,A_M)$). 
Thus to define the normalizing factors, it is enough to define the functions
$r_\alpha(\pi_v,s)$ for any root $\alpha$ of $(G,A_M)$ and any $\pi_v\in
\Pi(M(\Q_v))$. To this end note that $\pi_v$ is equivalent to a representation
$\pi_{1,v}\otimes\cdots\otimes \pi_{r,v}$ with 
$\pi_{i,v}\in\Pi(\GL_{n_i}(\Q_v))$ and the root $\alpha$ 
corresponds to an ordered pair $(i,j)$ of distinct integers between 1 and $r$.
Fix a nontrivial additive character $\psi_v$ of $\Q_v$. Let 
$L(s,\pi_{i,v}\times\widetilde \pi_{j,v})$ and 
$\epsilon(s,\pi_{i,v}\times\widetilde\pi_{j,v},\psi_v)$ be the Rankin-Selberg
$L$-function  and the $\epsilon$-factor attached to $(\pi_{i,v}, 
\widetilde\pi_{j,v})$ and $\psi_v$. 
Set
\begin{equation}\label{4.2}
r_\alpha(\pi_v,s)=\frac{L(s,\pi_{i,v}\times\widetilde \pi_{j,v})}
{L(1+s,\pi_{i,v}\times\widetilde \pi_{j,v})\epsilon(s,\pi_{i,v}\times\widetilde
\pi_{j,v},\psi_v)}.
\end{equation}
It follows from Theorem 6.1 of \cite{Sh1} that there are Haar measures on the 
group $N_Q(\Q_v)\cap N_{\overline P}(\Q_v)$, depending on $\psi_v$, such that
the normalizing factors (\ref{4.1}) have all the right properties
(see \cite[\S 4]{A7}, \cite[p.87]{AC}). 
Now suppose that $\pi\in\Pi_{\di}(M(\A)).$ Then the global normalizing factor
$r_{Q|P}(\pi,\lambda)$ is defined by the infinite product
$$r_{Q|P}(\pi,\lambda)=\Pi_v r_{Q|P}(\pi_v,\lambda),$$
which converges in a certain chamber. By (\ref{4.1}) it follows that there
exist meromorphic functions $r_\alpha(\pi,s)$ of one complex variable such 
that
\begin{equation}\label{4.3}
r_{Q|P}(\pi,\lambda)=\prod_{\alpha\in\sum_P\cap\sum_{\overline Q}}
r_\alpha(\pi,\lambda(\check{\alpha})).
\end{equation}

Let $\pi=\pi_1\otimes\cdots\otimes\pi_r.$ If $\alpha$ corresponds to $(i,j)$
then by (\ref{4.2}) we have
\begin{equation}\label{4.4}
r_\alpha(\pi,s)=\frac{L(s,\pi_i\times\widetilde\pi_j)}
{L(1+s,\pi_i\times\widetilde\pi_j)\epsilon(s,\pi_i\times\widetilde\pi_j)},
\end{equation}
where $L(s,\pi_i\times\widetilde\pi_j)$ and 
$\epsilon(s,\pi_i\times\widetilde\pi_j)$ are the global $L$-function and
the $\epsilon$-factor, respectively, considered in the previous section.

The main goal of this section is to study the multidimensional logarithmic 
derivatives of the normalizing factors that occur on the spectral side of
the trace formula \cite{A4}. By (\ref{4.3}) this problem is reduced
to the investigation of the logarithmic derivatives of the analytic
functions $r_\alpha(\pi,s).$ Furthermore, by (\ref{4.4}) each 
$r_\alpha(\pi,s)$ may be regarded as the normalizing factor attached to a 
standard maximal parabolic subgroup in $\GL_m$ with $m\le n.$
So let $m_1,m_2\in\N$ with $m_1+m_2\le n.$ Given $\pi_i \in\Pi_{\di}(\GL_
{m_i}(\A)),$ $i=1,2,$ set

\begin{equation}\label{4.5}
r(\pi_1\otimes\pi_2,s)=\frac{L(s,\pi_1\times\widetilde\pi_2)}{L(1+s,\pi_1\times
\widetilde\pi_2)\epsilon(s,\pi_1\times\widetilde\pi_2)}.
\end{equation}

We shall now study the logarithmic derivatives of these functions.
 For this purpose we need some preparation. Suppose that $\pi_i,$
$i=1,2,$ is given in the form (\ref{3.17}) and assume that $k_1\le k_2.$ 
Set $k=k_1+k_2-2.$ For $j=0,\ldots,k$ let the integers $a_j$ be defined by
\begin{equation}\label{4.6}
a_i=\left\{\begin{array}{r@{\quad:\quad}l}
i+1& i\le k_1-1;\\ k_1& k_1-1\le i \le k_2-1;\\ k-i+1& i\ge k_2-1.
\end{array}\right.
\end{equation}
Note that $a_i=a_{k-i}$, $i=0,...,k$. 
It follows from (\ref{3.17}) that
\begin{equation}\label{4.7}
L(s,\pi_1\times\widetilde\pi_2)=\prod_{i=0}^kL(s+k/2-i,\delta_1\times\widetilde\delta_2)^{a_i}.
\end{equation}
Define a polynomial of one variable $x$ by
$$p(x)=\prod_{i=0}^k\left((x+k/2-i)(1-x-k/2+i)\right)^{a_i}.$$
Then $p(x)$ has real coeefficients and satisfies $p(x)=p(1-x)$. Let $a$ be the
order of the pole of $L(s,\delta_1\times\widetilde\delta_2)$ at $s=1$. Note 
that $a\le 1$. Set
\begin{equation}\label{4.8}
\Lambda(s)=p(s)^a N(\pi_1\times\widetilde\pi_2)^{s/2}
L(s,\pi_1\times\widetilde\pi_2).
\end{equation}
Then $\Lambda(s)$ satisfies the functional equation
\begin{equation}\label{4.9}
\Lambda(s)=W(\pi_1\times\widetilde\pi_2)N(\pi_1\times\widetilde\pi_2)^{1/2}\overline{\Lambda(1-\overline{s})}.
\end{equation}
Furthermore $\Lambda(s)$ is an entire function of order 1. 
Therefore it can be written as Weierstrass product of the form
$$\Lambda(s)=e^{A+Bs}\prod_\rho(1-s/\rho)e^{s/\rho}$$
with $A,B\in\C$ and $\rho$ runs over the zeros of $\Lambda(s)$. Taking the
logarithmic derivative and applying the functional equation (\ref{4.9}) to the
right hand side, we get
\begin{equation*}
\begin{split}
\left(\frac{\Lambda(s)}{\Lambda(s+1)}\right)^\prime&\cdot
\frac{\Lambda(s+1)}{\Lambda(s)}=
\frac{\Lambda^\prime(s)}{\Lambda(s)}
+{\frac{\ov{\Lambda^\prime(-\overline{s})}}{\ov{\Lambda(-\overline{s})}}}\\
&=2\Re(B)+2\Re\sum_\rho\frac{1}{\rho}
+\sum_\rho\left\{\frac{1}{s-\rho}-
\frac{1}{s+\overline{\rho}}\right\}.
\end{split}
\end{equation*}
By (\ref{3.27}) it follows that the first two terms on the right hand side 
cancel and hence we get
\begin{equation}\label{4.9a}
\left(\frac{\Lambda(s)}{\Lambda(s+1)}\right)^\prime\cdot
\frac{\Lambda(s+1)}{\Lambda(s)}=
2\sum_\rho\frac{\Re(\rho)}{(s-\rho)(s+\overline{\rho})}.
\end{equation}
Therefore, combining (\ref{3.15}), (\ref{4.5}) and (\ref{4.8}), we obtain
\begin{equation*}
\begin{split}
\frac{r^\prime(\pi_1\otimes\pi_2,s)}{r(\pi_1\otimes\pi_2,s)}=&\log N(\pi_1\times\widetilde
\pi_2)\\
&+a\left\{\frac{p^\prime(s+1)}{p(s+1)}-\frac{p^\prime(s)}{p(s)}\right\}
+2\sum_\rho\frac{\Re(\rho)}{(s-\rho)(s+\overline{\rho})}.
\end{split}
\end{equation*}
In particular, if $s=i\lambda$, $\lambda\in\R$, then it follows from the
definition of $p(s)$ that
\begin{equation*}
\begin{split}
\frac{r^\prime(\pi_1\otimes\pi_2,i\lambda)}{r(\pi_1\otimes\pi_2,i\lambda)}&=
\log N(\pi_1\times\widetilde
\pi_2)\\
&+2a\sum_{i=0}^k\left\{\frac{a_i(k/2-i+1)}{\lambda^2+(k/2-i+1)^2}-
\frac{a_i(k/2-i-1)}{\lambda^2+(k/2-i-1)^2}\right\}\\
&+2\sum_{\rho}\frac{\be}{\be^2+(\gam-\lambda)^2}.
\end{split}
\end{equation*}

\begin{prop}\label{p4.1}
There exists $C>0$ such that
$$\int_{-T}^T\bigg|\frac{r^\prime(\pi_1\otimes\pi_2,i\lambda)}
{r(\pi_1\otimes\pi_2,i\lambda)}\bigg|
\;d\lambda
\le CT\log(T+\nu(\pi_1\times\widetilde\pi_2))$$
for all $T>0$ and $\pi_i\in\Pi_{\di}(\GL_{m_i}(\A))$, $i=1,2$.
\end{prop}
\begin{proof}
By the above formula it suffices to estimate the integral
$$\int_{-T}^T\sum_{\rho}\frac{|\be|}{\be^2+(\gam-\lambda)^2}\,d\lambda.$$
We split the series as follows
$$\sum_\rho=\sum_{|\Im(\rho)|\le T+1}+\sum_{|\Im(\rho)|>T+1}.$$
To estimate the integral of the first sum, observe that for all 
$\beta\in\R^+$ and $\gamma\in\R$ we have
$$\int_{-T}^T\frac{\beta}{\beta^2+(\gamma-\lambda)^2}\,d\lambda\le
\int_{-\infty}^\infty\frac{d\lambda}{1+\lambda^2}=\pi.$$
Hence by Proposition \ref{p3.5} we get
\begin{equation*}
\begin{split}
\int_{-T}^T\sum_{|\Im(\rho)|\le T+1}\frac{|\be|}{\be^2+(\gam-\lambda)^2}\,d\lambda&\le \pi N(T+1,\pi_1,\widetilde\pi_2)\\
&\le CT\log(T+\nu(\pi_1\times\widetilde\pi_2)).
\end{split}
\end{equation*}

It remains to consider the integral of the second sum. Observe that by
(\ref{4.7}) the zeros $\rho$ of $\Lambda(s)$ satisfiy $|\Re(\rho)|\le k/2+1$.
Set
$$\sigma=k+3,\quad C=2(k+2)^2.$$
Then the following inequality holds for all $\lambda\in\R$ with 
$|\lambda|\le T$, all $\beta\in\R^\times$ with $|\beta|\le k/2+1$ and all
$\gamma\in\R$ with $|\gamma|>T+1$:
$$\frac{|\beta|}{\beta^2+(\gamma-\lambda)^2}\le C\left\{
\frac{\sigma-\beta}{(\sigma-\beta)^2+(\gamma-T)^2}+
\frac{\sigma-\beta}{(\sigma-\beta)^2+(\gamma+T)^2}\right\}.$$
Thus we get
\begin{equation*}
\begin{split}
\sum_{|\Im(\rho)|>T+1}&\frac{|\be|}{\be^2+(\gam-\lambda)^2}\\
&\hskip40pt\le C\Biggl\{\sum_{\rho}\frac{\sigma-\be}{(\sigma-\be)^2
+(\gam-T)^2}\\
&\hskip60pt+\sum_{\rho}\frac{\sigma-\be}{(\sigma-\be)^2+(\gam+T)^2}\Biggr\}.
\end{split}
\end{equation*}
Combining (\ref{4.7}) and (\ref{3.28}), it follows that for $\sigma=k+3$ 
there exists $C_1>0$ such that
$$\sum_{\rho}\frac{\sigma-\be}{(\sigma-\be)^2+(\gam-T)^2}
\le C_1\log(|T|+\nu(\pi_1\times\widetilde\pi_2))$$
for all $T\in\R$ and $\pi_i\in\Pi_{\di}(\GL_{m_i}(\A))$, $i=1,2$. Combining
these observation we
 get
$$\int_{-T}^T\sum_{|\Im(\rho)|>T+1}\frac{|\be|}{\be^2+(\gam-\lambda)^2}
\le CT\log(T+\nu(\pi_1\times\widetilde\pi_2)).$$
This completes the proof of the proposition.
\end{proof}

The next proposition will be important for the determination of the asymptotic
behaviour of the spectral side.

\begin{prop}\label{p4.2}
There exists $C>0$ such that
\begin{equation*}
\begin{split}
\int_{-\infty}^\infty|r^\prime(\pi_1\otimes\pi_2,i\lambda)
r(\pi_1\otimes\pi_2,&i\lambda)^{-1}|e^{-t\lambda^2}\;d\lambda\\
&\le C\log(1+\nu(\pi_1\times\widetilde\pi_2))
\frac{1+|\log t|}{\sqrt{t}}
\end{split}
\end{equation*}
for all $0<t\le 1$ and $\pi_i\in\Pi_{\di}(\GL_{m_i}(\A))$, $i=1,2$.
\end{prop}
\begin{proof}
By Proposition \ref{p4.1} it follows that we have
$$\int_0^\lambda|r^\prime(\pi_1\otimes\pi_2,iu)
r(\pi_1\otimes\pi_2,iu)^{-1}|du\le 
C\lambda^2$$
as $|\lambda|\to\infty$. Hence, using integration by parts, it follows that
the integral on the left hand side of the claimed inequality equals
$$2t\int_{-\infty}^\infty\int_0^\lambda|r^\prime(\pi_1\otimes\pi_2,iu)
r(\pi_1\otimes\pi_2,iu)^{-1}|du\;\lambda e^{-t\lambda^2}\;d\lambda.$$
Applying Proposition \ref{p4.1} we get
\begin{equation*}
\begin{split}
\int_{-\infty}^\infty|r^\prime(\pi_1\otimes\pi_2,i\lambda)
r(\pi_1&\otimes\pi_2,i\lambda)^{-1}|e^{-t\lambda^2}\;d\lambda\\
&\le Ct\int_{-\infty}^\infty \log\left(|\lambda|+\nu(\pi_1\times
\widetilde\pi_2)\right)\lambda^2e^{-t\lambda^2}d\lambda\\
&\le C_1\log(1+\nu(\pi_1\times\widetilde\pi_2))\frac{1+|\log t|}{\sqrt{t}}
\end{split}
\end{equation*}
for all $0<t\le 1$ and $\pi_i\in\Pi_{\di}(\GL_{m_i}(\A))$, $i=1,2$.
\end{proof}

Let $M\in\L$ and let $Q,P\in\cP(M)$. 
Our next goal is to estimate the corresponding integrals involving the
 generalized logarithmic derivatives of the
global normalizing factors $r_{Q|P}(\pi,\lambda)$.
For this purpose we will 
use the notion of a $(G,M)$ family introduced by Arthur in Section 6 of 
\cite{A5}. For the convenience of the reader we  recall the  definition
of a $(G,M)$ family and explain some of its properties.

For each $P\in\cP(M)$, let $c_P(\lambda)$ be a smooth function on $i\af_M^*$. 
Then the set 
$$\{c_P(\lambda)\mid P\in\cP(M)\}$$
is called a $(G,M)$ family if the following holds: Let $P,P^\prime\in\cP(M)$
be adjacent parabolic groups and suppose that $\lambda$ belongs to the 
hyperplane spanned by the common wall of the chambers of $P$ and $P^\prime$.
Then
$$c_P(\lambda)=c_{P^\prime}(\lambda).$$
Let
\begin{equation}\label{4.10}
\theta_P(\lambda)=\vol\left(\af_P^G/\Z(\Delta_P^\vee)\right)^{-1}
\prod_{\alpha\in\Delta_P}\lambda(\alpha^\vee),\quad\lambda\in i\af_P^*,
\end{equation}
where $\Z(\Delta_P^\vee)$ is the lattice in $\af_P^G$ generated by the 
co-roots
$$\{\alpha^\vee\mid \alpha\in\Delta_P\}.$$
Let $\{c_P(\lambda)\}$ be a $(G,M)$ family. Then by Lemma 6.2 of \cite{A5},
the function
\begin{equation}\label{4.11}
c_M(\lambda)=\sum_{P\in\cP(M)}c_P(\lambda)\theta_P(\lambda)^{-1}
\end{equation}
extends to a smooth function on $i\af_M^*$. The value of $c_M(\lambda)$ at
$\lambda=0$ is of particular importance in connection with the spectral side
of the trace formula. It can be computed as follows. Let $p=\dim(A_M/A_G)$.
Set $\lambda=t\Lambda$, $t\in\R$, $\Lambda\in\af_M^*$, and let $t$ tend to 0.
Then
\begin{equation}\label{4.12}
c_M(0)=\frac{1}{p!}\sum_{P\in\cP(M)}\left(\lim_{t\to0}\left(\frac{d}{dt}\right)^p c_P(t\Lambda)\right)\theta_P(\Lambda)^{-1}
\end{equation}
\cite[(6.5)]{A5}. This expression is of course  independent of $\Lambda$.

For any $(G,M)$ family $\{c_P(\lambda)\mid P\in\cP(M)\}$ and any $L\in\cL(M)$
there is associated a natural $(G,L)$ family which is defined as follows.
 Let $Q\in\cP(L)$ and suppose that $P\subset Q$. 
The  function
$$\lambda\in i\af_L^*\mapsto c_P(\lambda)$$
depends only on $Q$. We will  denote it by $c_Q(\lambda)$. Then
$$\{c_Q(\lambda)\mid Q\in\cP(L)\}$$
is a $(G,L)$ family. We write
$$c_L(\lambda)=\sum_{Q\in\cP(L)}c_Q(\lambda)\theta_Q(\lambda)^{-1}$$
for the corresponding function (\ref{4.11}). 

Let $Q\in\cP(L)$ be fixed. If $R\in\cP^L(M)$, then 
$Q(R)$ is the unique group in $\cP(M)$ such that $Q(R)\subset Q$ and
$Q(R)\cap L=R$. Let $c_R^Q$ be the function on $i\af_M^*$ which is defined by
$$c_R^Q(\lambda)=c_{Q(R)}(\lambda).$$
Then $\{c_R^Q(\lambda)\mid R\in\cP^L(M)\}$ is an $(L,M)$ family. Let 
$c_M^Q(\lambda)$ be the function (\ref{4.11}) associated to this $(L,M)$ 
family.

We consider now special $(G,M)$ families defined by the global 
normalizing factors. Fix $P\in\cP(M)$, $\pi\in\Pi_{\di}(M(\A))$ and 
$\lambda\in i\af_M^*$. Define
\begin{equation}\label{4.13}
\nu_Q(P,\pi,\lambda,\Lambda):=r_{Q|P}(\pi,\lambda)^{-1}
r_{Q|P}(\pi,\lambda+\Lambda),\quad Q\in\cP(M).
\end{equation}
This set of functions is a $(G,M)$ family \cite[p.1317]{A4}. 
It is of a special form. 
By (\ref{4.3}) we have
$$\nu_Q(P,\pi,\lambda,\Lambda)=\prod_{\alpha\in\Sigma_Q\cap\Sigma_{\ov P}}
 r_\alpha(\pi,\lambda(\alpha^\vee))^{-1}
 r_\alpha(\pi,\lambda(\alpha^\vee)+\Lambda(\alpha^\vee)).$$
Suppose that $L\in\cL(M)$, $L_1\in\cL(L)$ and $S\in\cP(L_1)$. Let 
$$\{\nu_{Q_1}^S(P,\pi,\lambda,\Lambda)\mid Q_1\in\cP^{L_1}(L)\}$$ 
be the associated $(L_1,L)$ family and let $\nu_L^S(P,\pi,\lambda,\Lambda)$
be the function (\ref{4.11}) defined by this family. Set
$$\nu_L^S(P,\pi,\lambda):=\nu_L^S(P,\pi,\lambda,0).$$
If $\alpha$ is any root in $\Sigma(G,A_M)$, let $\alpha_L^\vee$ denote the 
projection of $\alpha^\vee$ onto $\af_L$. If $F$ is a subset of 
$\Sigma(G,A_M)$, let $F_L^\vee$ be the disjoint union of all the vectors
$\alpha_L^\vee$, $\alpha\in F$. 
Then by Proposition 7.5 of \cite{A4} we have
\begin{equation}\label{4.14}
\begin{split}
\nu_L^S(P,\pi,\lambda)=&\sum_F\vol(\af_L^{L_1}/\Z(F_L^\vee))\\
&\cdot\left(\prod_{\alpha\in F} r_\alpha(\pi,\lambda(\alpha^\vee))^{-1}
r_\alpha^\prime(\pi,\lambda(\alpha^\vee))\right),
\end{split}
\end{equation}
where $F$ runs over all subsets of $\Sigma(L_1,A_M)$ such that $F_L^\vee$ 
is a basis of $\af_L^{L_1}$. Let $t>0$. Then by (\ref{4.14}) we get
\begin{equation*}
\begin{split}
\int_{i\af_L^*/\af_G^*}|&\nu_L^S(P,\pi,\lambda)|
e^{-t\parallel\lambda\parallel^2}\,d\lambda\le
\sum_F\vol(\af_L^{L_1}/\Z(F_L^\vee))\\
&\cdot\int_{i\af_L^*/\af_G^*}\prod_{\alpha\in F}\Big| r_\alpha\pi,\lambda(\alpha\vee))^{-1}
r_\alpha\prime(\pi,\lambda(\alpha\vee))\Big|e^{-t\parallel\lambda\parallel^2}
\;d\lambda.
\end{split}
\end{equation*}
Fix any subset $F$ of $\Sigma(L_1,A_M)$ such that $F_L^\vee$ is a basis of 
$\af_L^{L_1}$.  Let 
$$\{\tilde\omega_\alpha\mid \alpha\in F\}$$
be the basis of $(\af_L^{L_1})^*$ which is dual to $F_L^\vee$. We
 can write $\lambda\in i\af_L^*/i\af_G^*$ as
$$\lambda=\sum_{\alpha\in F}z_\alpha\tilde\omega_\alpha+\lambda_1,\quad z_\alpha\in
i\R,\quad\lambda_1\in i\af_{L_1}^*/i\af_G^*.$$
Observe that $\lambda(\alpha^\vee)=z_\alpha$. Let $l_1=\dim(A_{L_1}/A_G)$.
 Then there exists $C>0$, independent of $\pi$,  such that for all $t>0$ we
 have
\begin{equation}\label{4.15}
\begin{split}
\int_{i\af_L^*/\af_G^*}\prod_{\alpha\in F}&\big| r_\alpha(\pi,\lambda(\alpha^\vee))^{-1}
r_\alpha^\prime(\pi,\lambda(\alpha^\vee))\big|e^{-t\parallel\lambda\parallel^2}\,d\lambda\\
&\le C t^{-l_1/2}\prod_{\alpha\in F}\int_{i\R}\big| r_\alpha(\pi,z_\alpha)^{-1}
r_\alpha^\prime(\pi,z_\alpha)\big|e^{-tz_\alpha^2}\,dz_\alpha.
\end{split}
\end{equation}
Suppose that $M=\GL_{n_1}\times\cdots\times\GL_{n_r}$. Then $\pi=\pi_1\otimes
\cdots\otimes\pi_r$ with $\pi_i\in\Pi_{\di}(\GL_{n_i}(\A))$. Now recall that a
given root $\alpha\in\Sigma(G,A_M)$ corresponds to an ordered pair $(i,j)$
of distinct integers $i$ and $j$ between 1 and $r$. Then it folows from 
(\ref{4.4}) and (\ref{4.5}) that $r_\alpha(\pi,s)=r(\pi_i\otimes\pi_j,s)$.
Let $l=\dim(A_L/A_G)$ and $k=\dim(A_L/A_{L_1})$. Then by Proposition \ref{p4.2}
and (\ref{4.15}) it follows that there exists $C>0$ 
\begin{equation}\label{4.16}
\begin{split}
\int_{i\af_L^*/\af_G^*}|\nu_L^S(P,\pi,&\lambda)|
e^{-t\parallel\lambda\parallel^2}\,d\lambda\\
&\le C \prod_{i,j}\log(1+\nu(\pi_i\times\widetilde\pi_j))
\frac{(1+|\log t|)^k}{t^{l/2}}
\end{split}
\end{equation}
for all $0<t\le1$ and all $\pi\in\Pi_{\di}(M(\A))$.

Next we shall estimate the numbers $\nu(\pi_i\times\widetilde\pi_j)$. 
For $\pi_\infty\in\Pi(\GL_m(\R))$, 
let the complex numbers $\mu_j(\pi_\infty)$, $j=1,...,m$, be defined by 
(\ref{4.14n}) and set
$$c(\pi_\infty)=\left(\sum^{m}_{j=1}|\mu_j(\pi_\infty)|^2\right)^{1/2}.$$
Given an open compact subgroup $K_f$ of $\GL_m(\A_f)$, set
$$\Pi(\GL_m(\A))_{K_f}:=
\{\pi\in\Pi(\GL_m(\A))\mid\H_{\pi_f}^{K_f}\neq0\},$$
where $\pi=\pi_\infty\otimes\pi_f$ and $\H_{\pi_f}$ denotes the Hilbert space
 of the representation $\pi_f$.
\begin{lem}\label{l4.3}
Let $K_{f,i}\subset\GL_{m_i}(\A_f),$ $i=1,2,$ be two open compact 
subgroups. There exists $C>0$ such that
$$\nu(\pi_1\times\pi_2)\le C(1+c(\pi_{1,\infty})
+c(\pi_{2,\infty}))$$
for all $\pi_i\in\Pi(\GL_{m_i}(\A))_{K_{f,i}}$, $i=1,2$.
\end{lem}
\begin{proof}
First consider $c(\pi_1\times\pi_2)$ which is defined by 
(\ref{3.18}). It follows from Lemma \ref{l3.2} that there exists 
$C>0$ such that
$$c(\pi_1\times\pi_2)\le C(c(\pi_{1,\infty})+c(\pi_{2,\infty}))$$
for all $\pi_i\in\Pi(\GL_{m_i}(\A)),$ $i=1,2.$  
It remains to estimate $N(\pi_1\times\pi_2)$. For this we first observe that,
as the epsilon factor is a product of local epsilon factors, we can factor
$N(\pi_1\times\pi_2)$ as
$$N(\pi_1\times\pi_2)=\prod_{p}N(\pi_{1,p}\times\pi_{2,p}),$$
where $p$ runs over the finite places of $\Q$. This is a finite product. In 
fact, if $p$ is unramified for both $\pi_1$ and $\pi_2$, then 
$N(\pi_{1,p}\times\pi_{2,p})=1$. Moreover there is an integer 
$f(\pi_{1,p}\times\pi_{2,p})$ such that
$$N(\pi_{1,p}\times\pi_{2,p})=p^{f(\pi_{1,p}\times\pi_{2,p})}$$
(see e.g. \cite{MS}). Since we fix the ramification, there is a finite set 
$S$ of finite places of $\Q$, such that
$$N(\pi_1\times\pi_2)=\prod_{p\in S}p^{f(\pi_{1,p}\times\pi_{2,p})}$$
for all $\pi_i\in\Pi(\GL_{m_i}(\A))_{K_{f,i}}$, $i=1,2$. This reduces our 
problem to the 
estimation of $f(\pi_{1,p}\times\pi_{2,p})$. Let $f(\pi_{i,p})$ be the
 conductor
of $\pi_{i,p}$, $i=1,2$. Then $f(\pi_{i,p})\ge 0$ and by Theorem 1 of
 \cite{BH} and Corollary (6.5) of \cite{BHK} we have
\begin{equation}\label{5.17a}
0\le f(\pi_{1,p}\times\pi_{2,p})\le m_1 f(\pi_{1,p})+m_2 f(\pi_{2,p}).
\end{equation}
Let $m\in\N$ and let $K_p$ be an open compact subgroup of $\GL_m(\Q_p)$.
By Lemma 2.2 of \cite{MS} there exists $C_p>0$ such that $f(\pi_p)\le C_p$
for all $\pi_p\in\Pi(\GL_m(\Q_p))$ with $\pi_p^{K_p}\not=0$. Together with
(\ref{5.17a}) this implies that there exists $C>0$ such that
$$N(\pi_1\times\pi_2)\le C$$
for all $\pi_i\in\Pi(\GL_{m_i}(\A))_{K_{f,i}}$, $i=1,2$.
This completes the proof of the lemma.
\end{proof}

We continue with the estimation of $c(\pi_\infty)$. Given $\pi_\infty\in
\Pi(\GL_m(\R),\xi_0)$, let $\lambda_{\pi_\infty}$ be the Casimir eigenvalue of
the restriciton of $\pi_\infty$ to $\GL_m(\R)^1$. Furthermore for $\sigma\in
\Pi(\rO(m))$ let $\lambda_\sigma$ denote the Casimir eigenvalue of $\sigma$.
We note that if $[\pi_\infty|_{\rO(m)}:\sigma]>0$, then 
$-\lambda_{\pi_\infty}+\lambda_\sigma\ge 0$ \cite[Lemma 2.6]{DH} .

\begin{lem}\label{l4.3a}
There exists $C>0$ such that
$$c(\pi_\infty)\le C(1-\lambda_{\pi_\infty}+\lambda_\sigma)^{1/2}$$
for all 
$\pi_\infty\in\Pi(\GL_m(\R),\xi_0)$ and  $\sigma\in\Pi(\rO(m))$
 with $[\pi_\infty|_{\rO(m)}:\sigma]>0$.
\end{lem}
\begin{proof}
Write $\pi_\infty$ as Langlands quotient $\pi_\infty=J^{\GL_m}_R(\tau,\bs)$,
where $\tau$ is a discrete series representation of $M_R(\R)$ and the 
parameters $s_1,...,s_r\in\C$ 
satisfy $\Re(s_1)\ge\Re(s_2)\ge\cdots\ge\Re(s_r)$.
We may assume that the central character of $\tau$ is trivial on $A_R(\R)^0$
and hance, we can regard $\tau$ as a discrete series representation of
$M_R(\R)^1$. Let $\mf_R^1$ denote the Lie algebra of $M_R(\R)^1$. Note that
$\mf_R^1$ is the direct sum of a finite number of copies of 
${\mathfrak sl}(2,\R)$. Let $\tf\subset\mf_R^1$ be the standard compact
Cartan subalgebra equipped with the canonical norm. Then $\hf=\tf\oplus\af_R$
is a Cartan subalgebra of $\gl_m(\R)$. Let $\Lambda_\tau\in i\tf^*$ be the
Harish-Cahndra parameter of $\tau$. It follows from the definition of the
parameters $\mu_j(\pi_\infty)$ in terms of the Langlands parameters that there
exists $C>0$ such that
$$c(\pi_\infty)^2\le
C(\parallel\Lambda_\tau\parallel^2+\parallel\bs\parallel^2)$$
for all $\pi_\infty\in\Pi(\GL_m(\R),\xi_0)$. Let $\gamma:Z(\gl_m(\C))\to 
I(\hf_\C)$ be the Harish-Chandra homomorphism. By Proposition 8.22 of
\cite{Kn} the infinitesimal character $\chi$ of the induced representation
$I^{\GL_m}_R(\tau,\bs)$ with respect to $\hf$ is given by 
$\chi(Z)=(\Lambda_\tau+\bs)(\gamma(Z))$, $Z\in Z(\gl_m(\C))$. 
Since $\pi_\infty$ is an irreducible quotient
of $I^{\GL_m}_R(\tau,\bs)$, it has the same infinitesimal character. 
Let $H_1,...,H_r$ be an orthonormal basis of $\af_R$ and $H_{r+1},...,H_m$ an
orthonormal basis of $\tf$. Then
$$\gamma(\Omega)=\sum_{i=1}^rH_i^2-\sum_{j=r+1}^m
H_j^2-\parallel\rho\parallel^2$$
\cite[p.168]{Wa1}. Hence, the Casimir eigenvalue $\lambda_\pi$ of $\pi_\infty$ 
is given by
$$\lambda_{\pi_\infty}=(\Lambda_\tau+\bs)(\gamma(\Omega))=\sum_{i=1}^rs_i^2+
\parallel\Lambda_\tau\parallel^2-\parallel\rho\parallel^2.$$
Since $\pi_\infty$ is unitary, it follows from Theorem 3.3 of Chapter XI of
 \cite{BW}
that there exists $C>0$, independent of $\pi_\infty$,  such that 
$\parallel\Re(\bs)\parallel\le C$. Hence there exists $C_1>0$ such that
$$\parallel\Lambda_\tau\parallel^2+\parallel\bs\parallel^2
\le C_1-\lambda_{\pi_\infty}+\parallel\Lambda_\tau\parallel^2$$
for all $\pi_\infty\in\Pi(\GL_m(\R),\xi_0)$. 
Now let $\sigma\in\Pi(\rO(m))$ and suppose that 
$\pi_\infty\in\Pi(\GL_m(\R),\xi_0)$ is such that 
$[\pi_\infty|_{\rO(m)}:\sigma]>0$. Since
$\sigma$ occurs in $\pi_\infty$, it also occurs in $I^{\GL_m}_R(\tau,\bs)$. 
Using
Frobenius reciprocity as in \cite[p.208]{Kn}, it follows that there exists
$\omega\in\Pi(\rO(m)\cap M_R(\R))$ such that
$$[\tau|_{\rO(m)\cap M_R(\R)}:\omega]>0\quad\text{and}\quad
[\sigma|_{\rO(m)\cap M_R(\R)}:\omega]>0.$$
Let $\lambda_\sigma$ and $\lambda_\omega$ denote the Casimir eigenvalues of
$\sigma$ and $\omega$, respectively. By \cite[(5.15)]{Mu2}, the second
inequality implies $\lambda_\omega\le\lambda_\sigma$. On the other hand,
by \cite[p.398]{Wa2}, the first inequality implies
$$\parallel\Lambda_\tau\parallel^2\le \lambda_\omega
+\parallel\rho_R\parallel^2.$$
Combining our estimations the lemma follows.
\end{proof}

Now let $K_f$ be an open compact subgroup of $G(\A_f)$. Set
$$K_{M,f}=K_f\cap M(\A_f).$$
Then $K_{M,f}$ is an open compact subgroup of 
$M(\A_f)$. There exist open compact subgroups $K_{f,i}$ of 
$\GL_{n_i}(\A_f)$, $i=1,...,r$,
such that $K_{f,1}\times\cdots\times K_{f,r}$ is a subgroup of finite index of
$K_f$. Set
$$\Pi(M(\A),\xi_0)_{K_f}=\{\pi\in\Pi(M(\A),\xi_0)\,\big|\,\H_{\pi_f}^{K_{M,f}}
\not=\{0\}\},$$
where $\pi=\pi_\infty\otimes\pi_f$.
 Let $\pi\in\Pi(M(\A),\xi_0)_{K_f}$. Then $\pi=\pi_1\otimes\cdots\otimes\pi_r$
and $\pi_i$ belongs to $\Pi(\GL_{n_i}(\A),\xi_0)_{K_{f,i}}$ and by 
Lemma \ref{l4.3} it follows that there exists $C>0$ such that
\begin{equation}\label{4.17}
\prod_{i,j}\log(1+\nu(\pi_i\times\widetilde\pi_j))\le
C \prod_{i,j}\log(2+c(\pi_{i,\infty})+c(\pi_{j,\infty}))
\end{equation}
for all $\pi=\pi_1\otimes\cdots\otimes\pi_r\in\Pi(M(\A),\xi_0)_{K_f}$. 

Let $K_{M,\infty}=\rO(n_1)\times\cdots\times\rO(n_r)$ be the standard maximal
compact subgroup of $M(\R)$. Let $\sigma\in\Pi(\rO(n))$. For
$\pi\in\Pi(M(\A),\xi_0)$ set
$$[\pi_\infty:\sigma]=\sum_{\tau\in\Pi(K_{M,\infty})}
[\pi_\infty|_{K_{M,\infty}}:\tau][\sigma|_{K_{M,\infty}}:\tau].$$
Put
$$\Pi(M(\A),\xi_0)_{K_f,\sigma}=
\{\pi\in\Pi(M(\A),\xi_0)_{K_f}\;|\;[\pi_\infty:\sigma]>0\}
$$
and
$$\Pi_{\di}(M(\A),\xi_0)_{K_f,\sigma}
=\Pi_{\di}(M(\A),\xi_0)\cap\Pi(M(\A),\xi_0)_{K_f,\sigma}.$$
Suppose that $\pi\in\Pi(M(\A),\xi_0)_{K_f,\sigma}$. Let 
$\tau\in\Pi(K_{M,\infty})$ be such that $[\sigma|_{K_{M,\infty}}:\tau]>0$
and $[\pi_\infty|_{K_{M,\infty}}:\tau]>0$.

Let $\lambda_{\pi_\infty}$ and $\lambda_\tau$
denote the Casimir eigenvalues of the restriction
of $\pi_\infty$ to $M(\R)^1$ and of $\tau$, respectively. Note that
$\lambda_{\pi_\infty}=\sum_i\lambda_{\pi_{i,\infty}}$ and $\lambda_\tau=
\sum_i\lambda_{\tau_i}$, where $\tau=\otimes_i\tau_i$.
Then it follows from (\ref{4.17}) and  Lemma \ref{l4.3a} that there exists
 $C>0$ such that
\begin{equation*}
\begin{split}
\prod_{i,j}\log(1+\nu(\pi_i\times\widetilde\pi_j))&\le
C\prod_{i,j}\log(2-\lambda_{\pi_{i,\infty}}+\lambda_{\tau_i}
-\lambda_{\pi_{j,\infty}}+\lambda_{\tau_j})\\
&\le C \bigl(\log(2-\lambda_{\pi_\infty}+\lambda_\tau)\bigr)^{r^2}
\end{split}
\end{equation*}
for all $\pi\in\Pi(M(\A),\xi_0)_{K_f,\sigma}$. Since there are only finitely
many $\tau$ that occur in $\sigma|_{K_{M,\infty}}$, we get
\begin{equation}\label{4.18}
\prod_{i,j}\log(1+\nu(\pi_i\times\widetilde\pi_j))\le C_1
\bigl(\log(2+|\lambda_{\pi_\infty}|)\bigr)^{r^2}
\end{equation}
for all $\pi\in\Pi(M(\A),\xi_0)_{K_f,\sigma}$.
Combining (\ref{4.16})--(\ref{4.18}) we obtain
\begin{prop}\label{p4.5}
Let $M\in\L$, $L\in\L(M)$ and $P\in\cP(M)$. Let $l=\dim(A_L/A_G)$.
 Let $K_f$ be an open compact
subgroup of $\GL_n(\A_f)$ and let $\sigma\in\Pi(\rO(n))$.
There exists $C>0$ such that
$$\int_{i\af_L^*/\af_G^*}|\nu^S_L(P,\pi,\lambda)|
e^{-t\parallel\lambda\parallel^2}\;d\lambda\le 
C\left(\log(2+|\lambda_{\pi_\infty}|)\right)^{n^2}
\frac{(1+|\log t|)^l}{t^{l/2}}$$
for all $0<t\le1$ and $\pi\in\Pi_{\di}(M(\A),\xi_0)_{K_f,\sigma}$.
\end{prop}

\section{The spectral side}
\setcounter{equation}{0}

We shall use the noninvariant trace formula of Arthur \cite{A1}, \cite{A2},
applied to the heat kernel, to determine the growth of the discrete spectrum.
To begin with, we explain the general structure of the spectral side of the 
Arthur trace formula. The spectral side is a sum of distributions 
$$\sum_{\chi\in\mX} J_\chi(f),\quad \chi\in C^\infty_0(G(\A)^1).$$
Here $\mX$ is the set of cuspidal datas which consists of Weyl group orbits
of pairs $(M_B,\rho_B)$, where $M_B$ is the Levi component of a parabolic 
subgroup and $\rho_B$ is a cuspidal automorphic representation of $M_B(\A)$.
The distributions $J_\chi$ are described by Theorem 8.2 of \cite{A4}. Let
$\Co^1(G(\A)^1)$ be the space of integrable rapidly decreasing functions
on $G(\A)^1$ \cite[\S 1.3]{MS}. In \cite[Theorem 0.1]{MS} it has been proved
that the spectral side of the trace formula for $\GL_n$ is absolutely 
convergent for all
 $f\in\Co^1(G(\A)^1)$. In this case the expression of the
spectral side simplifies.

To describe this in more detail, we need to introduce some notation. Let
$M\in\cL$ and  $P,Q\in\cP(M)$.  Let $\cA^2(P)$ and $\cA^2(Q)$ be the 
corresponding spaces of automorphic functions (see \S1.5).
Let $W(\af_P,\af_Q)$ be the set of all linear
isomorphisms from $\af_P$ to $\af_Q$ which are restrictions of elements of the
Weyl group $W(A_0)$. The theory of Eisenstein series associates to each 
$s\in W(\af_P,\af_Q)$ an intertwining operator
$$M_{Q|P}(s,\lambda):\cA^2(P)\to\cA^2(Q),\quad \lambda\in\af_{P,\C}^*,$$
which, for $\Re(\lambda)$ in a certain chamber, can be defined by an 
absolutely convergent integral and admits an analytic continuation to a
 meromorphic 
function of $\lambda\in\af_{P,\C}^*$ \cite{La}. Set
$$M_{Q|P}(\lambda):=M_{Q|P}(1,\lambda).$$
Fix $P\in\cP(M)$  and $\lambda\in i\af_M^*$. 
For $Q\in\cP(M)$ and $\Lambda\in i\af_M^*$ define
$$\mM_Q(P,\lambda,\Lambda)=M_{Q|P}(\lambda)^{-1}M_{Q|P}(\lambda+\Lambda).$$
Then
\begin{equation}\label{5.1a}
\{\mM_Q(P,\lambda,\Lambda)\mid \Lambda\in i\af_M^*,\;Q\in\cP(M)\}
\end{equation}
is a $(G,M)$ family with values in the space of operators on $\cA^2(P)$
\cite[p.1310]{A4}. 
Let $L\in\cL(M)$. Then, as explained in the previous section,
  the $(G,M)$ family (\ref{5.1a}) has an associated $(G,L)$ family
$$\{\mM_{Q_1}(P,\lambda,\Lambda)\mid \Lambda\in i\af_L^*,\;Q_1\in\cP(L)\}$$
and
$$\mM_L(P,\lambda,\Lambda)=\sum_{Q_1\in\cP(L)}\mM_{Q_1}(P,\lambda,\Lambda)
\theta_{Q_1}(\Lambda)^{-1}$$
extends to a smooth function on $i\af_L^*$. Put
$$\mM_L(P,\lambda)=\mM_L(P,\lambda,0).$$
This operator depends only on the intertwining operators. It equals
\begin{equation*}
\begin{split}
\mM_L&(P,\lambda)=\\
&\lim_{\Lambda\to0}\left(\sum_{Q_1\in\cP(L)}
\vol(\af_{Q_1}^G/\Z(\Delta^\vee_{Q_1}))M_{Q|P}(\lambda)^{-1}
\frac{M_{Q|P}(\lambda+\Lambda)}{\prod_{\alpha\in\Delta_{Q_1}}
\Lambda(\alpha^\vee)}\right),
\end{split}
\end{equation*}
where $\lambda$ and $\Lambda$ are constrained to lie in $i\af_L^*$, and for 
each $Q_1\in\cP(L)$, $Q$ is a group in $\cP(M_P)$ which is contained in $Q_1$.
Then $\mM_L(P,\lambda)$ is an unbounded operator which acts
on the Hilbert space $\ov \cA^2(P)$.  For $\pi\in\Pi(M(\A)^1)$ let 
$\cA^2_\pi(P)$ be the subspace of $\cA^2(P)$ determined by $\pi$ (see \S1.5). 
Let
$\rho_\pi(P,\lambda)$ be the induced representation of $G(\A)$ in 
$\ov \cA^2_\pi(P)$. Let $W^L(\af_M)_{\reg}$ be the set of
elements $s\in W(\af_M)$ such that $\{H\in\af_M\mid sH=H\}=\af_L$.
For any function $f\in\Co^1(G(\A)^1)$ and 
$s\in W^L(\af_M)_{\reg}$ set
\begin{equation}\label{5.1}
\begin{split}
J^L_{M,P}&(f,s)\\
&=\sum_{\pi\in\Pi_{\di}(M(\A)^1)}
\int_{i\af_L^*/i\af_G^*}\tr(\mM_L(P,\lambda)M_{P|P}(s,0)
\rho_\pi(P,\lambda,f))\;d\lambda.
\end{split}
\end{equation}

By Theorem 0.1 of \cite{MS} this integral-series is absolutely convergent with 
respect to the trace norm.
 Furthermore for $M\in\cL$ and $s\in W^L(\af_M)_{\reg}$ set
$$a_{M,s}=|\cP(M)|^{-1}|W^M_0||W_0|^{-1}|\det(s-1)_{\af^L_M}|^{-1}.$$
Then for any  $f$ in $\Co^1(G(\A)^1)$, the spectral  side $J_{\spec}(f)$
of the Arthur trace formula is given by
\begin{equation}\label{5.2}
J_{\spec}(f)=\sum_{M\in\cL}\sum_{L\in\cL(M)}\sum_{P\in\cP(M)}
\sum_{s\in W^L(\af_M)_{\reg}}
 a_{M,s} J^L_{M,P}(f,s).
\end{equation}
Note that all sums in this expression are finite. 

We shall now evaluate the spectral side at a function $\phi_t$, $t>0$, which 
is given in terms of the heat kernel of a Bochner-Laplace operator. Then our 
main purpose is to determine the behaviour of $J_{\spec}(\phi_t)$ as $t\to0$. 

Let $G(\R)^1=G(\A)^1\cap G(\R)$. By definition $G(\R)^1$ consisits of all
$g\in G(\R)$ with $|\det(g)|=1$. Hence $G(\R)^1$ is semisimple and
$$G(\R)=G(\R)^1\cdot A_G(\R)^0.$$
Let 
$$X=G(\R)^1/K_\infty$$
be the associated Riemannian symmetric space. Given $\sigma\in\Pi(K_\infty)$, 
let $\widetilde E_\sigma\to X$ be the associated homogeneous vector bundle. 
Let $\Omega_{G(\R)^1}$ be the Casimir element of $G(\R)^1$ and let
$\widetilde\Delta_\sigma$ be the operator in $L^2(\widetilde E_\sigma)$
which is induced by $-R(\Omega_{G(\R)^1})\otimes\Id$. Let 
\begin{equation}\label{5.2a}
H_t^\sigma\in\left(\Co^1(G(\R)^1)\otimes
\End(V_\sigma)\right)^{K_\infty\times K_\infty}
\end{equation}
be the kernel of the heat operator $e^{-t\widetilde\Delta_\sigma}$ where
$\Co^1(G(\R))$ is Harish-Chandra's space of integrable rapidly decreasing
functions.  Set
$$h_t^\sigma=\tr H_t^\sigma.$$
We extend $h_t^\sigma$ to a function on $G(\R)$ by
$$h_t^\sigma(g\cdot z)=h_t^\sigma(g),\quad g\in G(\R)^1,\; z\in A_G(\R)^0.$$
Then $h_t^\sigma$ satisfies
$$h_t^\sigma(gz)=h_t^\sigma(g),\quad g\in G(\R),\;z\in A_G(\R)^0.$$
Let $\chi_\sigma$ be the character of $\sigma$. Then  $h_t^\sigma$ also 
satisfies
$$h_t^\sigma=\chi_\sigma\ast h_t^\sigma\ast\overline{\chi_\sigma}.$$
Let $K_f$ be an open compact subgroup of $G(\A_f)$ and let ${\bf 1}_{K_f}$ be
the charcteristic function of $K_f$ in $G(\A_f)$. Set
$$\chi_{K_f}=\vol(K_f)^{-1}{\bf 1}_{K_f}.$$
Define the function $\phi_t$ on $G(\A)$ by
\begin{equation}\label{5.3}
\phi_t(g)=h_t^\sigma(g_\infty)\chi_{K_f}(g_f)
\end{equation}
for any point 
$$g=g_\infty g_f,\quad g_\infty\in G(\R),\; g_f\in G(\A_f),$$
in $G(\A)$. Then $\phi_t$ satisfies $\phi_t(gz)=\phi_t(g)$ for 
$z\in A_G(\R)^0$, $g\in G(\A)$. It follows from (\ref{5.2a}) and the 
definition of $\Co^1(G(\A)^1)$  that the
 restriction $\phi_t^1$ of $\phi_t$ to $G(\A)^1$ belongs to $\Co^1(G(\A)^1)$.

 Let $\pi$ be any unitary representation of
 $G(\A)$ which is trivial on $A_G(\R)^0$. Then we can define
$$\pi(\phi_t)=\int_{G(\A)/A_G(\R)^0} \phi_t(g)\pi(g)\,dg.$$
Suppose that $\pi=\pi_\infty\otimes\pi_f$, where $\pi_\infty$ and $\pi_f$ are
unitary representations of $G(\R)$ and $G(\A_f)$, respectively.
Then $\pi_\infty$ is trivial on $A_G(\R)^0$. So we can set
$$\pi_\infty(\phi_t)=\int_{G(\R)/A_G(\R)^0}\pi_\infty(g_\infty)
h_t^\sigma(g_\infty)\,dg_\infty.$$
Let $\Pi_{K_f}$ denote the orthogonal projection of the Hilbert space
$\H_{\pi_f}$ of $\pi_f$ onto the subspace $\H_{\pi_f}^{K_f}$ of $K_f$-invariant
vectors. Then
$$\pi(\phi_t)=\pi_\infty(h_t^\sigma)\otimes\Pi_{K_f}.$$
Now let $\pi\in\Pi(M(\A)^1)$. We identify $\pi$ with a representation of 
$M(\A)$ which is trivial on $A_M(\R)^0$. Let $I_P^G(\pi_\lambda)$, 
$\lambda\in\af_{M,\C}^*$, be the induced representation of $G(\A)$. 
 Let $\pi=\pi_\infty\otimes
\pi_f$. Then
$$I_P^G(\pi_\lambda)=I_P^G(\pi_{\infty,\lambda})\otimes I_P^G(\pi_{f,\lambda}).
$$
Let $\H_P(\pi_\infty)_\sigma$
denote the $\sigma$-isotypical subspace of the Hilbert space 
$\H_P(\pi_\infty)$ 
of the induced representation.  Then $\H_P(\pi_\infty)_\sigma$ is an invariant
subspace of $I_P^G(\pi_{\infty,\lambda},h_t^\sigma)$. 
Let $\lambda_{\pi}$ be  the Casimir eigenvalue
of the restriction of $\pi_\infty$ to $M(\R)^1$.
By  Proposition 8.22 of \cite{Kn} it follows that
$$I_P^G(\pi_{\infty,\lambda},h_t)\upharpoonright\H_P(\pi_\infty)_\sigma=
e^{-t(-\lambda_{\pi}+\parallel\lambda\parallel^2)}
\Id.$$
Now observe that there is a canonical isomorphism
$$j_P:\H_P(\pi)\otimes\Hom_{M(\A)}
\left(\pi,I_{M(\Q)A_M(\R)^0}^{M(\A)}(\xi_0)\right)\to\overline{\cA}_\pi^2(P),$$
which intertwines the induced representations. Let $\Pi_{K_f,\sigma}$ denote 
the orthogonal projection of $\ov\cA^2_\pi(P)$ onto 
$\cA^2_\pi(P)_{K_f,\sigma}$. 
Then it follows that
\begin{equation}\label{5.3a}
\rho_\pi(P,\lambda,\phi_t)=e^{-(-\lambda_\pi+\parallel\lambda\parallel^2)}
\Pi_{K_f,\sigma}.
\end{equation}
Suppose that $\lambda\in(\af^G_P)^*_\C$. Then $\rho_\pi(P,\lambda,g)$ is
trivial on $A_G(\R)^0$. This implies $\rho_\pi(P,\lambda,\phi_t)=
\rho_\pi(P,\lambda,\phi_t^1)$, where $\phi_t^1$ is the restriciton of 
$\phi_t$ to $G(\A)^1$. 
Together with (\ref{5.3a}) we get
\begin{equation}\label{5.4}
\begin{split}
J^L_{M,P}(\phi_t^1,s)=&\sum_{\pi\in\Pi_{\di}(M(\A)^1)}
e^{t\lambda_\pi}\\
&\cdot\int_{i\af_L^*/\af_G^*}e^{-t\parallel\lambda\parallel^2}
\tr(\mM_L(P,\lambda)M_{P|P}(s,0)\Pi_{K_f,\sigma})\;d\lambda.
\end{split}
\end{equation}
To study this integral-series, we introduce the normalized intertwining
operators
\begin{equation}\label{5.5}
N_{Q|P}(\pi,\lambda):=r_{Q|P}(\pi,\lambda)^{-1}M_{Q|P}(\pi,\lambda),\quad
\lambda\in\af_{M,\C}^*,
\end{equation}
where $r_{Q|P}(\pi,\lambda)$ are the global normalizing factors considered in
the previous section. Let $P\in\cP(M)$ and $\lambda\in i\af_M^*$ be fixed.
For $Q\in\cP(M)$ and  $\Lambda\in i\af_M^*$ define
\begin{equation}\label{5.6}
\mN_Q(P,\pi,\lambda,\Lambda)=N_{Q|P}(\pi,\lambda)^{-1}N_{Q|P}(\pi,\lambda+\Lambda),
\end{equation}
Then as functions of $\Lambda\in i\af_M^*$, 
$$\{\mN_Q(P,\pi,\lambda,\Lambda)\mid Q\in\cP(M)\}$$
is  a $(G,M)$ family. The verification is the same as in the case of
the unnormalized intertwining operators \cite[p.1310]{A4}. For $L\in\cL(M)$, 
let
$$\{\mN_{Q_1}(P,\pi,\lambda,\Lambda)\mid \Lambda\in i\af_L^*,\;Q_1\in\cP(L)\}
$$
be the associated $(G,L)$ family. 

Let $\mM_{Q_1}(P,\pi,\lambda,\Lambda)$ be the restriction of 
$\mM_{Q_1}(P,\lambda,\Lambda)$ to $\ov\cA^2_\pi(P)$. 
Then by (\ref{5.5}) and (\ref{4.13}) it 
follows that
\begin{equation}\label{5.7}
\mM_{Q_1}(P,\pi,\lambda,\Lambda)=\mN_{Q_1}(P,\pi,\lambda,\Lambda)
\nu_{Q_1}(P,\pi,\lambda,\Lambda)
\end{equation}
for all $\Lambda\in i\af_L^*$ and all $Q_1\in\cP(L)$. 

For $Q\supset P$ let
$\hat L_P^Q\subset \af_P^Q$ be the lattice generated by 
$\{\tilde\omega^\vee\mid \tilde\omega\in\hat\Delta_P^Q\}$. Define
$$\hat\theta_P^Q(\lambda)=\vol(\af_P^Q/\hat L_P^Q)^{-1}
\prod_{\tilde\omega\in\hat\Delta_P^Q}\lambda(\tilde\omega^\vee).$$
For $S\in\cF(L)$ put
\begin{equation}\label{5.8}
\begin{split}
\mN_S^\prime&(P,\pi,\lambda)\\
&=\lim_{\Lambda\to0}\sum_{\{R\mid R\supset S\}}(-1)^{\dim(A_S/A_R)}
\hat\theta_S^R(\Lambda)^{-1}\mN_R(P,\pi,\lambda,\Lambda)\theta_R(\Lambda)^{-1}.
\end{split}
\end{equation}
Let $\mM_L(P,\pi,\lambda)$ be the restriction of 
$\mM_L(P,\lambda)$ to $\ov\cA^2_\pi(P)$. 
Then by (\ref{5.7}) and  Lemma 6.3 of \cite{A5} we get
\begin{equation}\label{5.9}
\mM_L(P,\pi,\lambda)=\sum_{S\in\cF(L)}\mN_S^\prime(P,\pi,\lambda)
\nu_L^S(P,\pi,\lambda).
\end{equation}
Let $\mN_S^\prime(P,\pi,\lambda)_{K_f,\sigma}$ denote the restriction of
$\mN_S^\prime(P,\pi,\lambda)$ to $\cA^2_\pi(P)_{K_f,\sigma}$. Then by
(\ref{5.4}) we get
\begin{equation}\label{5.10}
\begin{split}
J^L_{M,P}&(\phi_t^1,s)=\sum_{\pi\in\Pi_{\di}(M(\A)^1)}
e^{t\lambda_\pi}\\
&\hskip-5pt\cdot\sum_{S\in\cF(L)}\int_{i\af_L^*/\af_G^*}
e^{-t\parallel\lambda\parallel^2}
\nu_L^S(P,\pi,\lambda)\tr(M_{P|P}(s,0)
\mN_S^\prime(P,\pi,\lambda)_{K_f,\sigma})\;d\lambda.
\end{split}
\end{equation} 
Next we shall estimate the norm of $\mN_S^\prime(P,\pi,\lambda)_{K_f,\sigma}$.
For a given place $v$ of $\Q$ let $J_{Q|P}(\pi_v,\lambda)$ be the  intertwining
 operator between the
induced representations $I_P^G(\pi_{v,\lambda})$ and $I_Q^G(\pi_{v,\lambda})$.
Let
$$R_{Q|P}(\pi_v,\lambda)=r_{Q|P}(\pi_v,\lambda)^{-1}J_{Q|P}(\pi_v,\lambda),
\quad\lambda\in\af_{M,\C}^*,$$
be the normalized local intertwining operator. These operators satisfy the
conditions $(R_1)-(R_8)$ of Theorem 2.1 of \cite{A7}. Assume that
$K_f=\prod_{p<\infty}K_p$. For any place $v$ denote by $\H_P(\pi_v)$ the
Hilbert space of the induced representation $I_P^G(\pi_{v})$. If $p<\infty$
let $R_{Q|P}(\pi_p,\lambda)_{K_p}$ be the restriction of  
$R_{Q|P}(\pi_p,\lambda)$ to the subspace of $K_p$-invariant vectors
$\H_P(\pi_p)^{K_p}$ in $\H_P(\pi_p)$. Let
$R_{Q|P}(\pi_\infty,\lambda)_\sigma$ denote the restriction of 
$R_{Q|P}(\pi_\infty,\lambda)$ to the $\sigma$-isotypical subspace of 
$I_P^G(\pi_\infty)$ in $\H_P(\pi_\infty)$. It was proved in \cite[(6.24)]{Mu2}
that there exist a finite set of places $S_0$, including the Archimedean
one, and constants $C>0$ and $q\in\N$, such that
\begin{equation*}
\begin{split}
\parallel \mN_S^\prime(P,\pi,\lambda)_{K_f,\sigma}\parallel\le C
\bigg(\sum_{p\in S_0\setminus\{\infty\}}\sum_{k=1}^q
&\parallel D_\lambda^k R_{Q|P}(\pi_p,\lambda)_{K_p}\parallel\\
&\sum_{k=1}^q\parallel D_\lambda^k R_{Q|P}(\pi_\infty,\lambda)_\sigma\parallel
\bigg)
\end{split}
\end{equation*}
for all $\lambda\in i\af_{M}^*$, $\sigma\in\Pi(K_\infty)$ and $\pi\in
\Pi(M(\A))$. By Proposition 0.2 of \cite{MS} it follows that there exists
$C>0$ such that
\begin{equation}\label{5.11}
\parallel \mN_S^\prime(P,\pi,\lambda)_{K_f,\sigma}\parallel\le C
\end{equation}
for all $\lambda\in i\af_{M}^*$ and $\pi\in\Pi_{\di}(M(\A)^1)$.
Observe that $M_{P|P}(s,0)$ is unitary. Let $l=\dim( A_L/A_G)$.
 Using (\ref{5.10}), (\ref{5.11}) and
Proposition \ref{p4.5} it follows that there exists $C>0$ such that
\begin{equation}\label{5.11a}
\begin{split}
|J^L_{M,P}&(\phi_t^1,s)|\le C\frac{(2+|\log t|)^l}{t^{l/2}}\\
&\cdot\sum_{\pi\in\Pi_{\di}(M(\A)^1)}
\dim\cA^2_\pi(P)_{K_f,\sigma}
\bigl(\log(1+|\lambda_\pi|)\bigr)^{n^2}
e^{t\lambda_\pi}
\end{split}
\end{equation}
for all $0<t\le 1$. The series can be estimated using Proposition \ref{p2.5}.
Let $X_M=M(\R)/K_{M,\infty}^\prime$ and let $m=\dim X_M$. 
It follows from  Proposition \ref{p2.5}  that for every $\epsilon>0$ 
there exists $C>0$ such that the series is bounded by $C t^{-m/2-\epsilon}$ 
for $0<t\le 1$. Together with (\ref{5.11a}) we obtain the following
 proposition.
\begin{prop}\label{p5.1}
Let $m=\dim X_M$ and $l=\dim A_L/A_G$. For every $\epsilon>0$ there exists $C>0$ such that
$$|J^L_{M,P}(\phi_t^1,s)|\le C t^{-(m+l)/2-\epsilon}$$
for all $0< t\le 1$.
\end{prop}

Now we distinguish two cases. First assume that $M=G$. Then $L=P=G$ and $s=1$.
Let $R^1_{\di}$ be the restriciton of the regular representation $R^1$ of
$G(\A)^1$ in $L^2(G(\Q)\ba G(\A)^1)$ to the discrete subspace. Then
$J^G_{G,G}(\phi_t^1,1)=\Tr R^1_{\di}(\phi_t^1).$
Let $R_{\di}$ be the regular representation of $G(\A)$ in 
$L^2_{\di}(A_G(\R)^0G(\Q)\ba G(\A))$. Then the operator $R_{\di}(\phi_t)$ is 
isomorphic to
$R^1_{\di}(\phi_t^1)$. Thus
$$J^G_{G,G}(\phi_t^1,1)=\Tr R_{\di}(\phi_t).$$
Given $\pi\in\Pi_{\di}(G(\A),\xi_0)$, let $m(\pi)$ denote the multiplicity 
with which $\pi$ occurs in the regular representation of $G(\A)$ in 
$L^2(\A_G(\R)^0G(\Q)\ba G(\A))$. Then  we get
\begin{equation}\label{5.12}
\begin{split}
J^G_{G,G}&(\phi_t^1,1)\\
&=\sum_{\pi\in\Pi_{\di}(G(\A),\xi_0)}
m(\pi)\dim\bigl(\H_{\pi_f}^{K_f}\bigr)\dim\bigl(\H_{\pi_\infty}(\sigma)\bigr)
e^{t\lambda_\pi}.
\end{split}
\end{equation}
Now assume that $M\not=G$ is a proper Levi subgroup. Let $P=MN$. Let
$X=G(\R)^1/K_\infty$. Then
$$X\cong X_M\times A_M(\R)^0/A_G(\R)^0\times N(\R).$$
Since $l=\dim A_L/A_G\le \dim A_M/A_G$, it follows that $m+l\le\dim X-1$.
 Thus together with Proposition \ref{p5.1} we get
\begin{theo}\label{th5.2}
Let $d=\dim X$. For every open compact subgroup $K_f$ of $G(\A_f)$ and every
$\sigma\in\Pi(\rO(n))$ the spectral side of the trace formula, evaluated at
 $\phi_t^1$, satisfies
\begin{equation}\label{5.13}
\begin{split}
J_{\spec}&(\phi_t^1)\\
&=\sum_{\pi\in\Pi_{\di}(G(\A),\xi_0)}m(\pi)
\dim\bigl(\H_{\pi_f}^{K_f}\bigr)\dim\bigl(\H_{\pi_\infty}(\sigma)\bigr)
e^{t\lambda_\pi}\\
&\hskip20pt+O(t^{-(d-1)/2})
\end{split}
\end{equation}
as $t\to0^+$.
\end{theo}
This theorem can be restated in a slightly different way as follows.
There exist arithmetic subgroups $\Gamma_i\subset G(\R)$, $i=1,...,m,$
such that
$$A_G(\R)^0G(\Q)\ba G(\A)/K_f\cong \bigsqcup_{i=1}^m(\Gamma_i\ba G(\R)^1)$$
(cf. \cite[section 9]{Mu1}).
 Let $\Delta_{\sigma,i}$ be the operator induced by the negative of the
Casimir operator in $C^\infty(\Gamma_i\ba G(\R)^1,\sigma)$, $i=1,...,m$.  Let
$$\lambda_0\le\lambda_1\le \lambda_2\le\cdots$$
be the $L^2$-eigenvalues of $\Delta_\sigma=\oplus_{i=1}^m\Delta_{\sigma,i}$,
where each eigenvalue is counted with its multiplicity. Let $d=\dim X$.
If we proceed in the same way as in the proof of Lemma \ref{l2.2}, 
then it follows that (\ref{5.13}) is equivalent to
\begin{equation}\label{5.14}
J_{\spec}(\phi_t^1)=\sum_i e^{-t\lambda_i}+O(t^{-(d-1)/2})
\end{equation}
as $t\to0^+$. 

Let  $\Gamma(N)\subset \SL_n(\Z)$ be the
principal congruence subgroup of level $N$. Let $\mu_0\le\mu_1\le\cdots$ be
 the eigenvalues, counted with multiplicity, of
$\Delta_\sigma$ acting in $L^2(\Gamma(N)\ba \SL_n(\R),\sigma)$. Then it follows
from (\ref{5.14}) and (\ref{2.10}) that
\begin{equation}\label{5.15}
J_{\spec}(\phi_t^1)=\varphi(N)\sum_i e^{-t\mu_i}+O(t^{-(d-1)/2})
\end{equation}
as $t\to0^+.$

Our next purpose is to study $J_{\spec}$ as a functional on the Schwartz space.
Let $K_f$ be an open compact subgroup of $G(\A_f)$ and let
 $\sigma\in\Pi(K_\infty)$. Denote by $\Co^1(G(\A)^1;K_f,\sigma)$ the set of
 all $h\in\Co^1(G(\A)^1)$ which are bi-invariant under $K_f$ and transform
 under $K_\infty$ according to $\sigma$. Let $\Delta_G$ be the Laplace
operator of $G(\R)^1$. 
Then we have
\begin{prop}\label{p5.3}
For every open compact subgroup $K_f$ of $G(\A_f)$ and every $\sigma\in
\Pi(K_\infty)$ there exist $C>0$ and $k\in\N$ such that
$$|J_{\spec}(f)|\le C\parallel (\Id+\Delta_G)^kf\parallel_{L^1(G(\A)^1)}$$
for all $f\in\Co^1(G(\A)^1;K_f,\sigma)$.
\end{prop}

\begin{proof}
This follows essentially from the proof of Theorem 0.2 in \cite{Mu2} combined 
with Proposition 0.2 of \cite{MS}. We include some details. 
Let $M\in\cL$, $L\in\cL(M)$ and $P\in\cP(M)$. By (\ref{5.2}) it suffices to
estimate $J^L_{M,P}(f,s)$. Since $M_{P|P}(s,0)$ is unitary, it follows
from (\ref{5.1}) that
\begin{equation*}
|J^L_{M,P}(f,s)|\le\sum_{\pi\in\Pi_{\di}(M(\A)^1)}\int_{i\af_L^*/\af_G^*}
\parallel\mM_L(P,\lambda)\rho_{\pi}(P,\lambda,f)\parallel_1d\lambda,
\end{equation*}
where $\parallel\cdot\parallel_1$ denotes the trace norm for operators 
in the Hilbert space $\ov\cA^2_\pi(P)$. Using (\ref{5.9}) it follows that the
right hand side is bounded by
\begin{equation*}
\sum_{\pi\in\Pi_{\di}(M(\A)^1)}\int_{i\af_L^*/i\af_G^*}
\parallel\mN'_S(P,\pi,\lambda)\rho_\pi(P,\lambda,f)\parallel_1|\nu^S_L(P,\pi,
\lambda)|\;d\lambda
\end{equation*}
The function $\nu^S_L(P,\pi,\lambda)$ can be estimated by Theorem 5.4 of
\cite{Mu2}. This reduces our problem to the estimation of the trace norm of 
the operator
$\mN'_S(P,\pi,\lambda)\rho_\pi(P,\lambda,f)$. 
Let $K_f$ be an open compact subgroup of $G(\A_f)$ and let 
$\sigma\in\Pi(K_\infty)$. Denote by $\Pi_{K_f,\sigma}$ the orthogonal 
projection of the Hilbert
space $\ov \cA_\pi^2(P)$ onto the finite-dimensional subspace 
$\cA_\pi^2(P)_{K_f,\sigma}$.
 Let $f\in\Co^1(G(\A)^1;K_f,\sigma)$. Then
\begin{equation*}
\rho_\pi(P,\lambda,f)=\Pi_{K_f,\sigma}\circ\rho_\pi(P,\lambda,f)
\circ\Pi_{K_f,\sigma}
\end{equation*}
for all $\pi\in\Pi(M(A)^1)$. Let 
$$D=\Id+\Delta_G.$$ 
For any $k\in\N$ let $\rho_{\pi}(P,\lambda,D^{2k})_{K_f,\sigma}$ denote 
the restriction of the operator $\rho_{\pi}(P,\lambda,D^{2k})$ to the subspace
$\cA_\pi^2(P)_{K_f,\sigma}$. Then we get
\begin{equation}\label{5.17}
\begin{split}
\parallel \mN_S^\prime&(P,\pi,\lambda)\rho_{\pi}(P,\lambda,f)\parallel_1\\
&\le\parallel \mN_S^\prime(P,\pi,\lambda)_{K_f,\sigma}
\parallel\cdot
\parallel\rho_{\pi}(P,\lambda,D^{2k})^{-1}_{K_f,\sigma}
\parallel_1\\
&\hskip5truecm\cdot\parallel\rho_{\pi}(P,\lambda,D^{2k}f)\parallel,
\end{split}
\end{equation}
By (6.9) of \cite{Mu2} we get
\begin{equation}\label{5.18}
\parallel\rho_\pi(P,\lambda,D^{2k})^{-1}_{K_f,\sigma}\parallel\le
C\frac{\dim\cA^2_\pi(P)_{K_f,\sigma}}{(1+\parallel\lambda\parallel^2+\lambda_\pi^2)^k},
\end{equation}
and since $\rho_\pi(P,\lambda)$ is unitary, we have
\begin{equation}\label{5.19}
\parallel \rho_\pi(P,\lambda,D^{2k}f)\parallel\le \parallel D^{2k}f
\parallel_{L^1(G(\A)^1)}.
\end{equation}
Together with (\ref{5.11}) it follows that there exists $C>0$ such that
\begin{equation}\label{5.20}
\begin{split}
\parallel \mN_S^\prime&(P,\pi,\lambda)\rho_{\pi}(P,\lambda,f)\parallel_1\\
&\le C\parallel D^{2k}f\parallel_{L^1(G(\A)^1)}
(1+\parallel\lambda\parallel)^{-k/2}\frac{\dim\cA^2_\pi(P)_{K_f,\sigma}}
{(1+\lambda_\pi^2)^{k/2}}
\end{split}
\end{equation}
for all $\lambda\in i\af_M^*$ and $\pi\in\Pi_{\di}(M(\A)^1)$. Let 
$d=\dim G(\R)^1/K_\infty$. By Theorem 5.4 of \cite{Mu2} there exists
$k_0\in\N$ such that for $k\ge k_0$ we have
\begin{equation}\label{5.21}
\int_{i\af_L^*/\af_G^*}|\nu^S_L(P,\pi,\lambda)|
(1+\parallel\lambda\parallel^2)^{-k/2}\,d\lambda\le
C_k(1+\lambda_\pi^2)^{8d^2}
\end{equation}
for all $\pi\in\Pi_{\di}(M(\A)^1)$ with $\cA^2_\pi(P)_{K_f,\sigma}\not=0$.
 Furthermore, by Proposition
\ref{p2.4} we have
\begin{equation}\label{5.22}
\sum_{\pi\in\Pi_{\di}(M(\A)^1)}\frac{\dim\cA^2_\pi(P)_{K_f,\sigma}}
{(1+\lambda_\pi^2)^{k/2}}<\infty
\end{equation}
for $k>m/2+1$, where $m=\dim M(\R)^1/K_{M,\infty}$. Combining (\ref{5.20})-
(\ref{5.22}), it follows that for each $k>m/2+16d^2+1$ there exists $C_k>0$
such that
\begin{equation*}
\begin{split}
\sum_{\pi\in\Pi_{\di}(M(\A)^1)}\int_{i\af_L^*/i\af_G^*}
\parallel\mN'_S(P,\pi,\lambda)\rho_\pi(P,\lambda,f)&\parallel_1|\nu^S_L(P,\pi,
\lambda)|\;d\lambda\\
&\le C_k \parallel D^{2k}f\parallel_{L^1(G(\A)^1)}.
\end{split}
\end{equation*}
This completes the proof.
\end{proof}

Now we return to the function $\phi_t$ defined by (\ref{5.3}). It follows 
from the definition that the restriction $\phi_t^1$ of $\phi_t$ belongs to 
$\Co^1(G(\A)^1,K_f,\sigma)$. We shall now modify $\phi_t$ in the following
way.
Let $\varphi\in C^\infty_0(\R)$ be such that $\varphi(u)=1$, if $|u|\le 1/2$,
and $\varphi(u)=0$, if $|u|\ge 1$. Let $d(x,y)$ denote the geodesic distance 
of $x,y\in X$ and set 
$$r(g_\infty):=d(g_\infty K_\infty,K_\infty).$$
 Given $t>0$, let $\varphi_t\in C^\infty_0(G(\R)^1)$ be defined by
$$\varphi_t(g_\infty)=\varphi(r^2(g_\infty)/t^{1/2}).$$
Then $\supp\varphi_t$ is contained in the set $\bigl\{g_\infty\in G(\R)^1\mid
r(g_\infty)<t^{1/4}\bigr\}$. Extend $\varphi_t$ to $G(\R)$ by
$$\varphi_t(g_\infty z)=\varphi_t(g_\infty),\quad g_\infty\in G(\R)^1,\;
z\in A_G(\R)^0,$$ 
and then to a function on $G(\A)$ by multiplying $\varphi_t$ by the
characteristic function of $K_f$. Put
\begin{equation}\label{5.22a}
\widetilde\phi_t(g)=\varphi_t(g)\phi_t(g),\quad g\in G(\A).
\end{equation}
Then the restriction $\widetilde\phi_t^1$ of $\widetilde\phi_t$ to $G(\A)^1$
 belongs to $C_c^\infty(G(\A)^1)$. 
\begin{prop}\label{p5.4}
There exist $C,c>0$ such that
$$|J_{\spec}(\phi_t^1)-J_{\spec}(\widetilde\phi_t^1)|\le C e^{-c/\sqrt{t}}$$
for $0<t\le1$. 
\end{prop}

\begin{proof}
Let $\psi_t=\phi_t-\widetilde\phi_t$ and $f_t=1-\varphi_t$. Let $\psi_t^1$ 
denote the restriction of $\psi_t$ to $G(\A)^1$.
 Then by Proposition \ref{p5.3} there exists $k\in\N$ such that
$$|J_{\spec}(\phi_t^1)-J_{\spec}(\widetilde\phi_t^1)|=|J_{\spec}(\psi_t^1)|
\le C_k\parallel (\Id+\Delta_G)^{k}\psi_t^1\parallel_{L^1(G(\A)^1)}.$$
In order to estimate the $L^1$-norm of $\psi_t^1$, recall that by definition
$$\psi_t(g_\infty g_f)=f_t(g_\infty)h_t^\sigma(g_\infty)\chi_{K_f}(g_f).$$
Hence
$$\parallel (\Id+\Delta_G)^k\psi_t^1\parallel_{L^1(G(\A)^1)}=\parallel
(\Id+\Delta_G)^k(f_th_t^\sigma)\parallel_{L^1(G(\R)^1)}.$$
Let $\gf(\R)^1$ be the Lie algebra of $G(\R)^1$ and let 
$X_1,...,X_a$ be an orthonormal basis of $\gf(\R)^1$. Then 
$\Delta_G=-\sum_i X_i^2$. Denote by $\nabla$ the canonical connection on
$G(\R)^1$. Then it  follows  that there exists $C>0$ such that
$$|(\Id+\Delta_G)^kf(g)|\le C\sum_{l=0}^{2k}\parallel\nabla^{l}f(g)\parallel,\quad g\in
G(\R)^1,$$
for all $f\in C^\infty(G(\R)^1)$.
By Proposition \ref{p1.1}  there exist constants $C,c>0$ such that
\begin{equation}\label{5.23}
\parallel\nabla^j h_t^\sigma(g)\parallel\le C t^{-(a+j)/2}e^{-cr^2(g)/t},\quad
g\in G(\R)^1,
\end{equation}
for  $j\le 2k$ and $0<t\le 1$. Let $\chi_t$ be the characteristic function of
 the set $\R-(-t^{1/4},t^{1/4})$. Recall that
 $f_t(g)=(1-\varphi)(r^2(g)/t^{1/2})$ and $(1-\varphi)(u)$ is constant for
$|u|\ge1$.
This implies that there exist constants $C,c>0$ such that
\begin{equation}\label{5.24}
\parallel\nabla^{j}f_t(g)\parallel\le Ct^{-k}\chi_t(r(g)),\quad 
g\in G(\R)^1,
\end{equation}
for $j\le 2k$ and $0<t\le 1$. 
Combining (\ref{5.23}) and (\ref{5.24}) we obtain
\begin{equation*}
\begin{split}
\sum_{l=0}^{2k}\parallel\nabla^{l}(f_th_t^\sigma)(g)\parallel&\le C_1 
t^{-a/2-2k}\chi_t(r(g))e^{-cr^2(g)/t}\\
&\le C_2 e^{-c_1/\sqrt{t}}e^{-c_1r^2(g)} 
\end{split}
\end{equation*}
for all $g\in G(\R)^1$ and $0<t\le 1$. Finally note that for every
$c>0$, $e^{-cr^2(g)}$ is an integrable function on $G(\R)^1$. This
finishes the proof. 
\end{proof}

\section[Proof of the main theorem]{Proof of the main theorem}
\setcounter{equation}{0}
In this section we 
 evaluate the geometric side of the trace formula at
the function $\widetilde\phi_t^1$ and investigate its asymptotic behaviour
as $t\to 0$. Then we compare the geometric and the spectral side and prove
our main theorem.
 
Let us briefly recall the structur of the geometric side $J_{\geo}$
of the trace formula \cite{A1}. The coarse $\ho$-expansion of $J_{\geo}(f)$
is a sum of distributions 
$$J_{\geo}(f)=\sum_{\ho\in \cO}J_{\ho}(f),\quad f\in C^\infty_c(G(\A)^1),$$
which are parametrized by the set $\cO$ of conjugacy classes of semisimple
elements in $G(\Q)$.  The distributions $J_{\ho}(f)$ are defined in 
\cite{A1}. We shall use the fine ${\ho}$-expansion of the spectral side 
\cite{A10} which expresses the distributions $J_{\ho}(f)$ in terms of weighted 
orbital integrals $J_M(\gamma,f)$. To describe the fine $\ho$-expansion
we have to introduce some notation. Suppose that $S$ is a finite set of
 valuations of $\Q$.  Set
$$G(\Q_S)^1=G(\Q_S)\cap G(\A)^1,$$
where
$$\Q_S=\prod_{v\in S}\Q_v.$$
Suppose that $\omega$ is a compact neighborhood of $1$ in $G(\A)^1.$
There is a finite set $S$ of valuations of $\Q$, which contains the 
Archimedean place, such that $\omega$ is the product of a compact neighborhood
of $1$ in $G(\Q_S)^1$ with $\prod_{v\notin S}K_v$.
Let $S^{0}_\omega$ be the minimal such set. Let $C^\infty_\omega(G(\A)^1)$ 
denote the space of functions in $C^\infty_c(G(\A)^1)$ which are supported on
$\omega.$ For any finite set $S\supset S^{0}_\omega$ set 
$$C^\infty_\omega(G(\Q_S)^1)=C^\infty_\omega(G(\A)^1)\cap C^\infty_c
(G(\Q_S)^1).$$
Let us recall the notion of $(M,S)$-equivalence \cite[p.205]{A10}. For any
$\gamma\in M(\Q)$ denote by $\gamma_s$ (resp. $\gamma_u$) the semisimple 
(resp. unipotent) Jordan component of $\gamma$. Then two elements $\gamma$
and $\gamma^\prime$ in $M(\Q)$ are called $(M,S)$-equivalent if there exists
$\delta\in M(\Q)$ with the following two properties.
\begin{enumerate}
\item[(i)] $\gamma_s$ is also the semisimple Jordan component of
$\delta^{-1}\gamma^\prime\delta$.
\item[(ii)] $\gamma_u$ and $(\delta^{-1}\gamma^\prime\delta)_u$, regarded
as unipotent elements in $M_{\gamma_s}(\Q_S)$, are 
$M_{\gamma_s}(\Q_S)$-conju\-gate.
\end{enumerate}
Denote by $(M(\Q))_{M,S}$ the set of $(M,S)$-equivalence classes in $M(\Q)$.
Note that $(M,S)$-equivalent elements $\gamma$ and $\gamma^\prime$ in $M(\Q)$
are, in particular, $M(\Q_S)$-conjugate. Given $\gamma\in M(\Q)$, let
$$J_M(\gamma,f),\quad f\in C^\infty_c(G(\Q_S)^1),$$
be the weighted orbital integral associated to $M$ and $\gamma$ \cite{A11}.
We observe that $J_M(\gamma,f)$ depends
only on the $M(\Q_S)$-orbit of $\gamma$. Then by Theorem 9.1 of \cite{A10} 
there exists a finite set $S_\omega \supset S^0_\omega$ of valuations 
of $\Q$ such that for all $S\supset S_\omega$ and any 
$f\in C^\infty_\omega(G(\Q_S)^1)$, we have
\begin{equation}\label{7.1}
J_{\geo}(f)=\sum_{M\in\cL}|W^M_0| |W^G_0|^{-1}\sum_
{\gamma\in(M(\Q))_{M,S}} a^M(S,\gamma)J_M(\gamma,f).
\end{equation}
This is the fine $\ho$-expansion of the geometric side of the trace formula.
The interior sum is finite.

Recall that the restriction $\widetilde\phi^1_t$ of $\widetilde\phi_t$ 
to $G(\A)^1$  belongs to $C^\infty_c(G(\A)^1)$
and hence, $J_{\geo}$ can be evaluated at $\widetilde\phi^1_t$. By
construction of  $\widetilde\phi^1_t$ there exists a compact neighborhood 
$\omega$ of 1 in $G(\A)^1$ and a finite set $S\supset S_\omega$ of valuations
of $\Q$ such that
$${\widetilde\phi}^1_t\in C^\infty_\omega(G(\Q_S)^1),\quad 0<t\le1.$$
Hence we can apply (\ref{7.1}) to evaluate 
$J_{\geo}({\widetilde\phi}^1_t)$. In this way our problem is reduced to the
investigation of the weighted orbital integrals 
$J_M(\gamma,{\widetilde\phi}^1_t)$. Actually for $\gamma\in M(\Q)$ we may
replace $\widetilde\phi^1_t$ by $\widetilde\phi_t$.

To begin with we establish some auxiliary results. Given $h\in G(\R)$, let
$$C_h=\{ g^{-1} hg\mid g\in G(\R)\}$$
be the conjugacy class of $h$ in $G(\R)$. 
\begin{lem}\label{l7.1}
Let $k\in K_\infty.$ Then $C_k\cap K_\infty$ is the $K_\infty$-conjugacy
class of $k.$
\end{lem}
\begin{proof} 
Let ${\gf}$ and $\kf$ denote the Lie algebras of  $G(\R)$ and $K_\infty$,
respectively. Let $\theta$ be a Cartan involution of $\gf$ with fixed point 
set $\kf$ and let ${\mathfrak p}$ be the $(-1)$-eigenspace of $\theta$.
Then the map
$$(k^\prime,X)\in K_\infty\times{\mathfrak p}\longmapsto k^\prime\exp(X)\in
G(\R)$$ 
is an analytic isomorphism of analytic manifolds. If $k_1\in K_\infty,$ then
$k_1$ is a $\theta$-invariant semisimple element. Therefore, its
centralizer $G_{k_1}$ is a reductive subgroup and the restriction of
$\theta$ to $G_{k_1}$ is a Cartan involution. Thus the restriction of the 
above Cartan decomposition to the centralizer of $k_1$ yields a Cartan
decomposition of $G_{k_1}(\R)$. 
Let $g\in G(\R)$ such that $g^{-1}kg\in K_\infty.$ Write $g=k^\prime\exp(X)$
with $k^\prime\in K_\infty$ and $X\in{\mathfrak p}.$ Since $g^{-1}kg$ is
$\theta$-invariant, we get
$$\exp(-X){k^\prime}^{-1} k k^\prime\exp(X)=\exp(X){k^\prime}^{-1}
k k^\prime \exp(-X).$$
Hence $\exp(2X)\in G_{ {k^\prime}^{-1} k k^\prime}(\R).$ 
From the Cartan decomposition of the latter group we conclude that 
$\exp(2X)=\exp(Y)$ for some 
$Y\in{\mathfrak p}_{ {k^\prime}^{-1} k {k^\prime}},$ and hence 
$X\in{\mathfrak p}_{ {k^\prime}^{-1} k {k^\prime}}.$
This implies that $g^{-1}kg=k^{\prime-1}kk'.$
\end{proof}
It follows from Lemma \ref{l7.1} that $C_k\cap K_\infty$ is a submanifold of 
$C_k.$

\begin{lem} \label{l7.2} Let $k\in K_\infty-\{\pm1\}.$ Then $C_k\cap K_\infty$
is a proper submani\-fold of $C_k.$
\end{lem}
\begin{proof} Let the notation be as in the previous lemma. First note that
the tangent space of $C_k$ at $k$ is given by
$$T_kC_k\cong ({\Ad}(k)-\Id)(\gf).$$
Furthermore
$${\Ad}(k)(\kf)\subset \kf,\quad {\Ad}(k)(\pg)\subset\pg.$$
Hence we get
$$T_k(C_k\cap K_\infty)=T_kC_k\cap{\kf}=({\Ad}(k)-\Id)(\kf),$$
and so the normal space $N_k$ to $C_k\cap K_\infty$ in $C_k$ at $k$ is given by
$$N_k\cong ({\Ad}(k)-\Id)(\pg).$$
Suppose that $\Ad(k)=\Id$ on $\pg$. Since $\kf=[\pg,\pg]$, it follows that
$\Ad(k)=\Id$ on $\gf$. Hence $k$ belongs to the center of $G^0$, which implies
that $k=\pm1$. Thus if $k\not=\pm1$, we have $\dim N_k>0$.
\end{proof}

Next we recall the notion of an induced space of orbits  \cite[p.255]{A11}.
Given an element $\gamma\in M(\Q_S)$, let $\gamma^G$ be the union of those
conjugacy classes in $G(\Q_S)$ which for any $P\in\cP(M)$ intersect 
$\gamma N_P(\Q_S)$ in an open set. There are only finitely many such 
conjugacy classes. 

\begin{prop}\label{p7.3} Let $d=\dim G(\R)^1/K_\infty$. Let $M\in\cL$ and 
$\gamma\in M(\Q)$. Then
$$\lim_{t\to0} t^{d/2} J_M(\gamma,\widetilde\phi_t)=0$$
if either $M\not=G$, or $M=G$ and $\gamma\not=\pm1$.
\end{prop}
\begin{proof}
By Corollary 6.2 of \cite{A11} the distribution $J_M(\gamma,\widetilde\phi_t)$ 
is given by the integral of $\widetilde\phi_t$ over $\gamma^G$ with respect 
to a measure $d\mu$ on $\gamma^G$ which is absolutely continuous with
respect to the invariant measure class. Thus $J_M(\gamma,\widetilde\phi_t)$
is equal to a finite sum of integrals of the form
$$\int_{G_{\gamma n}(\Q_S)\ba G(\Q_S)}\widetilde\phi_t(g^{-1}\gamma ng)
d\mu(g),$$
where $n\in N_P(\Q_S)$ for some $P\in\cP(M)$. Now recall that by (\ref{5.3})
and (\ref{5.22a}),
$\widetilde\phi_t(g)$ is the product of $\varphi_t(g_\infty)
h_t^\sigma(g_\infty)$ with
$\chi_{K_f}(g_f)$ for any $g=g_\infty g_f$. Hence our problem is reduced to
the investigation of the integral
$$\int_{G_{\gamma n_\infty}(\R)\ba G(\R)}(\varphi_th_t^\sigma)
(g^{-1}_\infty\gamma n_\infty g_\infty)\;d\mu(g_\infty).$$
Furthermore, by Proposition \ref{p1.1} there exists $C>0$ such that
$$|h_t^\sigma(g_\infty)|\le C t^{-d/2},\quad 0<t\le 1.$$
Hence it suffices to show that
\begin{equation}\label{7.2}
\lim_{t\to0}\int_{G_{\gamma n_\infty}(\R)\ba G(\R)}\varphi_t(g^{-1}_\infty
\gamma n_\infty g_\infty)\;d\mu(g_\infty)=0
\end{equation}
if either $M\not=G$, or $M=G$ and $\gamma\not=\pm1$. 

By definition of 
$\gamma^G$, the conjugacy class of $\gamma n$ in $G(\Q_S)$ has to intersect
$\gamma N_P(\Q_S)$ in an open subset. This implies that $\gamma n_\infty
\not=\pm1$, if either $M\not=G$, or $M=G$ and $\gamma\not=\pm1$. Then it
follows from Lemma \ref{l7.2} that $C_{\gamma n_\infty}\cap K_\infty$ is
a proper submanifold of $C_{\gamma n_\infty}$. The measure $d\mu(g_\infty)$
is of the form $f(g_\infty)dg_\infty$ for some smooth function $f$ on
$G(\R)$. Hence being a proper submanifold, $C_{\gamma n_\infty}\cap K_\infty$
is a subset of $C_{\gamma n_\infty}$ with measure zero with respect to 
$d\mu$. Next observe that
$$\int_{G_{\gamma n_\infty}(\R)\ba G(\R)}\varphi_t(g^{-1}_\infty
\gamma n_\infty g_\infty)\;|f(g_\infty)|\;dg_\infty<\infty.$$
Since $\supp\varphi_{t^\prime}\subset\supp\varphi_t$ for $t^\prime< t$, and
$0\le\varphi_t\le1$ for all $t>0$, it follows that there exists $C>0$ such that
$$\bigg|\int_{G_{\gamma n_\infty}(\R)\ba G(\R)}\varphi_t(g^{-1}_\infty
\gamma n_\infty g_\infty)\;d\mu(g_\infty)\bigg|\le C$$
for all $0<t\le 1$. Furthermore by definition of $\varphi_t$ we have
$$\lim_{t\to0}\varphi_t(x)=0$$
for all $x\in C_{\gamma n_\infty}-(C_{\gamma n_\infty}\cap K_\infty)$.
Since $C_{\gamma n_\infty}\cap K_\infty$ has measure zero with respect to
$d\mu$, (\ref{7.2}) follows by the dominated convergence theorem.
\end{proof}

We can now state the main result of this section.
\begin{theo}\label{th7.4} 
Let $d=\dim G(\R)^1/K_\infty$, let $K_f$ be an open compact subgroup of
$G(\A_f)$ and let $\sigma\in\Pi(\SO(n))$. Then
$$\lim_{t\to0}t^{d/2}J_{\geo}(\widetilde\phi_t^1)=\frac{\dim(\sigma)}
{(4\pi)^{d/2}}\vol(G(\Q)\ba G(\A)^1/K_f)
\bigl(1+{\bf 1}_{K_f}(-1)\bigr).$$
\end{theo}
\begin{proof}
By (\ref{7.1}) and Proposition \ref{p7.3} if follows that
$$\lim_{t\to0}t^{d/2} J_{\geo}(\widetilde\phi_t^1)=\lim_{t\to0}t^{d/2}
(a^G(S,1)\widetilde\phi_t^1(1)+a^G(S,-)\widetilde\phi_t^1(-1)).$$
By Theorem 8.2 of \cite{A10} we have
$$a^G(S,\pm1)=\vol(G(\Q)\ba G(\A)^1).$$
Furthermore
$$\widetilde\phi_t^1(\pm1)=h_t^\sigma(\pm1)\chi_{K_f}(\pm1).$$
Since $\sigma$ satisfies $\sigma(-1)=1$, it follows from (\ref{1.6}) that
$h_t^\sigma(-1)=h_t^\sigma(1)$
Finally, by Lemma \ref{l1.2} we have
$$h_t^\sigma(\pm1)=\frac{\dim(\sigma)}{(4\pi)^{d/2}}t^{-d/2}+O(t^{-(d-1)/2})$$
as $t\to0$. Combined with $\chi_{K_f}(\pm1)={\bf 1}_{K_f}(\pm1)\vol(K_f)^{-1}$,
the theorem follows.
\end{proof}

We shall now use the trace formula to prove the main results of this paper.
Recall that the coarse trace formula is the identity
$$J_{\spec}(f)=J_{\geo}(f),\quad f\in C^\infty_c(G(\A)^1),$$
between distributions on $G(\A)^1$ \cite{A1}. Applied to $\widetilde\phi_t^1$ we get
the equality
$$J_{\spec}(\widetilde\phi_t^1)=J_{\geo}(\widetilde\phi_t^1),\quad t>0.$$
Put $\varepsilon_{K_f}=1$, if $-1\in K_f$ and $\varepsilon_{K_f}=0$ 
otherwise.
Combining Theorem \ref{th5.2}, Proposition \ref{p5.4} and Theorem \ref{th7.4},
 we obtain
\begin{equation}\label{7.3}
\begin{split}
\sum_{\pi\in\Pi_{\di}(G(\A),\xi_0)}&m(\pi)
\dim\bigl(\H_{\pi_f}^{K_f}\bigr)\dim\bigl(\H_{\pi_\infty}(\sigma)\bigr)
e^{t\lambda_\pi}\\
&\sim\frac{\dim(\sigma)}{(4\pi)^{d/2}}
\vol(G(\Q)\ba G(\A)^1/K_f)(1+\varepsilon_{K_f})t^{-d/2}
\end{split}
\end{equation}
as $t\to0$. Applying Karamat's  theorem \cite[p.446]{Fe}, we obtain
\begin{equation}\label{7.4}
\begin{split}
\sum_{\pi\in\Pi_{\di}(G(\A),\xi_0)}&m(\pi)
\dim\bigl(\H_{\pi_f}^{K_f}\bigr)\dim\bigl(\H_{\pi_\infty}(\sigma)\bigr)
\\
&\sim\dim(\sigma)\frac{\vol(G(\Q)\ba G(\A)^1/K_f)}{(4\pi)^{d/2}\Gamma(d/2+1)}
(1+\varepsilon_{K_f})\lambda^{d/2}
\end{split}
\end{equation}
as $\lambda\to\infty$. By Lemma \ref{l2.3} it follows that this 
asymptotic formula continues
to hold if we replace the sum over $\Pi_{\di}(G(\A),\xi_0)$ by the sum over
$\Pi_{\cu}(G(\A),\xi_0)$. Finally note that by \cite{Sk} we have
 $m(\pi)=1$ for all $\pi\in\Pi_{\cu}(G(\A),\xi_0)$. This completes the proof of
Theorem \ref{th0.2}.

Now suppose that $K_f$ is the congruence subgroup $K(N)$ and $\Gamma(N)\subset
\SL_n(\Z)$ the principal congruence subgroup of level $N$. Then by (\ref{2.10})
we have
$$\vol(G(\Q)\ba G(\A)^1/K(N))=\varphi(N)\vol(\Gamma(N)\ba \SL_n(\R)).$$
Furthermore, $\varepsilon_{K(N)}=1$ if and only if $-1\in\Gamma(N)$. If
$-1$ is contained in $\Gamma(N)$, then the fibre of the canonical map
$$\Gamma(N)\ba \SL_n(\R)\to \Gamma(N)\ba \SL_n(\R)/\SO(n)$$
is equal to $\SO(n)/\{\pm1\}$. Otherwise the fibre is equal to $\SO(n)$. We
normalize the Haar measure on $\SL_n(\R)$ so that $\vol(\SO(n))=1$. Then in
either case we have
$$\vol(\Gamma(N)\ba \SL_n(\R))(1+\varepsilon_{K(N)})=\vol(\Gamma(N)\ba 
\SL_n(\R)/\SO(n)).$$
Let $X=\SL_n(\R)/\SO(n)$ and let $\lambda_0\le\lambda_1\le\cdots$ be the
eigenvalues, counted with multiplicity, of the Bochner-Laplace operator
$\Delta_\sigma$ acting in $L^2(\Gamma(N)\ba \SL_n(\R),\sigma)$.

Combining (\ref{5.14}), Proposition \ref{p5.4}, Theorem \ref{th7.4} and the 
above observations, we get
$$\sum_i e^{-t\lambda_i}=\dim(\sigma)\frac{\vol(\Gamma(N)\ba X)}{(4\pi)^{d/2}}
t^{-d/2}+o(t^{-d/2})$$
as $t\to0$. Using again Karamata's theorem \cite[p.446]{Fe}, we get
$$N_{\di}^{\Gamma(N)}(\lambda,\sigma)=\dim(\sigma)
\frac{\vol(\Gamma(N)\ba X)}{(4\pi)^{d/2}\Gamma(d/2+1)}
\lambda^{d/2}+o(\lambda^{d/2})$$
as $\lambda\to\infty$. By Proposition \ref{2.6} it follows that the same
asymptotic formula holds if we  replace
$N_{\di}^{\Gamma(N)}(\lambda,\sigma)$ by $N_{\cu}^{\Gamma(N)}(\lambda,\sigma)$.
This is exactly the statement of Theorem \ref{0.1}.
%%%%%%%%%%%%%%%%%%%%%%%%%%%%%%%%%%%%%%%%%%%%%%%%%%%%%%%%%%%%%%%%%%%%%%%%%%%%%%

\end{document}